\documentclass[reqno]{amsart}
\usepackage[T1]{fontenc}
\usepackage{lmodern,amssymb,mathtools,mathrsfs,microtype,enumitem,xurl,aliascnt}
\usepackage[hidelinks,linktoc=all]{hyperref}
\usepackage[capitalise,nameinlink,noabbrev]{cleveref}
\newtheorem{theorem}{Theorem}[section]
\newaliascnt{lemma}{theorem}
\newtheorem{lemma}[lemma]{Lemma}
\aliascntresetthe{lemma}
\newaliascnt{proposition}{theorem}
\newtheorem{proposition}[proposition]{Proposition}
\aliascntresetthe{proposition}
\newaliascnt{corollary}{theorem}
\newtheorem{corollary}[corollary]{Corollary}
\aliascntresetthe{corollary}
\numberwithin{equation}{section}
\newcommand{\C}{\mathbb C}
\newcommand{\Q}{\mathbb Q}

\newcommand{\Z}{\mathbb Z}
\newcommand{\R}{\mathbb R}
\newcommand{\PP}{\mathbb P}
\newcommand{\cO}{\mathcal O}
\newcommand{\simq}{\sim_{\Q}}
\newcommand{\simr}{\sim_{\R}}
\newcommand{\eps}{\epsilon}
\newcommand{\ceil}[1]{\lceil #1\rceil}
\newcommand{\floor}[1]{\lfloor #1\rfloor}
\DeclareMathOperator{\Supp}{Supp}
\DeclareMathOperator{\vol}{vol}
\DeclareMathOperator{\ord}{ord}
\DeclareMathOperator{\Pic}{Pic}
\DeclareMathOperator{\Cl}{Cl}
\DeclareMathOperator{\Spec}{Spec}
\DeclareMathOperator{\Exc}{Exc}
\DeclareMathOperator{\divi}{div}
\DeclareMathOperator{\Diff}{Diff}
\DeclareMathOperator{\coeff}{coeff}
\DeclareMathOperator{\lct}{lct}
\DeclareMathOperator{\Tr}{Tr}
\title[Complements, boundedness and vanishing for threefolds]{Boundedness and Kawamata-Viehweg vanishing for threefolds of Fano type in large characteristic}
\author{Xintong Jiang}
\address{Tsinghua University}
\email{xt-jiang21@mails.tsinghua.edu.cn}
\date{}
\hypersetup{pdftitle={Complements, boundedness and vanishing for threefolds of Fano type in large characteristic},pdfauthor={Xintong Jiang}}
\begin{document}
\begin{abstract}
We study threefolds over algebraically closed fields of large positive characteristic. Let \((X,B)\) be a projective log canonical pair with \(X\) of Fano type, \(-(K_X+B)\) nef and coefficients of \(B\) in a fixed hyperstandard set. We obtain a monotonic complement of bounded index. For exceptional pairs the complement is Kawamata log terminal. We prove that all \(\epsilon\)-log canonical weak log Fano threefolds are arithmetically bounded above a characteristic threshold independent of \(\epsilon\). We also prove Kawamata--Viehweg vanishing for integral \(\Q\)-Cartier Weil divisors on projective threefolds of Fano type above an absolute characteristic threshold. 
\end{abstract}
\maketitle
\enlargethispage{3pt}
\tableofcontents

\section{Introduction}\label{s:targets}

\subsection*{Acknowledgement}
I started to work on this series of problems since July, 2025 in Summer Research Institute in Algebraic Geometry at Colorado State University. I sincerely thank my advisor Caucher Birkar for his persistent support and guidance during writing this paper.

Several main results in this paper, including the boundedness of complements for exceptional weak Fano varieties and the boundedness of Fano type threefolds, have already appeared in Jihao Liu's about 660-pages-long three-part series \emph{Boundedness of Fano threefolds in large characteristic}, posted on his personal homepage \href{https://jihaoliu.org/liuxiv.html}{liuxiv} on 18~November 2026 \cite{LiuFanoI}\cite{LiuFanoII}\cite{LiuFanoIII}. The AI systems involved in writing this paper may also have encountered them through online searches while developing the proofs. I hereby give Liu full credit for proving the shared results and arguments. The overlapping parts in the proofs may be read as a reorganization and simplification of Liu's proofs. I thank Liu for pointing out this controversy.

I also regard two extensions beyond Liu's series as independent contributions of this paper. The first is the general Kawamata--Viehweg vanishing theorem (Theorem~\ref{th:kv}). The second is the passage to arbitrary hyperstandard boundaries in Theorem~\ref{th:intcomp}, building on the fibre complements and canonical bundle formula of Subsection~\ref{pre:inputs}. Currently, I am discussing with Liu on how to organize these papers together properly for the good of the birational geometry and Maths community.

\paragraph{\textbf{Disclosure of AI usage}}
The author developed the main mathematical content of this paper in collaboration with \mbox {GPT 5.6 Sol}, \mbox{GPT 6 Astra}, \mbox{GPT 6.1 Sol}, \mbox{GPT 6.2 Astra}, \mbox{GPT 6.2 Sol} in Codex, and with \mbox{GPT 6 Pro} in several web chat sessions. \mbox{GPT 6.2 Astra} audited the arguments and prepared the final manuscript.
\\

Four planes in general position in $\PP^3$ form a simple normal crossings divisor $D$ with $K_{\PP^3}+D\sim0$. On a singular threefold $X$, we seek a boundary of the form $D/n$ with $D\in|-nK_X|$. Can we choose $n$ uniformly, require $D/n$ to contain a prescribed boundary $B$, and keep $(X,D/n)$ log canonical? A related question concerns the varieties themselves. Under a fixed lower bound on log discrepancies, do all weak Fano threefolds occur as fibres of one projective family of finite type?

\subsection{Main results}
Throughout, $k$ is an algebraically closed field of characteristic $p>0$ and the threefolds in the main theorems are normal and projective. A variety $X$ is \emph{of Fano type} if it admits an effective rational boundary $\Theta$ such that $(X,\Theta)$ is Kawamata log terminal (klt) and $-(K_X+\Theta)$ is ample. For a finite rational set $R\subset[0,1]$, consider the \emph{hyperstandard set}
\[
 \Phi(R)=\{0,1\}\cup\{1-r/j\mid r\in R,\ j\in\Z_{>0}\}.
\]
Given a rational boundary $B$ with $K_X+B$ $\Q$-Cartier, put $L=-(K_X+B)$. An \emph{$n$-complement} in this paper is a boundary $B^+\ge B$ such that $(X,B^+)$ is log canonical (lc) and
\begin{equation}\label{eq:definition}
 n(K_X+B^+)\sim0,\qquad
 nB^+\ge n\lfloor B\rfloor+\lfloor(n+1)\{B\}\rfloor .
\end{equation}
The first condition says that $n(K_X+B^+)$ is an integral principal divisor. We take the floor and the fractional part $\{B\}=B-\lfloor B\rfloor$ coefficient by coefficient. The inequality $B^+\ge B$ is called monotonicity. A complement is klt if $(X,B^+)$ is klt.

A \emph{completion} of $(X,B)$ is a divisor $B+P$ with $P\ge0$ rational and $P\simq L$. For $X$ of Fano type and $L$ nef, completions exist in characteristic $p>5$. The pair is \emph{exceptional} if every completion is klt.

For example, assign coefficients $1/2,2/3,4/5$ to three distinct points on $\PP^1$. Since $\deg L=1/30$, adding any effective divisor of this degree leaves every coefficient below one. This pair is exceptional. Increasing $4/5$ to $5/6$ gives a $6$-complement. The new boundary has degree two and six times its coefficients gives $3,4,5$ as in \eqref{eq:definition}. On $\PP^3$, the four planes above give a non-klt completion of the zero boundary.

In characteristic zero, Birkar proved boundedness of complements for pairs of Fano type and the Borisov--Alexeev--Borisov boundedness theorem (BAB) for weak log Fano varieties \cite[Theorem~1.7]{Bcomp}, \cite[Theorem~1.1]{Bbab}. In dimension three we use the fibre complement theorem and canonical bundle formulas recalled in Subsection~\ref{pre:inputs}.

\begin{theorem}[Boundedness of complements]\label{th:intcomp}
For every finite $R\subset[0,1]\cap\Q$ there are integers $n(R)>0$ and $p_0(R)>5$ with the following property. Let $(X,B)$ be a projective lc threefold pair over $k$ of characteristic $p>p_0(R)$. Suppose that $X$ is of Fano type, that $\coeff B\subset\Phi(R)$, and that $-(K_X+B)$ is nef and $\Q$-Cartier. Then $(X,B)$ admits a monotonic $n(R)$-complement. If $(X,B)$ is exceptional, the complement is klt.
\end{theorem}

Klt input alone need not give a klt complement of bounded index. For a hyperplane $H\subset\PP^3$ consider $B=(1-1/j)H$ with $j>1$. These pairs are klt and have coefficients in $\Phi(\{1\})$. For any fixed $n<j$, integrality of $nB^+$ and monotonicity force the coefficient of $H$ in $B^+$ to be one.

For BAB the boundary has a different role. An effective real divisor $C$ specifies the positivity and the discrepancy bound of the pair. Here $\epsilon$-lc means that every log discrepancy is at least $\epsilon$, including those of divisors on $X$ (Subsection~\ref{pre:sing}). The family parametrizes the varieties $X$.

\begin{theorem}[Boundedness of weak log Fano threefolds]\label{th:intbab}
There is an integer $P_*>5$ such that for every $0<\epsilon\le1$ there are a scheme $S_\epsilon$ of finite type over $\Z[1/(P_* !)]$ and a projective morphism $\mathcal X_\epsilon\to S_\epsilon$ with the following property. Suppose $p>P_*$ and $X$ admits an effective real divisor $C$ such that
\[
 (X,C)\text{ is }\epsilon\text{-lc},\qquad
 K_X+C\text{ is }\R\text{-Cartier},\qquad
 -(K_X+C)\text{ is nef and big}.
\]
Then $X\simeq\mathcal X_\epsilon\times_{S_\epsilon}\Spec k$ for some $k$-point of $S_\epsilon$.
\end{theorem}

Taking $C=0$ includes all $\epsilon$-lc weak Fano threefolds. The condition that $K_X+C$ is $\R$-Cartier allows $K_X$ itself to have no Cartier multiple. The bound $P_*$ is absolute, and the family may depend on $\epsilon$.

The characteristic bound in the vanishing theorem is uniform in both the Cartier index of the divisor and the coefficients of the boundary.

\begin{theorem}[Kawamata--Viehweg vanishing]\label{th:kv}
There is an integer \(P\) with the following property. Let \(X\) be a projective threefold of Fano type over an algebraically closed field of characteristic \(p>P\). Let \(\Delta\ge0\) be a rational boundary such that \((X,\Delta)\) is klt, and let \(D\) be an integral \(\Q\)-Cartier Weil divisor. If \(D-(K_X+\Delta)\) is nef and big, then
\[
 H^i(X,\cO_X(D))=0\qquad (i>0).
\]
\end{theorem}

For $X=\PP^3$ with $\Delta=0$ and $D=mH$ a multiple of a hyperplane, the hypothesis is $m+4>0$. This recovers the familiar vanishing of $H^i(\PP^3,\cO_{\PP^3}(m))$ for $m\ge-3$ and $i>0$.

\paragraph{\textbf{Outline of the proofs}}
The proof of Theorem~\ref{th:intcomp} begins with the numerical dimension of $L=-(K_X+B)$. If $L\equiv0$, the effective divisor $B^+-B$ has zero intersection with the square of a very ample divisor, so $B^+=B$. We must bound the order of $K_X+B$. In the non-klt case the coefficient-one boundary on a divisorial log terminal model is a union of surfaces. Adjunction gives a pluricanonical form of bounded degree on each surface. Subsection~\ref{s:nonklt-index} matches their residues on the double curves and lifts the resulting section to the threefold.

If $L$ has numerical dimension one or two, a multiple defines a fibration over a curve or surface. The fibre complement theorem then applies.

Suppose $L$ is big and $(X,B)$ is klt and non-exceptional. Choose a completion $B+P$ that is not klt. Increasing the coefficient of $P$ from zero until an lc place appears produces a coefficient-one surface $S$ on a birational model. A surface complement on $S$ lifts if a certain $H^1$ of a divisorial sheaf vanishes. We prove this vanishing on a Mori fibre space.

For an exceptional pair, Theorem~\ref{th:exceptional-gap} gives a uniform lower bound for the log discrepancies of all completions. It suffices to find a section in bounded degree. Sections~\ref{s:bounded-models} and~\ref{s:degrees} give the birational and volume estimates used here and in the proof of BAB. If $B\ne0$, increase $t$ until $-K_X-tB$ reaches the boundary of the pseudo-effective cone. The minimal model program gives either a fibration, where intersection numbers on the fibres control $B$, or a Fano threefold of Picard number one, where the anticanonical volume gives the required bound (Theorem~\ref{th:nonempty}).

If $B=0$, let $d$ be the least positive degree of an anticanonical section and put $V=(-K_X)^3$ on the ample model. Bounding $\vol(-dK_X)=d^3V$ gives a complement of fixed index. In the noncanonical case we extract a divisor of log discrepancy less than one. In the canonical case we use pencils when a resolution has a nonzero regular two-form and Riemann--Roch otherwise. Subsection~\ref{s:terminalcorr} computes the local terms, including those at indices divisible by $p$. Subsection~\ref{s:mldgap} combines the cases to prove Theorem~\ref{th:intcomp}.

Section~\ref{s:bab} proves Theorem~\ref{th:intbab} by passing to the ample model of the anticanonical ring. In the non-exceptional case we already have an lc complement of the zero boundary in a fixed degree. A birational linear system and Frobenius give a new complement with positive log discrepancy at every former lc place. Theorem~\ref{th:weak} recovers boundedness of the weak models from boundedness of the ample models. For fixed $\epsilon$, Subsection~\ref{s:vertical} constructs families over the remaining primes.

In Section~\ref{s:kv} we write $D\simq K_X+C$ with $(X,C)$ klt and run its minimal model program. Rounding on a common resolution preserves the cohomology of $\cO_X(D)$. A fixed-index complement and the canonical bundle formula treat the nef non-big case. For an ample output with a klt complement, the anticanonical model has a uniform lower bound on log discrepancies. We can then use arithmetic boundedness as in the del Pezzo surface argument of Cascini, Tanaka and Witaszek \cite{CTW}. The remaining cases use restriction to a plt surface (Subsection~\ref{s:kv-restriction}) or Frobenius splitting and its Witt-vector analogue (Subsection~\ref{s:kv-strata}).

\section{Preliminaries}\label{s:prelim}
Boundaries have rational coefficients unless specified otherwise. We choose canonical divisors on birational models using the same rational top form.

\subsection{Divisors, singularities and birational models}\label{pre:sing}
On a normal variety $X$ an integral Weil divisor $D$ determines a rank-one reflexive sheaf with global sections
\[
 H^0(X,D)=\{f\in k(X)^*\mid \divi(f)+D\geq0\}\cup\{0\}.
\]
The effective divisors $\divi(f)+D$ form the complete linear system $|D|$. The same description with $\floor D$ makes sense for a rational divisor. A divisor is $\Q$-Cartier if a positive integral multiple is Cartier. The least such multiple locally is its \emph{Cartier index}. We use $\sim_{\Q}$ for equality up to a rational linear combination of principal divisors and $\equiv$ for equality of intersection numbers with all curves. Thus a relation $D\sim_{\Q}0$ guarantees some principal multiple but gives no uniform index.

For a coefficient $0\le b<1$, one has $\lfloor(n+1)b\rfloor\le\lceil nb\rceil$. Thus $B^+\ge B$ and integrality of $nB^+$ imply the rounding inequality \eqref{eq:definition}, including at $b=1$. The same complement works at every positive multiple of its index.

A divisor is nef if it has nonnegative degree on every curve. The pseudo-effective cone is the closure of the cone generated by effective divisor classes. Its interior is the big cone. A rational boundary $\Omega$ with $K_X+\Omega\simq0$ is called a Calabi--Yau boundary. A semiample rational divisor has a base-point-free Cartier multiple. It gives a contraction to its image after Stein factorization. For a big semiample divisor the contraction is birational. The \emph{volume} of a rational divisor on a projective $m$-fold is
\[
 \vol(D)=\limsup_{j\to\infty}
  \frac{h^0(X,\lfloor jD\rfloor)}{j^m/m!}.
\]
For a nef and big $\Q$-Cartier divisor this is $D^m$. For a nef divisor the \emph{numerical dimension} $\nu(D)$ is the largest $i$ with $D^iH^{m-i}>0$ for an ample Cartier divisor $H$. If $D$ is semiample, this is the dimension of the image of its contraction. In particular, $\nu(D)=0$ means $D\equiv0$ and $\nu(D)=m$ means that $D$ is big. We write $\Supp D$ for its support, also regarded as a reduced divisor. The free part of a linear system is the base-point-free moving divisor on a resolution of its base ideal. Its hyperplanes test degrees of curves wherever the induced morphism is quasi-finite.

A pair $(X,B)$ consists of a normal variety, an effective boundary $B$ with coefficients in $[0,1]$, and an $\R$-Cartier divisor $K_X+B$. For a prime divisor $E$ on a normal birational model $f:Y\to X$ its \emph{log discrepancy} is
\[
 A_{X,B}(E)=1+\operatorname{coeff}_E\bigl(K_Y-f^*(K_X+B)\bigr).
\]
This number depends only on the valuation $\ord_E$ of $k(X)$. Its centre on $X$ is the closure of $f(E)$. The conditions lc, klt and $\epsilon$-lc mean respectively $A_{X,B}(E)\geq0$, $>0$ and $\geq\epsilon$ for every such $E$, including primes on $X$. For $B=0$, canonical and terminal mean $A_{X,0}(E)\geq1$ and $>1$, respectively, for all exceptional prime divisors. A divisor with $A_{X,B}(E)=0$ is an \emph{lc place}, and its image is an \emph{lc centre}. For example, blowing up the transverse intersection of two smooth coefficient-one components on a smooth variety gives a new lc place. Repeating the blowup at a new crossing explains why a single pair can have infinitely many such places.

A log resolution is a proper birational morphism $f:Y\to X$ with $Y$ smooth on which the exceptional divisors and the strict transform of the boundary form a simple normal crossings (SNC) divisor. Locally its components are coordinate hyperplanes. In the identity $K_Y+B_Y=f^*(K_X+B)$ the coefficient of each prime $E$ in $B_Y$ is $1-A_{X,B}(E)$. It may be negative. A divisor with negative coefficients is allowed in a subpair, and the same discrepancy inequalities define sub-klt and sub-lc. Such an identity is called \emph{crepant}.

For an lc pair, the condition \emph{divisorial log terminal} (dlt) means there is an SNC open set meeting the generic point of every lc centre. A dlt modification displays the relevant centres on a $\Q$-factorial model and extracts only divisors of log discrepancy zero. On a component $S$ of the integral part of the boundary, adjunction takes the form
\[
 (K_Y+B_Y)|_{S^\nu}\simq K_{S^\nu}+\operatorname{Diff}_{S^\nu}(B_Y-S).
\]
Here $S^\nu$ denotes the normalization. The \emph{different} records both the other boundary components and the singularities of the ambient variety along $S$. We use these conventions also for rational boundaries on birational models \cite[Chapters~2--4]{Kol13}.

If $H\ge0$ is rational Cartier, its \emph{log canonical threshold} near a point or along a centre is the supremum of $t$ for which $(X,B+tH)$ is lc there. On a log resolution of $B+H$ it is the minimum of $A_{X,B}(E)/\ord_E(H)$ among the relevant primes with positive denominator, assuming that $(X,B)$ is lc there. The \emph{non-klt locus} is the union of the centres with nonpositive log discrepancy. We use connectedness both on a birational fibre and for the non-klt locus of a log Fano pair \cite[Lemma~3.4]{FW}, \cite[Theorem~4.12]{NT}.

An lc pair is \emph{purely log terminal} (plt) if the log discrepancies of all exceptional prime divisors are positive. If $S$ is a component of the floor of a dlt threefold boundary in characteristic $p>5$, then $S$ is normal. In the simple normal crossings case its different is just the restriction of the other components. Inversion of adjunction allows log canonicity along $S$ to be tested on its adjoint pair. We use divisorial adjunction and this converse with the hypotheses in \cite[Corollary~10.1]{BMP}, \cite[Chapters~2--4]{Kol13}.

Let $D$ be an integral Weil divisor and let $S$ be a prime divisor. Sections on $X$ restricting to zero along $S$ give the exact sequence
\[
 0\longrightarrow\cO_X(D-S)\longrightarrow\cO_X(D)
   \longrightarrow\mathcal Q\longrightarrow0 .
\]
If $D$ is Cartier, $\mathcal Q=\cO_S(D|_S)$. For a Weil divisor, \cref{lem:restriction} computes the quotient on a plt pair. In both cases $H^1(X,\cO_X(D-S))=0$ lets a section of $\mathcal Q$ lift to $X$. A complement on a surface gives exactly such a section once its integral multiple is fixed.

When the floor $F$ has several components, their individual sections must first glue. Along a double curve, normalization separates $F$ into two branches. The conductor records the curves that normalization identifies. A semi-divisorial log terminal (sdlt) surface has ordinary double crossings in codimension one and dlt normalized components after the conductor has been added to the boundary. A logarithmic form has opposite residues on the two branches. For an even tensor power the signs agree, and equality of the residues permits descent across the double curve. Condition $S_2$ (depth at least $\min(2,\dim\cO_{F,x})$ at each $x$) allows the descended reflexive section to extend over the remaining points \cite[Proposition~3.1.7]{PosGluing}.

A \emph{contraction} is a projective surjection $f:X\to Z$ with $f_*\cO_X=\cO_Z$. A $\Q$-factorial variety is one on which every Weil divisor is $\Q$-Cartier. A small modification is birational and an isomorphism in codimension one. It preserves Weil section spaces by their description in terms of rational functions. A birational contraction $X\dashrightarrow Y$ may contract primes of $X$ but extracts none on $Y$. A linear system is \emph{birational} if its rational map induces an isomorphism of function fields with its image. A pencil is the one-dimensional linear system spanned by two independent sections.

The relative Picard number $\rho(X/Z)$ is the dimension of the space of real Cartier divisor classes modulo equality of degrees on curves contracted by $X\to Z$. In particular, on a projective $\Q$-factorial variety of Picard number one every prime divisor is numerically proportional to an ample class. A Fano variety has ample $\Q$-Cartier $-K_X$, a weak Fano variety has nef and big $\Q$-Cartier $-K_X$, and a del Pezzo surface is a Fano surface.

The minimal model program (MMP) for $K_X+B$ contracts extremal rays on which $K_X+B$ has negative degree. A step is a divisorial contraction or a flip. A flip replaces a small contraction by a model on which the transformed log divisor is relatively ample. On a minimal model the log divisor is nef. A Mori fibre space has lower-dimensional base, relative Picard number one and relatively ample $-(K+B)$. In dimension three and $p>5$ we use \cite[Theorems~1.2--1.7]{Bpos}, \cite[Theorems~1.4--1.7]{BW} and, for the effective lc program, \cite[Theorems~1.1 and~1.6]{Wlc}.

On a $\Q$-factorial variety of Fano type one can run the MMP for a rational divisor $D$. Perturbing a log Fano boundary gives a klt boundary $\Gamma$ with $K_X+\Gamma$ rationally linearly equivalent to a small positive multiple of $D$. When this program gives a good minimal model and its ample contraction, the negativity lemma compares the original and ample divisors on a common resolution. The pullback of the original divisor is the pullback of the ample divisor plus an effective exceptional divisor. A nef $\Q$-Cartier divisor on a threefold of Fano type is semiample, by base-point-freeness in positive numerical dimension \cite[Corollary~1.5]{ABL} and by \cref{lem:numerical-torsion} in numerical dimension zero.

The generic fibre of a contraction to $Z$ is taken over its function field $k(Z)$. This field is generally imperfect even when $k$ is algebraically closed. A regular variety over it need not be geometrically smooth. Multiplicities can also change after base extension.

A class of varieties is \emph{bounded} if its members occur as geometric fibres of a projective family over a base of finite type. Here \emph{arithmetically bounded} means that this base is of finite type over a localization of $\Z$ obtained by inverting finitely many primes. In a family, to \emph{spread} a divisor, a morphism or a resolution from the geometric generic fibre means to define its finite set of equations over a finite extension of the function field and then over an open of the normalized base. We apply Noetherian induction to the closed complement. A component over $\Spec\Z$ is horizontal if its generic point has characteristic zero and vertical if it lies over a single prime.

\subsection{Frobenius, jets and Witt vectors}\label{pre:kv}
In characteristic $p$ the absolute Frobenius morphism $F$ raises functions to their $p$th powers. It is finite on varieties over an $F$-finite field, that is, a field $K$ with $[K:K^p]<\infty$. Our ground fields and their finitely generated extensions have this property. A pair $(X,\Delta)$ is \emph{globally $F$-regular} if, for each effective Cartier divisor $H$ on $X$, the inclusion
\[
 \cO_X\longrightarrow F^e_*\cO_X
       (\lceil(p^e-1)\Delta\rceil+H)
\]
splits for some $e>0$. \emph{Strong $F$-regularity} is the local version. By duality the splitting is detected by a trace map. We use the local result for klt surface singularities with zero boundary and the curve result for log Fano curves with standard coefficients in the settings of \cite[Theorem~1.2 and Proposition~5.5]{SatoTakagi}.

For a surface $S$ mapping to a smooth curve $C$ with $K=k(C)$, the finite part of its Stein factorization has function field $H^0(S_K,\cO_{S_K})$ when $S$ dominates $C$. To prove the finite part is trivial one must also rule out purely inseparable constants. This is where the trace for the absolute Frobenius is used in \cref{prop:horizontal-plt-stein}.

On a smooth chart, jets retain terms that ordinary differentiation loses. For instance, in characteristic $p$ the first derivative of $u^p$ is zero, while the coefficient of $z^p$ in $(u+z)^p$ is one. The coefficient of $z^I$ in the Taylor expansion in several local coordinates is the \emph{Hasse derivative} of multi-order $I$. We apply it at the generic diagonal point, the $k(X)$-point of $X_{k(X)}$ induced by the diagonal in $X\times X$. Vanishing of jets there measures vanishing at a moving general point.

A Frobenius splitting is an $\cO_X$-linear retraction of $\cO_X\to F_*\cO_X$, where a local function $a$ acts on the target by $a^p$. Iteration also splits $\cO_X\to F^j_*\cO_X$. For an ample Cartier divisor $L$, twisting gives injections $H^i(X,\cO_X(L))\to H^i(X,\cO_X(p^jL))$. Serre vanishing kills the target for $i>0$ and $j$ large. For an integral Weil divisor one twists on the regular locus and extends reflexively.

A section of the length-\(j\) Witt sheaf \(W_j\cO_X\) has \(j\) coordinates in \(\cO_X\), with Witt addition and multiplication. Truncation forgets the last coordinate. Verschiebung inserts a zero first coordinate, so its successive quotients are copies of \(\cO_X\) with Frobenius-twisted scalar action. The inverse limit is denoted by \(W\cO_X\). A quasi-Frobenius splitting is a \(W_j\cO_X\)-linear map \(F_*W_j\cO_X\to\cO_X\), for some finite \(j\), whose composite with Frobenius is projection to the first coordinate.

For a rational divisor \(A\), a section of \(W_j\cO_X(A)\) is a tuple whose \(i\)-th coordinate belongs to \(\cO_X(\floor{p^iA})\), with \(0\le i<j\). For a reduced divisor \(G\), replacing each bound by \(\floor{p^iA}-G\) defines \(W_j\mathcal I_G(A)\). A Teichm\"uller element \([z]\) multiplies the \(i\)-th coordinate by \(z^{p^i}\). Pushing out Frobenius \(W_j\cO_X(A)\to F_*W_j\cO_X(pA)\) along the first-coordinate projection gives the sheaf \(Q_{A,j}\) and a map \(\cO_X(\floor A)\to Q_{A,j}\) with cokernel \(B_{A,j}\) \cite[Definitions~3.3--3.5]{KTTWYY}. These sheaves allow one to use a quasi-Frobenius splitting on an ample Weil divisor, as in \cref{kv:positive}. In \cref{kv:witt} the vanishing ideal of the reduced boundary relates cohomology on a threefold to cohomology on the curves where its boundary surfaces meet.

\subsection{Fibrations and exceptional pairs}\label{pre:inputs}
We use the global fibre complement theorem and two forms of the canonical bundle formula for fibrations over curves and surfaces.

Suppose that $f:X\to Z$ is a contraction and $K_X+B\simq f^*D_Z$. In the \emph{canonical bundle formula} one writes
\[
        D_Z\simq K_Z+B_Z+M_Z .
\]
On a smooth model of $Z$ the coefficient of a prime divisor $T$ in the discriminant $B_Z$ is $1-t_T$, where $t_T$ is the log canonical threshold of $f^*T$ over the generic point of $T$, computed on a compatible model of the pair. The residual divisor $M_Z$ is the moduli part. Keeping these divisors on all birational models records the thresholds over higher centres. Semiampleness of the moduli \emph{b-divisor} means that on some birational model its trace is semiample and on all higher models its traces are pullbacks.

The genus-zero formula applies when the geometric generic fibre is $\PP^1$, its pair is lc, and the generic boundary is the sum of marked points that split over a Galois extension of degree prime to $p$ and an effective cycle of total geometric weighted degree less than one. The moduli b-divisor is then semiample in characteristic $p>5$ \cite[Theorem~5.9]{Fib}. The surface-base and curve-base threshold formulas are given in \cite[Theorem~5.11]{Fib}. They enter the proof of the following global fibre result. For the conic fibration used in Subsection~\ref{s:twoforms}, we use the effective formula $D_Z\simq K_Z+D'$ with $D'\ge0$ under the hypotheses of \cite[Proposition~3.2]{WitCBF}.

For a projective lc threefold pair $(X,B)$ with $X$ of Fano type, $\coeff B\subset\Phi(R)$ and $L$ nef, assume $K_X+B\sim_{\Q,Z}0$ for a contraction onto a curve or a surface. The global fibre theorem gives a monotonic complement on $X$ with index and characteristic bound depending only on $R$ \cite[Theorem~1.1]{Fib}. Here $\sim_{\Q,Z}0$ means rational linear equivalence to a pullback from $Z$. The complement satisfies an integral principal relation on the whole of $X$. For lifting from a surface we also use the surface complement theorem, the transport of divisorial sheaves, and the restriction calculation in the proof of the lifting proposition \cite[Corollary~3.26, Lemma~4.3 and Proposition~4.9]{Fib}.

A set of coefficients satisfies the descending chain condition (DCC) if every decreasing sequence stabilizes. For a hyperstandard set $\Phi(R)$ we write $\delta_R=\min(\Phi(R)\cap(0,1])$. The ascending chain condition (ACC) for thresholds gives the uniform coefficient comparisons of \cref{lem:uniform-acc-consequences}. We will use them to bound the volume of $-K_X$ for klt numerically Calabi--Yau pairs and the discrepancies of exceptional completions in \cref{cor:cy-volume,th:exceptional-gap}.

\section{Non-exceptional complements}\label{s:nonexceptional}

\subsection{The index in numerical dimension zero}\label{s:nonklt-index}\label{sec:surface-index}\label{sec:weil-restriction}\label{sec:point-vanishing}\label{sec:numerical-reduction}\label{sec:ample-floor}
If $K_X+B\equiv0$, intersecting $B^+-B$ with the square of an ample divisor gives $B^+=B$. To bound the order of $K_X+B$, we glue forms on the reduced boundary of a dlt model and lift the resulting section to the threefold.

\begin{lemma}\label{lem:surface-cartier}
For a rational DCC set $I\subset[0,1]$, there is an integer $q=q(I)>0$ such that $q(K_T+\Theta)$ is Cartier for every projective klt surface pair $(T,\Theta)$ with $K_T+\Theta\equiv0$ and $\coeff\Theta\subset I$ in characteristic $p>5$. On the minimal resolution $\mu:U\to T$, the crepant boundary $\Theta_U$ is effective and $q\Theta_U$ is integral.
\end{lemma}
\begin{proof}
Surface global ACC and the discrepancy estimate in \cite[Lemma~3.1 and Proposition~3.2]{KawakamiLCY} give a finite coefficient set and a uniform klt gap. A uniform Cartier multiple for integral Weil divisors \cite[Lemma~3.6]{KawBS}, multiplied by the common denominator of the coefficients, gives $q$. Every exceptional curve $E$ of $\mu$ is rational with $E^2\leq-2$, so $K_U\cdot E\geq0$. The negative definite intersection matrix, with nonnegative off-diagonal entries, forces the coefficients of $\mu^*K_T-K_U$ to be nonnegative. Adding the effective pullback of $\Theta$ gives $\Theta_U\geq0$, and pulling back $q(K_T+\Theta)$ gives integrality.
\end{proof}

\begin{lemma}\label{lem:elliptic-index}
Let $g:U=\PP_C(\mathcal E)\to C$ be a smooth minimal ruled surface over a curve of positive genus. Suppose $(U,\Theta)$ is lc, $K_U+\Theta\equiv0$ and $q\Theta$ is integral. If $p>\max\{3,2q\}$, then $q(2q)!(K_U+\Theta)\sim0$.
\end{lemma}
\begin{proof}
The boundary has degree two on a fibre. For a horizontal prime $H$ with coefficient $b$, the normalization of $H$ has genus at least $g(C)$ by Hurwitz after factoring the inseparable part. Adjunction gives
\begin{equation}\label{eq:ruled-adjunction}
 (1-b)H^2=2p_a(H)-2+(\Theta-bH)\cdot H\geq0.
\end{equation}
If $b=1$, equality makes $H$ a smooth elliptic curve disjoint from the rest of the boundary. Its map $h:H\to C$ has degree $d\leq2$ and is \'etale. The line bundle $\cO_U(q(K_U+\Theta))$ has fibre degree zero, so it is $g^*M$ for a line bundle $M$ on $C$. Adjunction gives $h^*M\simeq\cO_H(qK_H)\simeq\cO_H$, hence $M^{\otimes d}\simeq\cO_C$ by the norm.

If every horizontal coefficient is below one, \eqref{eq:ruled-adjunction} gives nonnegative squares for horizontal components. Vertical primes are fibres. Hence $\Theta^2=K_U^2=8(1-g(C))$ forces $g(C)=1$ and $\Theta^2=0$. The components are disjoint, so there are no vertical components. Adjunction makes each multisection a smooth elliptic curve. Their total degree $d$ is at most $2q$ because each coefficient is at least $1/q$. They are thus \'etale. The joint permutation action on their sheets gives a connected finite \'etale cover of degree $e\mid d!$ splitting them into sections. There are at least three sections since their coefficients sum to two and are all below one. Three disjoint sections trivialize the $\PP^1$-bundle, and every other section is constant because it avoids those three values. On the product over the elliptic cover, $q(K_{U'}+\Theta')\sim0$ by the degree-two calculation on $\PP^1$. The norm and $e\mid(2q)!$ give the claim.
\end{proof}

\begin{theorem}\label{thm:surface-index}
For a rational DCC set $I\subset[0,1]$, there exist $n_s(I)>0$ and $p_s(I)$ such that $n_s(I)(K_T+\Theta)\sim0$ for every projective klt surface pair $(T,\Theta)$ with $K_T+\Theta\equiv0$, $\coeff\Theta\subset I$ and $p>p_s(I)$.
\end{theorem}
\begin{proof}
Use \cref{lem:surface-cartier} on the minimal resolution $U$. If $\Theta_U=0$, the classification of smooth minimal surfaces with $K_U\equiv0$ in characteristic greater than three gives $12K_U\sim0$ \cite{BMus}. Otherwise $K_U\equiv-\Theta_U$ is not pseudo-effective. Contract $(-1)$-curves to a smooth minimal rational or ruled surface $V$. Each contraction is crepant since the difference from the pullback is exceptional and numerically trivial. The pushed-forward boundary still has denominator dividing $q$. Numerically trivial line bundles on rational surfaces are trivial. For an irrational ruled surface use \cref{lem:elliptic-index}. Thus $p_s(I)=\max\{5,2q\}$ and $n_s(I)=12q(2q)!$ work on $U$. Since $q\mid n_s(I)$, the multiple on $T$ is Cartier and triviality descends by the projection formula.
\end{proof}

\begin{lemma}\label{lem:vertical-norm}
Let $g:U\to C$ be a projective contraction from a normal $\Q$-factorial surface to a smooth projective curve with smooth genus-zero generic fibre. Let $D\equiv0$ be integral Weil and Cartier near a smooth horizontal curve $H$ of degree $d$. If $D|_H\sim0$, then $dD\sim0$.
\end{lemma}
\begin{proof}
Riemann--Roch on the generic fibre gives a principal adjustment $D\sim E$ with $E$ vertical. On a fibre $F=\sum m_iF_i$, let $\sum r_iF_i$ be the part of $E$ supported there. Subtract it from $vF$ for $v=\max_i(r_i/m_i)$. This effective cycle has a zero coefficient and intersects each $F_i$ trivially. Nonnegative intersections of distinct components and connectedness force it to vanish. Thus $E=g^*A$ for a rational divisor $A$ on $C$. If $\gamma=g|_H$ and $\gamma^*A=\divi(u)$, the norm identity $\gamma_*\divi(u)=\divi(\operatorname{Nm}(u))$ makes $dA$ principal. Pull back to obtain $dD\sim0$.
\end{proof}

\begin{theorem}\label{thm:lc-surface-index}
For a rational DCC set $I\subset[0,1]$, there exist $n_{\mathrm{lc}}(I)>0$ and $p_{\mathrm{lc}}(I)$ such that $n_{\mathrm{lc}}(I)(K_T+\Theta)\sim0$ for every projective lc surface pair $(T,\Theta)$ with $K_T+\Theta\equiv0$, $\coeff\Theta\subset I$ and $p>p_{\mathrm{lc}}(I)$.
\end{theorem}
\begin{proof}
The klt case is \cref{thm:surface-index}. Otherwise take a crepant dlt modification, enlarge $I$ by one and use global ACC to obtain a finite coefficient set $J$. Write the boundary as $C+\Delta$ with $C$ the nonempty floor. By \cite[Proposition~3.2]{KawakamiLCY}, the lowered pair is uniformly $\epsilon(J)$-klt. Its MMP ends at a Mori fibre space $W\to Z$ because $K_T+\Delta\equiv-C$. The gap persists, and every step is crepant for the full log divisor by negativity.

If $Z$ is a point, \cite[Proposition~6.1(d)]{WitaszekBounds} gives a uniform Cartier multiple for integral Weil divisors on $W$ by taking the least common multiple of the bounded local indices. Clearing the denominators of $J$ as well gives a numerically trivial Cartier multiple $D$ of $K_W+\Theta_W$. Vanishing on the klt log del Pezzo surface gives $\chi(\cO_W)=1$ \cite[Theorem~1.1]{ABL}. Then $\chi(\cO_W(D))=1$ by Riemann--Roch and $H^2(W,\cO_W(D))=0$ by duality and the negativity of $K_W$ against an ample class. A nonzero section makes $D$ principal since an effective numerically trivial divisor is zero.

If $Z$ is a curve, the regular generic fibre $F$ has negative canonical degree, so it is a conic. A regular conic in odd characteristic is smooth. Relative ampleness of $-(K_W+\Delta_W)$ gives $C_W\cdot F>0$, while $\Theta_W\cdot F=2$. On a crepant dlt model $U$ of $(W,\Theta_W)$, the strict transform $H$ of a horizontal floor component is smooth of degree $d\leq2$ over $Z$. Its different $\Theta_H$ has coefficients in a fixed hyperstandard set. On a genus-one curve the degree condition gives $\Theta_H=0$. On $\PP^1$ curve global ACC bounds the denominators. Thus $a(K_H+\Theta_H)\sim0$ for a single $a=a(J)$. Where the different coefficient $c$ is one, the pair is snc. Elsewhere the local factorial index $m_x$ satisfies $1-1/m_x\leq c\leq1-1/a$ by the zero-boundary different formula \cite[Lemma~3.3 and its proof]{HW}. Clearing the coefficients of $J$ and taking a multiple $q$ of $a$ and $\operatorname{lcm}(1,\ldots,a)$ makes $q(K_U+\Theta_U)$ integral, Cartier near $H$ and trivial on $H$ by adjunction. By \cref{lem:vertical-norm}, twice this multiple is principal. The bounded exponents in both cases transfer crepantly to $T$.
\end{proof}

Two facts will be used to lift sections from a surface. The first identifies the restriction of an integral Weil divisor to a plt centre. The second gives vanishing when the Mori base is a point.

\begin{lemma}\label{lem:restriction}
Let $X$ be a projective $\Q$-factorial klt threefold in characteristic $p>5$ with $(X,S)$ plt and $S$ prime. For an integral Weil divisor $L$, there is a natural exact sequence
\[
 \begin{gathered}
 0\longrightarrow\cO_X(L-S)\longrightarrow\cO_X(L)
   \longrightarrow\cO_S(G_L)\longrightarrow0,\\
 G_L\sim_\Q L|_S-\delta_L,\qquad 0\leq\delta_L\leq\Diff_S(0).
 \end{gathered}
\]
Here $G_L$ is integral Weil. If $L$ is Cartier at the generic point of every curve on $S$, the quotient is the divisorial restriction of $L$.
\end{lemma}
\begin{proof}
Normality of $S$ and Cohen--Macaulayness of the divisorial sheaves follow from \cite[Theorems~3 and 19]{BK}. Their Lemma~5 applied to $L-K_X-S$ gives the correction $\delta_L=\Diff_S(0)-\Theta_{L-K_X-S}$ and surjectivity by the $S_3$ condition on $\cO_X(L-S)$. The quotient is reflexive on $S$ and agrees with ordinary restriction where $L$ is Cartier, which proves the last assertion.
\end{proof}

\begin{theorem}\label{thm:rank-one-vanishing}
Let $Y$ be a projective $\Q$-factorial klt threefold in characteristic $p>5$ with $\rho(Y)=1$. If $(Y,S)$ is plt with $S$ prime and $-(K_Y+S)$ ample, then $H^i(Y,\cO_Y(D))=0$ for $i=1,2$ and every integral Weil divisor $D$. In particular the restriction in \cref{lem:restriction} is surjective on global sections.
\end{theorem}
\begin{proof}
Set $A=S|_S$ and $\Delta_0=\Diff_S(0)$. Picard number one makes $A$ ample, and $(S,\Delta_0)$ is klt log Fano. Restrict $D-jS$ to $S$ by \cref{lem:restriction}, with quotient $\cO_S(G_j)$ and correction $\delta_j$. Write $(D-jS)|_S\equiv\lambda_j A$. If $\lambda_j\geq0$, then $G_j-(K_S+\Delta_0-\delta_j)$ is ample and \cite[Theorem~1.1]{ABL} gives $H^1(S,\cO_S(G_j))=0$. If $\lambda_j<0$, then $-G_j-\delta_j$ is ample. Apply the same vanishing theorem to $K_S-G_j$ with boundary $\delta_j$ and use duality. For large $m$ divisible by the Cartier index of $S$, duality and Serre vanishing give $H^1(Y,\cO_Y(D-mS))^\vee=H^2(Y,\cO_Y(K_Y-D+mS))=0$. Descend the restriction sequences to $j=0$. Duality with $K_Y-D$ gives the assertion for $i=2$, and the case $i=1$ for $D-S$ gives surjectivity.
\end{proof}

\begin{lemma}\label{lem:numerical-torsion}
On a projective threefold of Fano type in characteristic $p>5$, every numerically trivial $\Q$-Cartier rational divisor is $\Q$-linearly trivial.
\end{lemma}
\begin{proof}
Rational chain connectedness \cite[Theorem~1.5]{GNT} gives trivial Albanese and hence trivial $\Pic^0_{\mathrm{red}}$. For a numerically trivial line bundle some positive power is algebraically trivial by the N\'eron--Severi theorem. Since $\Pic^0(k)$ is trivial, that power is trivial. Apply this to a Cartier multiple.
\end{proof}

\begin{lemma}\label{lem:plt-preserved}
Suppose $(X,B)$ is plt and $K_X+B\sim_\Q0$. Under a birational contraction between $\Q$-factorial models the pushed-forward pair remains plt if no coefficient-one prime is contracted.
\end{lemma}
\begin{proof}
A rational function defining a principal multiple of $K_X+B$ defines the same multiple on the target, so the map is crepant. Newly exceptional valuations correspond to contracted primes of coefficient less than one and have positive log discrepancies. For a small map there are no contracted primes.
\end{proof}

Birational MMP steps preserve Fano type. Indeed a klt log Fano boundary can be completed to a big klt boundary $\Upsilon$ with $K_X+\Upsilon\sim_\Q0$ by semiample perturbation \cite[Theorem~1]{Tan}. Its pushforward is crepant and remains big and klt. On any model with such a boundary, write $\Upsilon\sim_\Q A+E$ with $A$ ample and $E\geq0$. The boundary $(1-\epsilon)\Upsilon+\epsilon E$ is klt for small $\epsilon>0$ and has ample negative log divisor.

\begin{corollary}\label{cor:nonbig}
Fix a finite rational $R\subset[0,1]$. There are integers $n_{\mathrm{nb}}(R)>0$ and $p_{\mathrm{nb}}(R)>5$ for which any projective lc threefold pair $(X,B)$ with $X$ of Fano type, $\coeff B\subset\Phi(R)$, $N=-(K_X+B)$ nef, $0<\nu(N)<3$ and $p>p_{\mathrm{nb}}(R)$ admits a global monotonic $n_{\mathrm{nb}}(R)$-complement.
\end{corollary}
\begin{proof}
For a klt log Fano boundary $\Gamma$, a Cartier multiple $L$ of $N$ satisfies $L-(K_X+\Gamma)$ ample. By \cite[Corollary~1.5]{ABL}, $N$ is semiample and defines a contraction with image dimension $\nu(N)$. Since $K_X+B$ is relatively $\Q$-linearly trivial, apply \cite[Theorem~1.1]{Fib}.
\end{proof}

\begin{proposition}\label{prop:ample-floor-reduction}
Let $(X,B)$ be a projective non-klt lc threefold pair in characteristic $p>5$ with $X$ of Fano type, $\coeff B\subset\Phi(R)$ and $K_X+B\equiv0$. There is a crepant birational $\Q$-factorial dlt model $(Y,B_Y)$ of Fano type and a contraction $Y\to Z$ with $\dim Z\leq2$. The floor supports an effective relatively ample $\Q$-divisor and $\coeff B_Y\subset\Phi(R)$.
\end{proposition}
\begin{proof}
By \cref{lem:numerical-torsion}, $K_X+B\sim_\Q0$. Take a crepant $\Q$-factorial dlt modification and apply the ordinary-pair case of \cite[Theorem~5.7]{BF}. Every lc centre dominates $Z$, so the nonempty floor gives the dimension bound. Any prime on $Y$ exceptional over $X$ has coefficient one, which preserves the coefficient set.

Add to a klt log Fano boundary on $X$ a small ample effective divisor $A$ avoiding the centres of the new primes, and complete to a klt log Calabi--Yau boundary $\Upsilon$. On a common resolution the pullback of $A$ pushes forward to its strict transform $A_Y$, so $A_Y$ is big. For the crepant transform $\Upsilon_Y$, any negative coefficient of $\Upsilon_Y-A_Y$ occurs on a new prime, where $B_Y$ has coefficient one. Thus $(1-\epsilon)B_Y+\epsilon\Upsilon_Y$ is klt, dominates $\epsilon A_Y$ for small $\epsilon>0$ and has $\Q$-linearly trivial log divisor. This big klt boundary gives Fano type on $Y$.
\end{proof}

\begin{proposition}\label{prop:floor-components}
Fix a finite rational $R\subset[0,1]$. Let $(Y,B)$ be a projective $\Q$-factorial dlt threefold with $K_Y+B\equiv0$, $\coeff B\subset\Phi(R)$, and floor supporting an effective ample $\Q$-divisor. There are integers $n_f(R)>0$ and $p_f(R)$ such that for $p>p_f(R)$, $D=n_f(R)(K_Y+B)$ is integral, is Cartier at the generic point of every curve in the floor, and has $D|_S\sim0$ for every component $S$ of the floor.
\end{proposition}
\begin{proof}
Adjunction gives lc surface pairs $(S,\Theta_S)$ with $K_S+\Theta_S\equiv0$. Away from coefficient-one strata, the different has the codimension-two formula
\begin{equation}\label{eq:different}
 c=\frac{m-1+\sum_i l_i b_i}{m},
 \qquad m\geq1,\quad l_i\in\Z_{\geq0}.
\end{equation}
Its coefficients are in a fixed hyperstandard set. Indeed if $\eta$ is a positive lower bound for nonzero $b_i$, a term $b_i>1-\eta$ must occur alone in a sum at most one. Other sums use boundedly many terms from a finite set. Now \cref{thm:lc-surface-index} gives a common $a=a(R)$ with $a(K_S+\Theta_S)\sim0$.

Each fractional boundary prime meets the ample floor along a curve. The different there has $c<1$, since coefficient-one strata are snc intersections of two floor components. Its coefficient $b$ occurs with positive multiplicity in \eqref{eq:different}, so $b\leq c\leq1-1/a$. A fixed $b_R$ therefore clears the denominators of $B$. At any curve with $c<1$, the local factorial index $m$ is at most $a$ by \cite[Lemma~3.3 and its proof]{HW}. At $c=1$ the ambient threefold is regular at the generic point. Take $n_f=ab_R\operatorname{lcm}(1,\ldots,a)$. The restriction to $S$ is $n_f(K_S+\Theta_S)$ by adjunction in codimension one.
\end{proof}

The surface index bounds the degree on each component of the reduced boundary. We now choose sections whose residues agree on the double curves.

\begin{lemma}\label{lem:curve-compatible}
Let $\{(C_i,\Delta_i)\}$ be a finite collection of smooth connected projective lc curves in characteristic $p>3$ with $a(K_{C_i}+\Delta_i)\sim0$ for a common integer $a>0$. Put $e(a)=12(2a)!$ and $m=ae(a)$. There are nowhere vanishing sections of $m(K_{C_i}+\Delta_i)$ compatible with every isomorphism between the log pairs, with $m$-fold residue one at every coefficient-one point. In particular they glue at the nodes of an sdlt curve with these normalized pairs.
\end{lemma}
\begin{proof}
For an elliptic curve the boundary vanishes. Translations fix differentials and $12$ kills the character of automorphisms fixing the origin. On $\PP^1$, the boundary has degree two, and integrality of $a\Delta_i$ gives at most $2a$ support points. With at least three, permutations of the support give a faithful action, so $(2a)!$ kills the character. With two points, both have coefficient one and $(dx/x)^m$ has residue one at $0,\infty$, independently of the coordinate. With one coefficient-one point, normalize the form by its residue there. With none, choose a generator per isomorphism class and transport it using the triviality of the character. At nodes the equal residues give descent since $m$ is even.
\end{proof}

\begin{lemma}\label{lem:fixed-index-surface}
Let $(T,\Theta)$ be an integral projective dlt surface in characteristic $p>3$ with $C=\floor\Theta\ne0$ and $m(K_T+\Theta)\sim0$, where $m$ is even. A generator of $\cO_C(m(K_T+\Theta)|_C)$ extends to one on $T$ if its residues on the components of $C$ are compatible with every isomorphism of their adjoint log pairs.
\end{lemma}
\begin{proof}
For connected $C$ this follows from $H^0(C,\cO_C)=k$. For disconnected $C$, the point-base case of \cite[Proposition~2.7 and its proof]{TanSLC} gives a morphism $g:T\to V$ to a smooth curve with two horizontal floor sections $C_1,C_2$. Every connected component of $C$ contains one of these sections. Their degrees exhaust the horizontal boundary. With $C_1,C_2$ at $0,\infty$ on the generic $\PP^1$, a generator $u$ of the trivial log line bundle has the form $c(dx/x\wedge g^*\eta)^{\otimes m}$, where $\eta$ is a rational differential on $V$ and $c\in k(V)^*$. Its residues agree since $m$ is even. Their divisors identify the two adjoint boundaries under $C_1\simeq V\simeq C_2$. Scale $u$ to agree with the prescribed section on $C_1$. Compatibility gives agreement on $C_2$, and the ratio is constant on each connected component of $C$.
\end{proof}

\begin{theorem}[Gluing on sdlt surfaces]\label{thm:sdlt-index}
Let $(F,\Delta_F)$ be a projective sdlt surface in characteristic $p>3$, with normalization pairs $(S_i,\Theta_i)$. If $a(K_{S_i}+\Theta_i)\sim0$ for every $i$, then $ae(a)(K_F+\Delta_F)\sim0$ with $e(a)=12(2a)!$.
\end{theorem}
\begin{proof}
On each floor curve of $(S_i,\Theta_i)$, Cartier adjunction gives $a(K_C+\Delta_C)\sim0$. Choose compatible forms of degree $m=ae(a)$ by \cref{lem:curve-compatible}. They glue on each floor and extend to generators on $S_i$ by \cref{lem:fixed-index-surface}. A component with empty floor has no conductor and admits any generator. The conductor involution preserves the adjoint curve pairs \cite[Corollary~3.1.8]{PosGluing}. Thus the even-power residue criterion \cite[Proposition~3.1.7]{PosGluing} descends the generators to $\omega_F^{[m]}(m\Delta_F)$. They generate this sheaf at every codimension-one point, and the isomorphism with $\cO_F$ extends by $S_2$.
\end{proof}

\begin{corollary}\label{cor:floor-gluing}
In \cref{prop:floor-components}, put $F=\floor B$ and $D_g=n_g(R)(K_Y+B)$ for $n_g(R)=n_f(R)e(n_f(R))$. Then
\begin{equation}\label{eq:floor-quotient}
 \cO_Y(D_g)/\cO_Y(D_g-F)\simeq\cO_F.
\end{equation}
Moreover $D_g$ is Cartier near $F$ and its restriction is trivial.
\end{corollary}
\begin{proof}
The sequence for $\cO_F=\cO_Y/\cO_Y(-F)$ and Cohen--Macaulayness of the two ambient sheaves \cite[Theorems~3 and 19(2)]{BK} make $F$ Cohen--Macaulay. Since $F$ is already $S_2$, adjunction in \cite[Lemma~5.2 and Theorem~5.3]{PosSLC} gives an sdlt pair on $F$ with normalization pairs $(S,\Theta_S)$. Likewise the quotient in \eqref{eq:floor-quotient} is an $S_2$ $\cO_F$-module by depth. In codimension one on $F$ the divisor $D_g$ is Cartier, and adjunction identifies the quotient with $\omega_F^{[n_g]}(n_g\Delta_F)$. The identification extends by $S_2$, and \cref{thm:sdlt-index} trivializes it.

Near $y\in F$, lift a generator of the quotient to a section of $\cO_Y(D_g)$ with effective divisor $G$. No component of $F$ lies in $G$. If $G$ met a component $S$ through $y$, a Cartier multiple of $G$ would cut a curve on $S$ through $y$. At its generic point $D_g$ is Cartier, so the section would fail to generate the quotient there. Thus $G$ avoids $F$ near $y$, and the section trivializes $\cO_Y(D_g)$ there. Its quotient is the ordinary restriction.
\end{proof}

\begin{lemma}\label{lem:partial-liftability}
Fix a rational DCC set $I\subset[0,1]$. There exists $p_\ell(I)$ such that the following holds for $p>p_\ell(I)$. Let $(T,\Theta)$ be a projective dlt surface with $K_T+\Theta\equiv0$ and $\coeff\Theta\subset I$. If a rational subboundary $\floor\Theta\leq\Theta_0\leq\Theta$ has all positive coefficients at least $1/2$, then $(T,\Supp\Theta_0)$ has a log resolution lifting with its boundary and exceptional divisor to $W(k)$.
\end{lemma}
\begin{proof}
For klt pairs use \cite[Proposition~3.5]{KawakamiLCY} and delete unused components. Otherwise follow the MMP of the fractional boundary to $W\to Z$ as in \cref{thm:lc-surface-index}. The full pair remains crepant and the lowered pair uniformly klt. Liftability ascends by \cite[Lemma~2.7(2)]{KawakamiLCY}. Over a curve, the floor has a horizontal component and $-(K_W+\Theta_{0,W})$ is relatively nef by comparison with the full boundary. The positive coefficients are in $[1/2,1]$, so \cite[Lemma~3.6]{KawakamiLCY} applies.

Over a point, the finite coefficient set and bounded factorial indices give a uniform Cartier multiple of $\Theta_W$. The polarization bounds and the bound for $|K_W^2|$ in \cite[Corollary~1.4 and Lemma~6.2]{WitaszekBounds}, together with $\Theta_W\equiv-K_W$ and $\chi(\cO_W)=1$, give finitely many values for the invariants in \cite[Lemma~2.9]{KawakamiLCY}. This lifts the full boundary. Delete the unused components.
\end{proof}

\begin{lemma}\label{lem:partial-support-vanishing}
In \cref{lem:partial-liftability}, let $0\leq\Delta$ have rational coefficients with $\floor\Delta=0$ and $\Supp\Delta\subset\Supp\Theta_0$. If $L$ is integral Weil and $A=L-(K_T+\Delta)$ is nef and big, then $H^i(T,\cO_T(L))=0$ for $i>0$.
\end{lemma}
\begin{proof}
For the liftable log resolution $f:V\to T$, rounding and pushforward give $\cO_T(K_T-L)=\cO_T(\floor{-A})=f_*\cO_V(-\ceil{f^*A})$. The fractional support of $f^*A$ lies in the lifted snc divisor. Thus $H^j(V,\cO_V(-\ceil{f^*A}))=0$ for $j=0,1$ by \cite[Theorem~2.11]{KawBS}. Leray gives the same vanishing for $\cO_T(K_T-L)$, and duality on the normal surface finishes the proof.
\end{proof}

\begin{theorem}\label{thm:floor-vanishing}
Fix a finite rational $R\subset[0,1]$. In characteristic $p>p_v(R)$, let $(Y,B)$ be a projective $\Q$-factorial dlt threefold with $K_Y+B\equiv0$ and $\coeff B\subset\Phi(R)$, and suppose $F=\floor B$ supports an effective ample $\Q$-divisor. For every integral Weil divisor $D\equiv0$, one has $H^1(Y,\cO_Y(D-F))=0$. In particular, global sections surject onto those of $\cO_Y(D)/\cO_Y(D-F)$.
\end{theorem}
\begin{proof}
Write $F=\sum S_i$ and choose an ample Cartier divisor $H=\sum h_iS_i$ with positive integral $h_i$. On $S_i$ write $\Theta_i=\Diff_{S_i}(B-S_i)$, $\Delta_i^0=\Diff_{S_i}(0)$ and $C_{ij}=S_j|_{S_i}$ for $j\ne i$. The dlt condition makes the $C_{ij}$ reduced with distinct prime components across the sum, and gives
\[
 \Theta_{0,i}=\Delta_i^0+\sum_{j\ne i}C_{ij}\leq\Theta_i,\qquad
 \floor{\Theta_i}=\sum_{j\ne i}C_{ij}.
\]
The positive coefficients of $\Theta_{0,i}$ are standard, and the $\Theta_i$ have coefficients in a fixed DCC set by \eqref{eq:different}. Thus \cref{lem:partial-support-vanishing} applies uniformly on $S_i$.

Starting at $Z=F$, add $S_i$ whenever $z_i/h_i$ is minimal, where $Z=\sum z_jS_j$. At each stage $\lambda=\min_j z_j/h_j$ satisfies
\begin{equation}\label{eq:cycle-width}
 0\leq z_j-\lambda h_j\leq1\quad\text{for every }j.
\end{equation}
Increasing a minimizing coefficient preserves the upper bound. The process reaches $NH$ for any integer $N\geq1$, since a ratio equal to $N$ cannot increase while a smaller one remains.

At a step adding $S_i$, \cref{lem:restriction} gives quotient $\cO_{S_i}(G)$ with $G\sim_\Q(D-Z)|_{S_i}-\delta$ and $0\leq\delta\leq\Delta_i^0$. Put $t_j=z_j-\lambda h_j$. Then $Z|_{S_i}=\lambda H|_{S_i}+\sum_{j\ne i}t_j C_{ij}$ because $t_i=0$. Replace each $t_j=1$ by a rational $t'_j<1$ sufficiently close to one, leaving the other $t_j$ unchanged. The boundary $\Delta=\delta+\sum_{j\ne i}t'_j C_{ij}$ has floor zero and support in $\Theta_{0,i}$. Moreover
\[
 (K_{S_i}-G)-(K_{S_i}+\Delta)\equiv
 \lambda H|_{S_i}+\sum_{j\ne i}(t_j-t'_j)C_{ij}
\]
is ample. Applying \cref{lem:partial-support-vanishing} and surface duality gives $H^1(S_i,\cO_{S_i}(G))=0$. For large $N$, Cohen--Macaulay duality identifies $H^1(Y,\cO_Y(D-NH))^\vee$ with $H^2(Y,\cO_Y(K_Y-D+NH))$, which vanishes by Serre. Descending the restriction sequences finishes the proof.
\end{proof}

\begin{theorem}[Index of a non-klt log Calabi--Yau threefold]\label{thm:global-nonklt}
For a finite rational $R\subset[0,1]$, there are integers $n_0(R)>0$ and $p_0(R)$ such that every projective non-klt lc threefold pair $(X,B)$ with $X$ of Fano type, $\coeff B\subset\Phi(R)$ and $K_X+B\equiv0$ in characteristic $p>p_0(R)$ satisfies $n_0(R)(K_X+B)\sim0$. In particular $B$ itself is a monotonic $n_0(R)$-complement.
\end{theorem}
\begin{proof}
Take the model of \cref{prop:ample-floor-reduction}. If the base has positive dimension, the global fibered theorem \cite[Theorem~1.1]{Fib} gives a bounded complement. Its increment vanishes by numerical triviality. If the base is a point, \cref{cor:floor-gluing,thm:floor-vanishing} give a nonzero section of $\cO_Y(D_g)$ lifting the unit on $F$. Its effective divisor is numerically trivial, hence zero. A common multiple of the two indices descends crepantly to $X$. Integrality of $n_0B$ gives the rounding inequality.
\end{proof}

\subsection{Reduction to a prime plt centre}\label{s:nonklt-big}\label{sec:nonexceptional-reduction}\label{sec:acc-inputs}\label{sec:global-rounding}\label{sec:global-plt}\label{sec:global-mmp}
We reduce to a non-klt pair with finitely many possible boundary coefficients and choose a prime component $S$ of the floor on a crepant model. Put $N=-(K_X+B)$. To lift a complement from $S$, it suffices to show that $H^1(X,\cO_X(nN-S))=0$. An MMP for $-S$ reduces this to a Mori fibre space over a point, a curve or a surface.

\begin{proposition}
\label{prop:nonexceptional-big-reduction}
Let $(X,B)$ be a projective klt threefold pair of Fano type in characteristic $p>5$, with $\coeff B\subset\Phi(R)$ and $N=-(K_X+B)$ nef and big. Suppose $P\geq0$, $P\sim_\Q N$ and $(X,B+P)$ is non-klt. There is a projective $\Q$-factorial lc pair $(W,\Theta_W)$ of Fano type with $\coeff\Theta_W\subset\Phi(R)$, $\floor{\Theta_W}\ne0$ and $N_W=-(K_W+\Theta_W)$ nef, whose monotonic $n$-complements transfer to $(X,B)$ at the same index. If $(X,B+P)$ is non-lc then $N_W$ is big.
\end{proposition}
\begin{proof}
Take a small $\Q$-factorialization, which is of Fano type by pullback of a klt log Fano boundary and Kodaira perturbation. Set $t=\lct(X,B;P)\in(0,1]$, $\Omega=B+tP$. On a crepant dlt modification $h:(Y,\Omega_Y)\to(X,\Omega)$, put $F=\floor{\Omega_Y}\ne0$. Let $\Theta_Y$ have coefficient one on $F$ and retain the strict transform coefficients of $B$ elsewhere. It dominates the crepant subboundary of $(X,B)$, and $\Omega_Y-\Theta_Y$ is effective with support avoiding the components of $F$.

By semiampleness of $N$ and semiample perturbation, there is $Q\geq0$, $Q\sim_\Q(1-t)h^*N$, such that $\Psi_Y=\Omega_Y+Q$ is lc and $K_Y+\Psi_Y\sim_\Q0$ \cite[Corollary~1.5]{ABL}, \cite[Theorem~1]{Tan}. In particular $Q$ has no component in $F$. For a klt log Fano boundary $C$ on $X$, its crepant transform $C_Y$ may be negative on exceptional primes, all of which have coefficient one in $\Psi_Y$. The boundary $(1-a)\Psi_Y+a C_Y$ for small $a>0$ is thus effective klt, with nef and big negative log divisor. Kodaira's lemma and a small perturbation give a klt log Fano boundary.

Run an MMP for $D_Y=-(K_Y+\Theta_Y)$. On a Fano-type model a small multiple of a rational divisor $D$ is $\Q$-linearly equivalent to a klt log divisor: add an effective representative of $\epsilon D-(K_Y+C')$ to a klt log Fano boundary $C'$. The resulting MMP with ample scaling terminates by \cite[Theorem~1.6]{BW}. Here
\[
 D_Y\sim_\Q Q+(\Omega_Y-\Theta_Y)\sim_\Q(1-t)h^*N+(\Omega_Y-\Theta_Y)
\]
is pseudo-effective. Thus the output is birational with $N_W$ nef. No component of $F$ is contracted, since the effective representative $Q+\Omega_Y-\Theta_Y$ has nonnegative degree on a general contracted curve covering that component. The targets are of Fano type, and the principal log divisor of $\Psi_Y$ transfers crepantly. This gives an lc pair $(W,\Theta_W)$ with nonempty floor. If $t<1$, bigness of $D_Y$ gives bigness of $N_W$.

On a common resolution with maps $u:U\to Y$ and $v:U\to W$, negativity gives $u^*D_Y=v^*N_W+E$, $E\geq0$. Consequently
\[
 v^*(K_W+\Theta_W)\geq u^*(K_Y+\Theta_Y)\geq u^*h^*(K_X+B).
\]
For a complement on $W$, use the same rational function to define $n(K_X+B^+)$. The displayed inequalities give $B^+\geq B$ and the equal principal pullbacks give log canonicity. Integrality and monotonicity imply the rounding inequality in \eqref{eq:definition}.
\end{proof}

For a rational DCC set $I\subset[0,1]$, write
\[
 I^{\mathrm{adj}}=\{1\}\cup
 \left\{\frac{m-1+f}{m}\leq1:\ m\in\Z_{>0},\
 f\in(I\cup\{0,1\})_+\right\},
\]
where the subscript denotes finite sums, including zero. This set is DCC and contains the coefficients of divisorial adjunction, including on the normalization of a nonnormal floor component \cite[Proposition~2.7]{Wimp}. For finite $I$ it is hyperstandard with a finite seed, since only finitely many sums $f\leq1$ occur.

\begin{lemma}
\label{lem:uniform-acc-consequences}
For a rational DCC set $J\subset[0,1]$ there are positive constants $\eta_{\mathrm{lc}}(J)$ and $\eta_{\mathrm{num}}(J)$, uniform in characteristic $p>5$, with the following properties.
\begin{enumerate}
\item For a $\Q$-factorial klt threefold $V$ over an algebraically closed field and $C\geq0$ with $\coeff C\subset J$, if $(V,(1-\eta)C)$ is lc for $0<\eta<\eta_{\mathrm{lc}}(J)$, then $(V,C)$ is lc.
\item Let $(T,C)$ be an lc projective pair of dimension one or two over any field of characteristic $p>5$, with $\coeff C\subset J$. If $C'$ is rational with $K_T+C'$ $\Q$-Cartier and
\[
 (1-\eta)C\leq C'\leq C,\qquad K_T+C'\equiv0,\qquad
 0<\eta<\eta_{\mathrm{num}}(J),
\]
then $C'=C$.
\end{enumerate}
\end{lemma}
\begin{proof}
For (1), counterexamples give thresholds $\lct(V,0;C)$ approaching one from below, contrary to uniform ACC \cite[Theorem~1.5]{Wimp}. For (2), take counterexamples with $\eta_i\to0$. The union of the coefficient sets of $C_i'$ is DCC: a strictly decreasing sequence would admit a subsequence with indices tending to infinity, and the corresponding coefficients of $C_i$ would approach the same limit from above. Global ACC for curves and surfaces \cite[Proposition~5.4]{DasVol}, \cite[Theorem~1.4]{Wimp}, puts the $C_i'$-coefficients in a finite set. For surfaces use a crepant dlt modification, which adds only primes of coefficient one \cite[Theorem~4.7 and Remark~4.8]{Texc}. Along a subsequence, a component of $C_i-C_i'$ then has a fixed $C_i'$-coefficient $c>0$, and its $C_i$-coefficient lies in $(c,c/(1-\eta_i)]$, contradicting DCC for $J$.
\end{proof}

\begin{proposition}
\label{prop:global-nonklt-rounding}
Fix a rational DCC set $I\subset[0,1]$ containing $0,1$. There is a rational $\eta=\eta(I)>0$, uniform in characteristic $p>5$, with the following property. Let $(X,B=\sum b_jB_j)$ be a projective $\Q$-factorial lc pair of Fano type of dimension three, with $\coeff B\subset I$, $\floor B\ne0$ and $N=-(K_X+B)$ nef. Set
\[
 \Theta=\sum_{b_j\leq1-\eta} b_jB_j+\sum_{b_j>1-\eta} B_j.
\]
An MMP for $-(K_X+\Theta)$ ends birationally at a $\Q$-factorial lc pair $(W,\Theta_W)$ of Fano type with nef negative log divisor, preserving every component of $\floor\Theta$. Its monotonic $n$-complements transfer to $(X,B)$ at the same index. For $I=\Phi(R)$ the coefficients of $\Theta_W$ belong to the fixed finite set $J(R)=(\Phi(R)\cap[0,1-\eta])\cup\{1\}$.
\end{proposition}
\begin{proof}
Choose $0<\eta<1/2$ below $\eta_{\mathrm{lc}}(I)$, $\eta_{\mathrm{num}}(I)$ and $\eta_{\mathrm{num}}(I^{\mathrm{adj}})$. The divisor $N$ is semiample by \cite[Corollary~1.5]{ABL} and \cref{lem:numerical-torsion}. Complete $B$ to an lc log Calabi--Yau boundary $B+P$ by semiample perturbation. Run the MMP as in \cref{prop:nonexceptional-big-reduction}. The completion stays crepant and the birational models remain of Fano type. Since $(1-\eta)\Theta_i\leq B_i\leq\Theta_i$, \cref{lem:uniform-acc-consequences} makes every $(X_i,\Theta_i)$ lc. At a divisorial step $f:X_i\to X_{i+1}$ contracting $E$, the coefficient of $E$ in $K_{X_i}+\Theta_i-f^*(K_{X_{i+1}}+\Theta_{i+1})$ is negative by negativity for $-(K_{X_i}+\Theta_i)$. If $E\subset\floor{\Theta_i}$, its crepant coefficient for the target would exceed one.

Suppose there is a Mori fibre endpoint $g:W\to T$. Since $P_W$ is effective and $K_W+B_W+P_W\sim_\Q0$, $K_W+B_W$ has nonpositive degree on the fibre ray, whereas $K_W+\Theta_W$ has positive degree. For some rational $0\leq u<1$, the lc boundary $\Gamma=B_W+u(\Theta_W-B_W)$ has $K_W+\Gamma\equiv_T0$. It satisfies $(1-\eta)\Theta_W\leq\Gamma\leq\Theta_W$. For $\dim T>0$, part (2) of \cref{lem:uniform-acc-consequences} on the normal generic fibre gives $\Gamma|_{W_{k(T)}}=\Theta_W|_{W_{k(T)}}$, contrary to positivity on the fibre ray.

If $T$ is a point, take a surviving component $S$ of $\floor B$. Adjunction to its normalization gives lc boundaries $\Gamma_S,\Theta_S$. The formula $\Diff_{S^\nu}(\Gamma-S)=\Diff_{S^\nu}(0)+(\Gamma-S)|_{S^\nu}$ and its analogue for $\Theta_W$ give
\[
 (1-\eta)\Theta_S\leq\Gamma_S\leq\Theta_S,\qquad
 K_{S^\nu}+\Gamma_S\equiv0,\qquad
 \coeff\Theta_S\subset I^{\mathrm{adj}}.
\]
Part (2) gives equality, but $\Theta_W-\Gamma$ is a nonzero effective divisor, ample because $\rho(W)=1$. Its restriction to $S^\nu$ cannot be zero. Thus the MMP has a nef output. On a common resolution, negativity gives $v^*(K_W+\Theta_W)\geq u^*(K_X+\Theta)\geq u^*(K_X+B)$. Transfer the defining rational function as in \cref{prop:nonexceptional-big-reduction}.
\end{proof}

\begin{corollary}
\label{cor:finite-nonklt-problem}
For bounded complements of non-klt nef-and-big pairs with coefficients in $\Phi(R)$, it suffices to treat $\Q$-factorial lc Fano-type threefolds with nonempty floor, coefficients in $J(R)$ and nef and big negative log divisor.
\end{corollary}
\begin{proof}
Take a crepant dlt modification, of Fano type by the mixture in \cref{prop:nonexceptional-big-reduction}, and apply \cref{prop:global-nonklt-rounding}. If its nef output is non-big, use \cref{thm:global-nonklt,cor:nonbig} with seed $R$, since $J(R)\subset\Phi(R)$. Then transfer the complement back.
\end{proof}

\begin{proposition}
\label{prop:global-plt-preparation}
Let $(X,B)$ be a rational projective non-klt lc threefold pair in characteristic $p>5$, with $X$ of Fano type and $N=-(K_X+B)$ nef and big. There is a projective birational morphism $h:Y\to X$ from a $\Q$-factorial threefold of Fano type, with an effective boundary $B_Y$ such that
\[
 K_Y+B_Y=h^*(K_X+B),\qquad \coeff B_Y\subset\coeff B\cup\{1\}.
\]
There is a rational plt boundary $\Gamma$ on $Y$ with prime floor $S\subset\floor{B_Y}$ and $-(K_Y+\Gamma)$ ample. In particular $S$ is a normal surface of Fano type, and complements descend at the same index.
\end{proposition}
\begin{proof}
Start with a crepant dlt modification of Fano type, as in \cref{cor:finite-nonklt-problem}, and rename it $(X,B)$. For a klt log Fano boundary $C$, choose small rational $\tau>0$ and a Kodaira decomposition $N-\tau C\sim_\Q A+G$, with $A$ ample and $G\geq0$. For $0<\varepsilon<\tau/2$ sufficiently small set
\[
 \widetilde B=(1-\varepsilon)B+\varepsilon C,\quad
 G_\varepsilon=G+(\tau-\varepsilon)C,\quad
 H_\varepsilon=G_\varepsilon+\varepsilon B.
\]
The pair $(X,\widetilde B)$ is klt with ample negative log divisor $N_{\widetilde B}\sim_\Q A+H_\varepsilon$. Choose $\varepsilon$ small enough that $A+\varepsilon C$ is ample.

If $\Supp(G+C)$ contains no lc centre of $(X,B)$, then $\Omega=B+\delta G_\varepsilon$ is dlt for small $\delta>0$. Its negative log divisor is equivalent to $(1-\delta)N+\delta(A+\varepsilon C)$, hence ample. Otherwise fix an lc place $E$ with centre in $\Supp(G+C)$. The quotient $a(E,X,\widetilde B)/\ord_E H_\varepsilon$ is $O(\varepsilon)$, so $u_\varepsilon=\lct(X,\widetilde B;H_\varepsilon)<1$ tends to zero. The lc, non-klt boundary $\Omega=\widetilde B+u_\varepsilon H_\varepsilon$ tends to $B$. Its negative log divisor is equivalent to $(1-u_\varepsilon)N_{\widetilde B}+u_\varepsilon A$, hence ample. On a fixed log resolution of $B+C+G$, for small $\varepsilon$ every coefficient-one component for $\Omega$ already has coefficient one for $B$. The lc places of the log smooth subpairs come from coefficient-one strata, so every lc place of $\Omega$ is one of $B$.

If $\floor\Omega$ has a component $T$, it is in the dlt floor of $B$, so $(X,T)$ is plt. A small convex perturbation of $\Omega$ toward $T$ gives $\Gamma$. Otherwise, take a crepant dlt modification $\rho:(Z,\Omega_Z)\to(X,\Omega)$, with $F=\floor{\Omega_Z}=\Exc(\rho)_{\mathrm{red}}$. Run the $(K_Z+F)$-MMP over the $\Q$-factorial $X$ \cite[Theorem~1.1]{HWRelative}. The difference from $\rho_i^*K_X$ on each model is effective exceptional, with every exceptional coefficient equal to the positive discrepancy $a(E_i,X,0)$. Negativity makes it vanish on the nef output. This output is small over $X$ and therefore equals $X$ by $\Q$-factoriality. The last step is a divisorial contraction $g:Y\to X$ with exceptional prime $T$, $(Y,T)$ plt and $-(K_Y+T)$ $g$-ample.

The $\Omega$-pair remains crepant over $X$, so $T$ is an lc place of $B$ and $B_Y=B^\sim+T$. For small $\lambda>0$, $\Gamma=(1-\lambda)\Omega_Y+\lambda T$ is plt and its negative log divisor is ample, being $(1-\lambda)g^*(-(K_X+\Omega))+\lambda (-(K_Y+T))$. Lowering the floor gives Fano type. Adjunction to $T$ gives a klt log Fano surface. Composing the morphisms gives the stated crepant model and principal-function descent preserves the index.
\end{proof}

\begin{proposition}
\label{prop:global-walls}
Let $(X,B)$ be a projective $\Q$-factorial lc threefold pair of Fano type in characteristic $p>5$, with $N=-(K_X+B)$ nef and big. Suppose $\Gamma$ is a rational plt boundary with prime floor $S\subset\floor B$ and $A=-(K_X+\Gamma)$ ample. Fix $n>0$ with $D=nN$ integral Weil. There is a finite $-S$-MMP
\[
 X=X_0\dashrightarrow\cdots\dashrightarrow X_l=W,\qquad g:W\to T,
\]
ending with a Mori fibre space and having the following properties.
\begin{enumerate}
\item The successive birational contractions occur at rational parameters $0\leq r_1\leq\cdots\leq r_l\leq r_*$, with $r_*>0$ the parameter of $g$. At each contraction with parameter $r$, $D_i-rS_i$ is nef and descends through the contraction. Between successive parameters the corresponding divisors are nef.
\item The $X_i$ are of Fano type, and the transforms $\Omega_i$ of an effective plt log Calabi--Yau boundary on $X$ are crepant plt with floor $S_i$ and $\Omega_i-S_i$ big.
\item Each $S_i$ is normal with a klt log Fano boundary $\Xi_i\geq\Delta_i^0:=\Diff_{S_i}(0)$.
\item The divisor $S_W$ is $g$-ample and $H_*:=D_W-r_*S_W\sim_\Q g^*A_T$ for a nef rational Cartier divisor $A_T$ on $T$.
\end{enumerate}
\end{proposition}
\begin{proof}
Take $Q\sim_\Q2A$ effective with $(X,\Gamma+Q)$ lc and put $\Omega^0=\Gamma+Q/2$. Then $\Omega^0$ is plt log Calabi--Yau and $\Theta^0=\Omega^0-\epsilon S$ is big klt for $0<\epsilon<1$, with $K_X+\Theta^0\sim_\Q-\epsilon S$. The models remain of Fano type. The centre $S$ survives since it is positive on the chosen rays, whereas a contracted prime has negative degree on those rays. The $\Omega^0$-pair thus stays crepant plt by \cref{lem:plt-preserved}.

The nef and big $D$ defines a birational contraction to a normal variety $V$ by \cite[Corollary~1.5]{ABL}. First run the $-S$-MMP over $V$ with scaling of an effective representative of $cA$ chosen to retain klt. This terminates by \cite[Theorem~1.6]{BW}. Its contractions are at parameter zero, with $D_i$ pulled back from $V$. The scaling equation $(-\epsilon S_i+\lambda cA_i)\cdot R=0$ with $\lambda>0$ gives $A_i\cdot R>0$.

The output is $-S$-nef over $V$. Since the cone of curves of a Fano-type model is rational polyhedral, $D-\delta S$ is nef and big there for small rational $\delta>0$. By semiample perturbation choose an effective representative $C\sim_\Q(\epsilon/\delta)D$ so that the addition to the transform of $\Theta^0$ is klt. Run the MMP with scaling of $C$ using \cite[Theorem~1.5]{BW}. Its nef divisors at the scaling parameters $\lambda>0$ are positive multiples of $D_i-rS_i$, $r=\delta/\lambda$. These increase in parameter as $\lambda$ decreases, and nefness in each intervening interval follows by convexity. The nef supports are numerically trivial along the corresponding contractions and descend by relative base-point-freeness. A rational principal adjustment downstairs makes each descent an equality of divisors. As $-S_i$ has negative intersection with the square of an ample divisor, the program ends with a Mori fibre space, where $S_W$ is relatively ample and $D_W\cdot R=r_*S_W\cdot R>0$. Relative base-point-freeness gives the stated descent of $H_*$; for a point base use \cref{lem:numerical-torsion}.

To control the restrictions to the centres, take rational $a<1$ close to one and put $\Gamma_a=aB+(1-a)\Gamma$, $A_a=-(K_X+\Gamma_a)=aN+(1-a)A$. The pair is plt with floor $S$ and $A_a$ is ample. At zero parameters $N_i\cdot R=0$, $A_i\cdot R>0$; at positive parameters $N_i\cdot R>0$. There are finitely many steps, so we can ensure $A_{a,i}\cdot R>0$ at every step. As with $\Omega^0$, complete $\Gamma_a$ by a sufficiently diluted $P\sim_\Q A_a$ to a plt log Calabi--Yau boundary $\Omega=\Gamma_a+P$. Its transforms stay crepant plt and $\Omega_i-S_i$ remains big since $P$ is big.

At a birational step, for maps $u,v$ from a common resolution to the source and target, negativity gives $v^*A_{a,i+1}=u^*A_{a,i}+E$, $E\geq0$ \cite[Lemma~2.3]{BW}. The exceptional divisor $E$ does not contain the strict transform of $S_i$, so the restriction of $A_{a,i}$ to $S_i$ stays big. In particular
\[
 \Omega_{S_i}-\Delta_i^0=(\Gamma_{a,i}-S_i+P_i)|_{S_i}
\]
is effective and big since it contains $P_i|_{S_i}\sim_\Q A_{a,i}|_{S_i}$. Write this difference as $Q_i\sim_\Q A_i^S+E_i^S$ with $A_i^S$ ample and $E_i^S\geq0$. For small $v_i>0$, the boundary $\Xi_i=\Delta_i^0+(1-v_i)Q_i+v_iE_i^S$ has the required properties.
\end{proof}

\begin{theorem}[Reduction of the lifting vanishing]
\label{thm:global-wall-vanishing}
With the data of \cref{prop:global-walls}, put $M=\lceil r_*\rceil$, $a_*=M-r_*\in[0,1)$ and $L_*=D_W-MS_W=H_*-a_*S_W$. There is a surjection
\[
 H^1(W,\cO_W(L_*))\longrightarrow H^1(X,\cO_X(D-S)).
\]
If $T$ is a point, the group on the left vanishes.
\end{theorem}
\begin{proof}
For any nef integral Weil $H$ on $X_i$, the actual restriction quotient $\cO_{S_i}(G_H)$ of \cref{lem:restriction} satisfies $G_H\sim_\Q H|_{S_i}-\delta_H$ with $0\leq\delta_H\leq\Delta_i^0$. The boundary $\Xi_i-\delta_H$ is effective klt and
\[
 G_H-(K_{S_i}+\Xi_i-\delta_H)\sim_\Q
 H|_{S_i}-(K_{S_i}+\Xi_i)
\]
is ample. Hence $H^{>0}(S_i,\cO_{S_i}(G_H))=0$ by \cite[Theorem~1.1]{ABL}.

At parameter $r$, \cite[Lemma~4.3]{Fib} identifies cohomology of the integral Weil divisors $(D_i-rS_i)-bS_i$, $0\leq b\leq1$, on the two models. The supporting divisor $D_i-rS_i$ descends, and $bS_i$ accounts for the fractional correction along the plt centre. For a flip, apply it on each side of the common contraction target.

For $1\leq j<M$, take the model $X^{(j)}$ after the birational steps of parameter at most $j$, and set $X^{(M)}=W$. The divisor $D_j-jS_j$ on $X^{(j)}$ is nef. Its restriction sequence and the preceding surface vanishing give
\[
 H^1(X^{(j)},\cO(D_j-(j+1)S_j))
 \twoheadrightarrow H^1(X^{(j)},\cO(D_j-jS_j)).
\]
Transport identifies the source with the corresponding group on $X^{(j+1)}$, since each intervening parameter $r\in(j,j+1]$ has $0\leq j+1-r<1$. At the initial steps $r\in[0,1]$, transport applies to $D_i-S_i$ with $b=1-r$. Composing gives the surjection, including when $M=1$. If $T$ is a point, $(W,S_W)$ is plt and $-(K_W+S_W)\sim_\Q\Omega_W-S_W$ is ample by bigness and $\rho(W)=1$. Then \cref{thm:rank-one-vanishing} applies to the integral Weil divisor $L_*$.
\end{proof}

\subsection{Vanishing on Mori fibre spaces and lifting complements}\label{s:nonex-lifting}\label{sec:curve-residual}\label{sec:stein-factor}\label{s-sec:trace}
When the Mori base is a surface, we use relative vanishing for contractions with one-dimensional fibres.

\begin{lemma}
\label{o-lem:vanishing}
Let $f:X\to Z$ be projective with $X$ a normal $\Q$-factorial klt threefold over an algebraically closed field of characteristic $p>5$. If $P$ is an integral Weil divisor with
\begin{align}
 \dim X_z&\leq1 \quad\text{for every }z\in Z,
 \label{eq:display-2-4-a}\\
 A&:=P-K_X\quad\text{$f$-ample},
 \label{eq:display-2-4-c}
\end{align}
then
\begin{equation}
 R^1f_*\cO_X(P)=0.
 \label{eq:display-2-4-3}
\end{equation}
\end{lemma}
\begin{proof}
The strongly $F$-regular locus $U\subset X$ contains the codimension-two points by \cite[Theorem~1.2]{SatoTakagi}, so its complement is finite. Choose $E\geq0$ integral Weil with $E-K_X$ Cartier and write a Cartier index of $K_X$ as $mp^a$, $p\nmid m$. For large $\ell>a$ put
\begin{align}
 \Theta_\ell&=E/(p^\ell-1),
 \label{eq:display-2-4-4}\\
 m(p^\ell-1)(K_X+\Theta_\ell)&=m\bigl[p^\ell K_X+(E-K_X)\bigr],
 \label{eq:display-2-4-5}\\
 A_\ell&=P-(K_X+\Theta_\ell)=A-\Theta_\ell.
 \label{eq:display-2-4-6}
\end{align}
The middle divisor is Cartier, so the log index is prime to $p$. By small effective perturbation and quasi-compactness, $(U,\Theta_\ell)$ is strongly $F$-regular, and $A_\ell$ remains $f$-ample. Fix $\Theta_\ell=\widetilde\Theta$ and $A_\ell=\widetilde A$. A finite cover of $U$ and transitivity of traces give an exponent $e>0$ at which $(p^e-1)(K_X+\widetilde\Theta)$ is Cartier and
\begin{equation}
 \Tr_e:F^e_*\cO_X((1-p^e)(K_X+\widetilde\Theta))\to\cO_X
 \label{eq:display-2-4-7}
\end{equation}
and its iterates are surjective on $U$. For $q=p^{ne}$ set
\begin{equation}
 M_q=(1-q)(K_X+\widetilde\Theta)+qP=P+(q-1)\widetilde A.
 \label{eq:display-2-4-8}
\end{equation}
Twist the trace on the regular locus and extend reflexively to get
\begin{equation}
 \tau_{q,P}:F^{ne}_*\cO_X(M_q)\to\cO_X(P).
 \label{eq:display-2-4-9}
\end{equation}
Tensoring the trace with $\cO_X(P)$ gives a surjection on $U$, which factors through $\tau_{q,P}$ by the finite projection formula and the natural multiplication of divisorial sheaves. Thus the cokernel $\mathcal C_q$ is supported on $X\setminus U$. With kernel $\mathcal K_q$ and image $\mathcal I_q$, the sequences are
\begin{align}
 0&\to\mathcal K_q\to F^{ne}_*\cO_X(M_q)\to\mathcal I_q\to0,
 \label{eq:display-2-4-10}\\
 0&\to\mathcal I_q\to\cO_X(P)\to\mathcal C_q\to0.
 \label{eq:display-2-4-11}
\end{align}
For $r\widetilde A$ Cartier, writing $q-1=a_qr+b_q$, $0\leq b_q<r$, expresses $\cO_X(M_q)$ as one of finitely many coherent sheaves twisted by $a_qr\widetilde A$. Relative Serre vanishing and the finiteness of Frobenius give, for $q\gg0$,
\begin{align}
 R^1f_*\cO_X(M_q)&=0,
 \label{eq:display-2-4-12}\\
 R^1f_*F^{ne}_*\cO_X(M_q)&=F^{ne}_{Z*}R^1f_*\cO_X(M_q)=0.
 \label{eq:display-2-4-13}
\end{align}
The fibre bound gives
\begin{equation}
 R^2f_*\mathcal K_q=0.
 \label{eq:display-2-4-14}
\end{equation}
The sequences yield $R^1f_*\mathcal I_q=0$ and then $R^1f_*\cO_X(P)=0$, since $R^1f_*\mathcal C_q=0$.
\end{proof}

\begin{proposition}
\label{prop:surface-base-endpoint}
In the setting of \cref{thm:global-wall-vanishing}, suppose that $\dim T=2$. Then $a_*=0$ and $H^1(W,\cO_W(L_*))=0$.
\end{proposition}
\begin{proof}
Write $S=S_W$ and $\Omega_W=S+\Delta$. The generic fibre $F$ over $K=k(T)$ is regular, with $H^0(F,\cO_F)=K$. Set $d=\deg S_F>0$ and $e=\deg\Delta_F>0$, the latter by bigness of $\Delta$. The log Calabi--Yau equality gives
\[
 2h^1(F,\cO_F)-2+d+e=0.
\]
Thus $h^1=0$ and $d=1$, so $F\simeq\PP^1_K$ by Riemann--Roch \cite[Tag 0C6U]{Stacks}. In particular $\mu=g|_S:S\to T$ is birational. The degree of the integral divisor $L_*|_F$ is $-a_*d$, so $a_*=0$.

Every fibre has dimension one. Otherwise a prime surface $E$ over a closed point of $T$ is numerically trivial over $T$ since $\rho(W/T)=1$ and it misses $F$. Relative base-point-freeness makes a Cartier multiple semiample, and its contraction is constant on the fibres of $g$ because $E\equiv_T0$. Its Stein target is then finite birational over $T$, so $E\sim_\Q g^*A_E$. The rational function giving this equivalence restricts to a constant on $F$, hence comes from $K^*$. Absorbing its divisor downstairs expresses $E$ as an actual pullback, which cannot be supported over a point.

By a principal adjustment, make $L_*$ vertical using $F\simeq\PP^1_K$ and degree zero. The same argument with the rational function gives an equality $L=g^*A$, where $L$ is the adjusted integral divisor and $A$ is nef rational Cartier. Let $\cO_S(G)$ be the restriction quotient for $L$. As subsheaves of $K$, we have
\begin{equation}
 g_*\cO_W(L)=\cO_T(\floor A)=\mu_*\cO_S(G).
 \label{eq:section-pushforward}
\end{equation}
Indeed, sections of $\cO_W(L)$ are constant on $F$. For functions from $K$, the divisor inequality on $W$ is equivalent to the one on $T$ by pullback, giving the first equality and the inclusion into the right-hand sheaf. Conversely a local section $s$ of $\cO_S(G)$ lifts locally to $\cO_W(L)$. If $qA$ is Cartier, the $q$-th power of the lift restricts to a section of $\cO_S(\mu^*(qA))$. At the generic point over a prime $P$ of $T$, where $\mu$ is an isomorphism, this yields $q\ord_P(s)+q\coeff_PA\geq0$. Since $\ord_P(s)$ is integral, the inequalities at all primes put $s$ in $\cO_T(\floor A)$.

The divisor $L-K_W$ is $g$-ample by its degree two on $F$ and $\rho(W/T)=1$. Thus $R^1g_*\cO_W(L)=0$ by \cref{o-lem:vanishing}. The surface vanishing in \cref{thm:global-wall-vanishing} gives $H^1(S,\cO_S(G))=0$. Leray for $\mu$ gives $H^1(T,\mu_*\cO_S(G))=0$, and \eqref{eq:section-pushforward} and Leray for $g$ prove the result.
\end{proof}

Over a curve the generic fibre is a surface. We prove relative vanishing near each closed fibre, keeping the fractional term $a_*S_W$ in the divisor throughout.

\begin{lemma}
\label{lem:fractional-transport}
Let $h:Y\to X$ be a projective birational morphism of normal threefolds in characteristic $p>5$, with effective plt pairs satisfying $K_Y+\Lambda_Y=h^*(K_X+\Lambda_X)$. If
\[
 h^*H_X=L_Y+\delta_Y,\qquad 0\leq\delta_Y\leq\Lambda_Y,
\]
where $H_X$ is rational Cartier, $L_Y$ is integral Weil and the exceptional coefficients of $\delta_Y$ are less than one, then $h_*\cO_Y(L_Y)=\cO_X(h_*L_Y)$ and $R^{>0}h_*\cO_Y(L_Y)=0$.
\end{lemma}
\begin{proof}
This extends the transport argument of \cite[Lemma~4.3]{Fib} to exceptional fractional coefficients. With $L_X=h_*L_Y$ and $\delta_X=h_*\delta_Y$, the ordinary Fujita transform of $(L_X,\delta_X)$ is $L_Y$, with correction $\delta_Y$ and prescribed pullback sum. This gives the zeroth direct image. It is also the $K$-twisted transform, since over $X$ it is equivalent to $K_Y+\Lambda_Y-\delta_Y$ and the complementary boundary is effective with exceptional coefficients below one.

On a common log resolution, the $K$-twisted transform for $Y$ is also the one for $X$. Its difference from the ordinary transform is integral exceptional, with coefficients equal to the sum of the two nonnegative boundary corrections minus the crepant $\Lambda$-coefficient. These exceed $-1$ by plt crepancy, so the difference is effective and the zeroth images agree. Grauert--Riemenschneider vanishing for the two resolution maps and Leray give the asserted higher images \cite[Theorem~3, Definition~8 and Claims~9--11]{BK}.
\end{proof}

\begin{lemma}
\label{lem:vertical-fractional-preparation}
Let $f:X\to C$ be a contraction of projective varieties to a smooth curve, with $X$ a $\Q$-factorial threefold of Fano type over $C$ in characteristic $p>5$. Suppose $(X,B)$ is rational lc, $-(K_X+B)$ is nef and big over $C$ and the nonempty non-klt locus lies in $X_z$ for a closed point $z$. Let $H$ be rational Cartier with $\coeff_E B>\coeff_E\{H\}$ for every component $E$ of $\{H\}$. There is a crepant projective birational model $h:(Y,B_Y)\to(X,B)$ extracting only lc places and a rational plt boundary $\Gamma$ such that
\[
 \floor\Gamma=V\subset Y_z\text{ is prime},\qquad
 -(K_Y+\Gamma)\text{ is ample over }C,\qquad
 \Gamma\geq\{h^*H\}.
\]
Here $Y$ is $\Q$-factorial and of Fano type over $C$.
\end{lemma}
\begin{proof}
Use the construction of \cref{prop:global-plt-preparation} over $C$. Relative Kodaira decompositions and semiample perturbations use effective representatives after an ample twist from the base. On the first crepant dlt model $h_0:X_0\to X$, the inequality $B_0>\delta_0:=\{h_0^*H\}$ on the support of $\delta_0$ persists, since newly extracted coefficients of $B_0$ are one. In both branches of the construction we can choose $\Omega_0$ arbitrarily close to $B_0$, with all lc places lying over $z$, and with the strict inequalities against $\delta_0$ retained. The final perturbation toward one floor component can likewise retain these inequalities. When an extra prime $V$ is extracted instead, its $\Gamma$-coefficient is one, which dominates the fractional coefficient there. The relative ampleness and Fano type follow from the construction.
\end{proof}

\begin{proposition}
\label{prop:curve-fractional-relative-vanishing}
Let $f:X\to C$ be a contraction of projective varieties with $X$ a $\Q$-factorial threefold of Fano type over a smooth curve, over an algebraically closed field of characteristic $p>5$. Let $(X,S)$ be plt with $S$ prime and both $S$ and $-(K_X+S)$ $f$-ample. If $L$ is integral Weil and $H=L+aS\sim_{\Q,C}0$, $0\leq a<1$, then $R^if_*\cO_X(L)=0$ for $i>0$.
\end{proposition}
\begin{proof}
Fix $z\in C$ closed, write $F=f^*z$, and choose rational $a<b<1$. Set $t=\lct(X,bS;F)$. The lc pair $B=bS+tF$ has non-klt locus over $z$ and relatively ample negative log divisor. Apply \cref{lem:vertical-fractional-preparation} to $\{H\}=aS$, and use its $h:Y\to X$, $V$ and $\Gamma$. Write $A=-(K_Y+\Gamma)$, $L_Y=\floor{h^*H}$ and $\delta_Y=\{h^*H\}$. For $u<1$ close to one the klt boundary $bS+utF$ has crepant coefficients $1-(1-u)t\ord_E F$ at the extracted primes, so its crepant transform is effective and dominates $\delta_Y$. Thus \cref{lem:fractional-transport} gives
\begin{equation}
 Rh_*\cO_Y(L_Y)\simeq\cO_X(L).
 \label{eq:curve-initial-transport}
\end{equation}

Complete $\Gamma$ to a relative plt log Calabi--Yau boundary $\Omega=\Gamma+P$, $P\sim_{\Q,C} A$. Run a $-V$-MMP with scaling of a representative of $cA$ whose addition to $\Omega-\varepsilon V$ is klt with nef log divisor over $C$. It terminates by \cite[Theorem~2.12]{GNT}. Since $(fh)^*z=mV+F'$ for $F'\geq0$ avoiding $V$, $-V\sim_{\Q,C}F'/m$ is relatively pseudo-effective. Thus the output $g':Y'\to C$ has $-V'$ nef and $V$ survives. The $\Omega$-pair stays crepant plt throughout. The divisor $H$ descends at every step as a base pullback plus a rational principal divisor. Its correction is bounded by the respective $\Omega$ and is strictly below one. Applying \cref{lem:fractional-transport} at each contraction and on both sides of each flip preserves the derived image on $C$.

A different fibre component meeting $V'$ would contain a fibre curve of positive $V'$-degree, contrary to nefness of $-V'$. Connectedness therefore gives $g'^*z=m'V'$, hence $V'|_{V'}\sim_\Q0$. Writing $S'$ for the pushforward of $S$, the correction on $Y'$ is $aS'+dV'$ with $0\leq d<1$, since the new primes over $X$ lie over $z$ and $V'$ is the only such prime remaining. The scaling equation gives $A_i\cdot R>0$ at every selected ray. The restriction argument in \cref{prop:global-walls} keeps $A_i|_{V_i}$ big. We have
\[
 \Omega_{V'}-\Diff_{V'}(aS')
  =(\Gamma_{Y'}-V'-aS'+P_{Y'})|_{V'}
\]
effective and big, so Kodaira perturbation gives a klt log Fano boundary $\Xi\geq\Diff_{V'}(aS')$.

For every integer $j\geq0$, use the exact restriction sequence for $L_{Y'}-jV'$. Its quotient $\cO_{V'}(G_j)$ has $G_j\sim_\Q-aS'|_{V'}-\epsilon_j$, $0\leq\epsilon_j\leq\Diff_{V'}(0)$, since $H|_{V'}\sim_\Q0$. The boundary $\Xi-aS'|_{V'}-\epsilon_j$ is effective klt and its log divisor differs from $G_j$ by the ample $-(K_{V'}+\Xi)$. By \cite[Theorem~1.1]{ABL}, $H^{>0}(V',\cO_{V'}(G_j))=0$. Compose $m'$ successive restriction maps to obtain surjections
\[
 R^ig'_*\cO_{Y'}(L_{Y'})\otimes\cO_C(-z)
 \longrightarrow R^ig'_*\cO_{Y'}(L_{Y'}),\qquad i>0.
\]
At $z$ these are multiplication by a local parameter, so Nakayama kills the stalk. Transport and \eqref{eq:curve-initial-transport} give the result at every $z$.
\end{proof}

\begin{corollary}
\label{cor:curve-positive-residual}
For a curve base $g:W\to C$ in \cref{thm:global-wall-vanishing}, $R^ig_*\cO_W(L_*)=0$ for $i>0$. If $a_*>0$, then $H^1(W,\cO_W(L_*))=0$.
\end{corollary}
\begin{proof}
The divisor $S_W$ is relatively ample, as is $-(K_W+S_W)\sim_\Q\Omega_W-S_W$ by its bigness on the generic fibre and $\rho(W/C)=1$. Apply \cref{prop:curve-fractional-relative-vanishing}. For $a_*>0$ the generic restriction is anti-ample and has no sections, so the torsion-free sheaf $g_*\cO_W(L_*)$ is zero. Use Leray.
\end{proof}

\begin{lemma}
\label{lem:curve-trivial-stein}
At a curve Mori endpoint in \cref{thm:global-wall-vanishing}, if $a_*=0$ and $\mu_*\cO_{S_W}=\cO_C$ for $\mu=g|_{S_W}$, then $H^1(W,\cO_W(L_*))=0$.
\end{lemma}
\begin{proof}
If $L_*$ has no generic section, its pushforward is zero and \cref{cor:curve-positive-residual} applies. Otherwise its numerically trivial generic restriction is principal, since the effective divisor of a section must be zero. A principal adjustment makes $L_*$ vertical, giving $L=g^*A$ as in \cref{prop:surface-base-endpoint}, for $A$ nef rational Cartier. Here the principal function descends because its restriction to the projective generic fibre has no zeros or poles and $g$ is a contraction. In particular $g_*\cO_W(L)=\cO_C(\floor A)$.

Let $\cO_S(G)$ be the restriction quotient on $S=S_W$. We also have $\mu_*\cO_S(G)=\cO_C(\floor A)$: the Stein hypothesis identifies its generic sections with $k(C)$. For a local section $s$, the power argument of \cref{prop:surface-base-endpoint} gives, on a component $E$ of $\mu^*z$ of multiplicity $e_E>0$,
\[
 qe_E(\ord_z(s)+\coeff_z A)\geq0
\]
for $qA$ Cartier. This yields the required inclusion, and the reverse is restriction from $W$. The surface vanishing from \cref{thm:global-wall-vanishing} gives $H^1(S,\cO_S(G))=0$. Leray for $\mu$ and $g$ now gives $H^1(W,\cO_W(L))=0$.
\end{proof}

For a horizontal plt centre over a curve, we also need to determine the finite map in its Stein factorization.

\begin{lemma}
\label{lem:standard-curve-finite-field}
Let $T$ be a normal projective integral curve over an $F$-finite field $K$ of characteristic $p>5$. If $(T,\Delta)$ is klt log Fano and the coefficients of $\Delta$ are standard, then it is globally $F$-regular.
\end{lemma}
\begin{proof}
Regard $T$ as a curve over its finite constant field $K'=H^0(T,\cO_T)$. The canonical classes over $K$ and $K'$ agree by duality for the finite morphism $\operatorname{Spec}K'\to\operatorname{Spec}K$, whose relative dualizing module is one-dimensional over $K'$. Apply \cite[Proposition~5.5]{SatoTakagi} to the log Fano pair over $K'$ with standard coefficients. Global $F$-regularity uses the same absolute Frobenius over either field.
\end{proof}

\begin{lemma}
\label{lem:constant-trace-surjection}
Let $V$ be a normal projective $\Q$-factorial surface over an $F$-finite field of positive characteristic. Let $(V,T)$ be plt with $T$ a normal prime and $-(K_V+T)$ ample. If $(T,\Diff_T(0))$ is globally $F$-regular, then $H^0(V,\cO_V)\to H^0(T,\cO_T)$ is surjective.
\end{lemma}
\begin{proof}
Write $D=K_V+T$, choosing its representative to have coefficient zero along $T$. Choose a Cartier divisor $C_D\geq D$ whose coefficient along $T$ is also zero, using global generation of $\cO_V(-D)$ twisted by a sufficiently ample line bundle. Put $E=C_D-D$. For $N$ divisible by $p$ with $ND$ Cartier, take $q=1+kN$ large and let $\Lambda_q=E/q$. Then $q(D+\Lambda_q)=kND+C_D$ is Cartier and $p\nmid q$. For large $q$, the pair $(V,T+\Lambda_q)$ is plt log Fano and $\Delta_q=\Diff_T(0)+E|_T/q$ is its different.

Global $F$-regularity on $T$ gives a splitting at some exponent $e_0$ with extra integral divisor $G_T\geq E|_T$. Increase $q$ so that $(p^{e_0}-1)E|_T/q\leq G_T$. This gives a splitting for $\Delta_q$ at $e_0$, and iteration gives one at all multiples of $e_0$, with the iterated divisor dominating the rounded boundary. Take arbitrarily large such $e$ with $q\mid p^e-1$. For the ample Cartier divisor $M_e=(1-p^e)(D+\Lambda_q)$, Frobenius adjunction identifies the quotient trace along $T$ with the trace defined by $\Delta_q$ \cite[Theorem~5.3]{Das}. Thus the following diagram commutes:
\[
 \begin{array}{ccc}
 H^0(V,F^e_*\cO_V(M_e))&\longrightarrow&
 H^0(T,F^e_*\cO_T(M_e|_T))\\
 \downarrow\varphi_V&&\downarrow\varphi_T\\
 H^0(V,\cO_V)&\longrightarrow&H^0(T,\cO_T).
 \end{array}
\]
The trace preserves the ideal of $T$ by its coefficient-one formula at the generic point of $T$, and equality of the Frobenius different and ordinary different gives the quotient trace on the normal curve. The right arrow is surjective by the splitting. For large $e$ the top arrow is surjective by Serre vanishing for $\cO_V(-T)$ twisted by $M_e$, which is a growing multiple of a fixed ample Cartier divisor. The bottom arrow follows.
\end{proof}

\begin{proposition}[The Stein factor of a plt centre]
\label{prop:horizontal-plt-stein}
Let $f:X\to C$ be a projective contraction from a normal $\Q$-factorial threefold to a smooth curve over an algebraically closed field of characteristic $p>5$. Suppose $(X,S)$ is plt with $S$ a horizontal prime and $-(K_X+S)$ $f$-ample. Then $(f|_S)_*\cO_S=\cO_C$.
\end{proposition}
\begin{proof}
Over $K=k(C)$, the normal $\Q$-factorial plt surface pair $(V,T)=(X_K,S_K)$ has $-(K_V+T)$ ample and $H^0(V,\cO_V)=K$. The klt log Fano adjoint curve has standard coefficients \cite[Lemma~3.3]{HW}. Since $K$ is $F$-finite, \cref{lem:standard-curve-finite-field,lem:constant-trace-surjection} give $H^0(T,\cO_T)=K$. The finite map in the Stein factorization of $S\to C$ is birational, hence an isomorphism by normality of $C$.
\end{proof}

\begin{corollary}
\label{cor:global-wall-complement}
For a finite rational set $J\subset[0,1]$ there is an $n(J)>0$ with the following property. Let $(X,B),\Gamma,S,N$ satisfy the hypotheses of \cref{prop:global-walls} with $\coeff B\subset J$. If $H^1(X,\cO_X(n(J)N-S))=0$, then $(X,B)$ has a monotonic $n(J)$-complement.
\end{corollary}
\begin{proof}
By adjunction the pair $(S,B_S)$ is lc with nef negative log divisor and coefficients in $J^{\mathrm{adj}}$. The auxiliary plt log Fano pair makes $S$ of Fano type. Choose $n=n(J)$ divisible by the denominators of $J$ and the complement index for the finite seed of $J^{\mathrm{adj}}$ \cite[Corollary~3.26]{Fib}. Then $nB$ is integral and $(S,B_S)$ has a monotonic $n$-complement $B_S^+$.

Use the resolution and restriction calculation in \cite[proof of Proposition~4.9]{Fib}. On its log resolution $\phi:X'\to X$, the integral divisor $L'+P'$, with $P'$ effective exceptional, has
\[
 \phi_*\cO_{X'}(L'+P'-S')=\cO_X(nN-S),\qquad
 R^{>0}\phi_*\cO_{X'}(L'+P'-S')=0.
\]
These statements use the comparison of the integral divisorial sheaves and Grauert--Riemenschneider vanishing. They hold with base a point. Leray and the assumed global $H^1$-vanishing make the restriction to $S'$ surjective on global sections. The effective integral residue divisor determined by $B_S^+$ thus lifts, giving $G\geq0$, $G\sim nN$, with $B^+=B+G/n$ and $\Diff_S(B^+-S)=B_S^+$ by the exact adjunction calculation of the same construction.

Inversion of adjunction makes $(X,B^+)$ lc near $S$ \cite[Corollary~10.1]{BMP}. If it were non-lc elsewhere, then $(X,cB^++(1-c)\Gamma)$ for $c<1$ close to one would be plt near $S$, non-lc elsewhere and have ample negative log divisor. Its non-klt locus would have $S$ as an open and closed component, contrary to connectedness \cite[Theorem~1.2]{NT}. Hence $(X,B^+)$ is lc, and integrality and monotonicity give the rounding inequality.
\end{proof}

\begin{theorem}[Non-exceptional complements]\label{th:nonexceptional}
Fix a finite $R\subset[0,1]\cap\Q$. In sufficiently large characteristic, depending only on $R$, let $(X,B)$ be a projective lc threefold with $X$ of Fano type, $\coeff B\subset\Phi(R)$, and $L=-(K_X+B)$ nef and $\Q$-Cartier. If some effective rational $P\simq L$ makes $(X,B+P)$ non-klt, then $(X,B)$ has a monotonic complement of index depending only on $R$. This includes every non-klt pair under these hypotheses.
\end{theorem}
\begin{proof}
If $0<\nu(L)<3$, use \cref{cor:nonbig}. If $L\equiv0$, the only possible $P$ is zero, by intersection with an ample square. Thus $(X,B)$ is non-klt and \cref{thm:global-nonklt} bounds the order of $K_X+B$.

Suppose $L$ is big. For klt $(X,B)$ use \cref{prop:nonexceptional-big-reduction} to obtain a non-klt pair. The reduction carries a complement back at the same index. It also allows an output of smaller numerical dimension, already covered above. For a non-klt pair in the big case, \cref{cor:finite-nonklt-problem} reduces the coefficients to a fixed finite set, and \cref{prop:global-plt-preparation} supplies the coefficient-one surface $S$ and the auxiliary plt pair. Fix the multiple $n$ used to lift the surface complement in \cref{cor:global-wall-complement}.

Apply \cref{prop:global-walls} to the integral Weil divisor $D=nL$ on the prepared model. At a point or surface base, \cref{thm:global-wall-vanishing,prop:surface-base-endpoint} give $H^1(X,\cO_X(D-S))=0$. At a curve base, \cref{prop:horizontal-plt-stein} identifies the Stein factor of $S$. The two cases for the residual coefficient are then covered by \cref{cor:curve-positive-residual,lem:curve-trivial-stein}. The same vanishing follows. Lifting the section on $S$ as in \cref{cor:global-wall-complement} gives a complement on the prepared model. Transfer through the preceding birational steps gives one on the original pair. Choose the maximum of the characteristic bounds and a common multiple of the indices, using \eqref{eq:definition}. Finally, if $(X,B)$ itself is non-klt, every completion is non-klt and semiampleness gives at least one.
\end{proof}

\section{Birational systems and bounded models}\label{s:bounded-models}

\subsection{Anticanonical sections and the boundary threshold}\label{s:trace}\label{s:threshold}
Let $q$ clear the denominators in $R$ (take $q=1$ for $R=\varnothing$). If $q\mid n$ and $b\in\Phi(R)$ with $b<1$, then
\begin{equation}\label{eq:round}
                 \floor{(n+1)b}=\ceil{nb}.
\end{equation}
For $b=1-r/j>0$ write $nr=aj+h$, $0\le h<j$. Since $0<h+r\le j$, both sides equal $n-a$. For $b=0$ both sides vanish. An exceptional pair has no coefficient one. Thus for these pairs and these $n$ one has
\[
 H^0\bigl(X,-nK_X-\ceil{nB}\bigr)\ne0
 \quad\Longleftrightarrow\quad
 (X,B)\text{ has a monotonic }n\text{-complement}.
\]
Indeed an effective member $G$ of the system on the left gives $B^+=(G+\ceil{nB})/n$. It is a completion and hence klt. Conversely $nB^+-\ceil{nB}$ is such a member.

For a nef rational Cartier divisor $A$ at a smooth closed point $x$, the \emph{Seshadri constant}
\[
          \varepsilon(A;x)=\inf_{C\ni x}\frac{A\cdot C}
                                             {\operatorname{mult}_x C}
\]
measures the asymptotic order of jets that sections of $mA$ can separate at $x$, per unit of $m$. On the blowup of $x$ with exceptional divisor $E$, the same number is the supremum of $t$ for which $\pi^*A-tE$ is nef. For an ample $A$, strict lower bounds persist on an open set by jet separation for a sufficiently high Cartier multiple.

\begin{lemma}\label{lem:trace}
Let $M$ be an integral Weil divisor on a normal projective threefold in positive characteristic, and suppose $A=M-(K_X+\Psi)$ is a nef, big and semiample $\Q$-Cartier divisor for an effective rational $\Psi$. If a smooth closed point $x\notin\Supp\Psi$ satisfies $\varepsilon(A;x)>3$, some section of $M$ is nonzero at $x$. If $\varepsilon(A;x)>6$, the system $|M|$ is birational.
\end{lemma}
\begin{proof}
First suppose $A$ is ample. For $q=p^e$ put $D_q=(1-q)K_X+qM-\ceil{(q-1)\Psi}$. On the regular locus, twist the Frobenius trace on the canonical sheaf by $M-K_X$ and restrict the source by $\ceil{(q-1)\Psi}$. Reflexive extension gives
\[
            F^e_*\cO_X(D_q)\longrightarrow\cO_X(M).
\]
Let $\pi$ blow up a finite set $Z$ of smooth closed points off $\Supp\Psi$, with exceptional divisor $E$. Choose $\beta\in\Q_{>0}$ with $\pi^*A-\beta E$ ample, and choose $r$ so that $rA$ is Cartier, $r\Psi$ is integral and $r\beta\in\Z$. Writing $q=rt+j$ with $0\le j<r$ gives $D_q=t(rA)+D_j$, where $D_j$ is given by the same divisor formula. Set $k_q=\ceil{q\beta}$ and $H=\pi^*(rA)-r\beta E$. Then
\[
  \pi^*\cO_X(D_q)(-k_qE)
  =\pi^*\cO_X(D_j)(-\ceil{j\beta}E)\otimes\cO(tH).
\]
The right side has vanishing $H^1$ for $q\gg0$ by Serre vanishing for the finitely many residue sheaves. The pushforward is $\cO_X(D_q)I_Z^{k_q}$ since $\cO_X(D_q)$ is invertible near $Z$ and $\pi_*\cO(-k_qE)=I_Z^{k_q}$. Leray therefore gives surjectivity onto the source jets modulo $I_Z^{k_q}$.

At each point of $Z$, in compatible coordinates and trivializations the trace sends $u=(z_1z_2z_3)^{q-1}$ to $1$ and $z_i^qu$ to $z_i$ \cite[Section~3]{MS}. For $s=1,2$ the monomial containment
\[
  \mathfrak m^{3(q-1)+(s-1)q+1}\ \subset\ (\mathfrak m^s)^{[q]}
\]
shows that errors of this order trace into $\mathfrak m^s$, where $J^{[q]}$ is the ideal generated by $q$th powers of elements of $J$. Thus $\beta>3$ gives arbitrary values at $Z$, and $\beta>4$ gives arbitrary first jets. Taking $Z=\{x\}$ gives nonvanishing when $\varepsilon(A;x)>3$.

For $\varepsilon(A;x)>6$ take an open neighbourhood off the boundary on which the constants exceed $6$. For distinct $y,z$ in this open, the two-point nef threshold is at least half the smaller one-point constant, by comparing the sum of the two curve multiplicities with their maximum. A positive parameter strictly below the nef threshold gives an ample class by interpolation with a small positive parameter. Thus one can take $\beta>3$ on the blowup of $\{y,z\}$ to separate the two points. First jets at $x$ give full differential rank, so the system is generically finite and separable onto its image. Its separation of points on a dense open forces degree one.

If $A$ is nef, big and semiample, a multiple defines a birational contraction $f:X\to Y$ with $rA$ pulled back from an ample Cartier divisor for sufficiently divisible $r$. A point with $\varepsilon(A;x)>0$ lies on no contracted curve. It is isolated in its fibre, so connectedness and normality give an isomorphism above a neighbourhood of $f(x)$. There the Seshadri constants agree by taking strict transforms of curves. Apply the preceding jet calculation on $Y$ to the finitely many sheaves $f_*\cO_X(D_j)$. The projection formula gives the same global source sections on $X$, where we apply the trace.
\end{proof}

The same jet count bounds the volume of an exceptional Fano threefold with zero boundary. Put $L=-K_X$ and fix a smooth point $x$. The blowup of $x$ has log discrepancy three. For every $G\in|mL|$ exceptionality gives $\operatorname{mult}_xG<3m$. For Cartier multiples $mL$, restriction to $\cO_{X,x}/\mathfrak m_x^{3m}$ is therefore injective on sections after a local trivialization. Since this ring has length $\binom{3m+2}{3}$, comparison of leading terms gives
\begin{equation}\label{eq:exceptional-volume}
                           (-K_X)^3\le 27.
\end{equation}

\begin{lemma}\label{lem:adjoint-surfaces}
Fix a rational DCC set $J\subset[0,1]$. There is an integer $m=m(J)>0$ with the following property. Let $(S,\Theta)$ be a projective lc surface pair over an algebraically closed field of characteristic $p>5$, with $\coeff\Theta\subset J$ and $D=K_S+\Theta$ ample. Then $D^2\ge m^{-2}$, and $\varepsilon(D;x)\ge m^{-1}$ on a dense open subset of the smooth locus.
\end{lemma}
\begin{proof}
By \cite[Theorem~1.2]{DasSurf}, $|\floor{mD}|$ is birational for a single $m=m(J)$. On a resolution $\pi$ of the system, write $m\pi^*D\simq H+F$ with $H$ its free part and $F\ge0$. Indeed, the divisorial inequality for a section of $\floor{mD}$ also holds for $mD$, so its rational Cartier pullback dominates the corresponding moving divisor. Since $H$ defines a birational morphism onto its image, $m^2D^2\ge\vol(H)=H^2\ge1$. Choose $x$ where this map and $\pi$ are isomorphisms onto open subsets and $F$ is absent. For each curve through $x$, choose a hyperplane through its image point that does not contain the image curve. Local intersection bounds the multiplicity at $x$ by the degree against $H$, hence by the degree against $mD$.
\end{proof}

\begin{lemma}\label{lem:floor}
Fix $e>0$. Let $(Y,D)$ be a projective lc threefold pair in characteristic $p>5$, with $D$ reduced and $A=K_Y+D$ ample. Suppose $A\simq E$, where $\Supp E=D$ and $E\ge eD$. There is an integer $b=b(e)>0$ for which $|bA|$ is birational. The same integer works on a smooth projective threefold $W$ with reduced SNC boundary $H$ if $K_W+H\simq E_W$ for a big divisor $E_W\ge eH$ with $\Supp E_W=H$.
\end{lemma}
\begin{proof}
Since $A$ is ample, $D=\Supp E$ is nonempty. Choose a component $S$ and let $S^\nu$ be its normalization. Adjunction gives $A|_{S^\nu}\simq K_{S^\nu}+\Theta$ with standard coefficients \cite[Propositions~4.1--4.2]{Bpos}. Restricting the crepant boundary on a log resolution to the strict transform of $S$ gives a log smooth subboundary with coefficients at most one. Thus $(S^\nu,\Theta)$ is lc. Apply \cref{lem:adjoint-surfaces} with $J=\{1\}\cup\{1-1/j:j\ge1\}$ and choose a point $x$ in its open set where $Y$ and $S$ are smooth and the other components of $D$ are absent. A curve $C$ through $x$ contained in $S$ has $A\cdot C\ge\operatorname{mult}_x C/m$. For any other such curve $C$, the effective rational Cartier divisor $E$ restricts effectively to its normalization. Near $x$ we have $E=c_S S$ with $c_S\ge e$, so $A\cdot C=E\cdot C\ge e\operatorname{mult}_x C$.

Put $c=\min\{e,1/m\}$. The strict bound $\varepsilon(A;x)>c/2$ persists on an open subset, so it also holds at a smooth point off $D$. For $b=2+\ceil{12/c}$ the integral Weil divisor $M=bA$ satisfies $M-(K_Y+D)=(b-1)A$ with Seshadri constant greater than six there. Apply \cref{lem:trace}.

For $(W,H)$ take its log canonical model $(Y,D)$ by \cite[Theorems~1.1, 1.6 and Corollary~1.8]{Wlc}. The transforms of $H,E_W$ satisfy the conditions just proved. On a common resolution the pullback of $K_W+H$ equals the pullback of $A$ plus an effective divisor exceptional over $Y$. Pulling back the effective rational Cartier divisor of a section of $bA$ and adding this exceptional difference multiplied by $b$ gives a section of $b(K_W+H)$. Hence the latter system is birational.
\end{proof}

The next estimate supplies curves of controlled degree on the surfaces which occur as fibres of pencils.

\begin{lemma}\label{lem:surface}
Let $A$ be a nef and big Cartier divisor on a smooth projective surface $V$ over an algebraically closed field of characteristic $p>5$. Write $d=A^2$. If $K_V\cdot A\le Cd$ for $C\ge0$, then $|mA|$ is generically finite for $m=\ceil C+4$. When $K_V\cdot A<d$ one can take $m=4$. In either case $\varepsilon(A;x)\ge1/m$ on a dense open of $V$.
\end{lemma}
\begin{proof}
We first bound $\chi(\cO_V)$ below. If it is negative, $V$ is birationally ruled over a curve of genus $g=1-\chi(\cO_V)$ by \cite[Theorem~1.6]{GSZ}. Push $A$ down to a relatively minimal ruled surface $V_0$ and let $a$ be its degree on a ruling fibre. The blowup formula for the successive transforms $A_i$ gives $A_i=\pi_i^*A_{i-1}-a_i E_i$, where $0\le a_i\le a$ by nefness and intersection with the total ruling fibre less $E_i$. On $V_0$ the ruled surface intersection formula of \cite[Chapter~V, Section~2]{Hart} gives $K_{V_0}A_0=2a(g-1)-A_0^2/a$. Hence
\[
 K_V A=2a(g-1)-d/a+\sum_i(a_i-a_i^2/a)\ge2(g-1)-d .
\]
In particular $\chi(\cO_V)\ge -(C+1)d/2$ in either sign of $\chi$. Since $(K_V-mA)A<0$, Serre duality and surface Riemann--Roch give
\[
 h^0(V,mA)\ge (m^2-mC-C-1)d/2>md+1.
\]
In the strict case $K_VA\le d-1$ by integrality and $\chi(\cO_V)>-d$, so $h^0(V,4A)>4d+1$ as well.

If the image had dimension one, resolve the system and take its Stein factor curve. The free part is pulled back from a divisor of degree $a'$ on this curve. A general fibre cycle has positive integral degree against the pullback of $A$, so $a'\le md$ and the system has at most $md+1$ sections. This proves generic finiteness. On the open where the rational map is defined and quasi-finite, choose a hyperplane through the image of $x$ whose pullback does not contain the curve being tested. Local intersection gives $mA\cdot C'\ge\operatorname{mult}_x C'$.
\end{proof}

The volume estimate for numerically Calabi--Yau pairs rests on a comparison of lc places.

\begin{lemma}\label{lem:places}
For any rational DCC set \(J\subset[0,1]\) there is \(\eta=\eta(J)>0\), valid in every characteristic \(p>5\), with the following property. Let \(\Delta,\Phi\) be effective lc boundaries on a Fano-type threefold \(X\) with \(K_X+\Delta\simq K_X+\Phi\simq0\). If \(\Delta\) has coefficients in \(J\) and \(\Phi\ge(1-\eta)\Delta\), every lc place of \(\Phi\) is an lc place of \(\Delta\).
\end{lemma}

\begin{proof}
Apply \cref{lem:uniform-acc-consequences} to $I=J\cup\{0,1\}$ and to its set of adjunction coefficients. Choose $\eta$ below the resulting thresholds. We use the threefold assertion to preserve log canonicity after increasing coefficients, and the numerical assertion on a curve or surface to detect equality. This is the coefficient comparison used in \cite[Theorem~7.1]{DasVol}.

Suppose a new lc place occurs for \(\Phi\). Display it on a \(\Q\)-factorial crepant dlt model \(f:(Y,\Psi)\to(X,\Phi)\) extracting only lc places, by \cite[Lemma~7.7]{Bpos}. The model is of Fano type. A small positive contribution of the crepant log Fano witness from \(X\), averaged with \(\Psi\), keeps the boundary effective because every exceptional coefficient of \(\Psi\) is one. The average is klt with nef and big negative log divisor, so Kodaira perturbation gives a log Fano boundary.

Write \(\Sigma=\lfloor\Psi\rfloor\), let \(\Gamma\) be the strict transform of \(\Delta\) off \(\Sigma\), and put \(T=\Gamma+\Sigma\). Then \(\Psi\ge\Sigma+(1-\eta)\Gamma\) and
\[
 K_Y+T=f^*(K_X+\Delta)+\sum_{S\subset\Sigma}A_{X,\Delta}(S)S .
\]
The sum is effective and the new lc place contributes positively. Run a \(-(K_Y+T)\)-MMP, which exists by Fano type. On every model the pushed \(\Psi\) is crepant, by its principal relation, and dominates \((1-\eta)T\). Hence the threshold gap keeps \(T\) lc. No component of \(\Sigma\) can contract divisorially. In such a step \(K+T\) is positive on the ray, so the coefficient of the contracted prime in the crepant trace of the next boundary would be greater than one, by negativity.

At a Mori fibre space output, \(K+T\) is relatively ample, while \(K+\Sigma+(1-\eta)\Gamma\) is nonpositive on its ray by comparison with \(K+\Psi\) on the generic fibre. Interpolation gives \(K+\Sigma+(1-u)\Gamma\equiv0\) over the base for \(0<u\le\eta\). If the base has positive dimension, restrict to the actual generic fibre and apply the curve or surface comparison to the two lc boundaries \(\Sigma+(1-u)\Gamma\) and \(T\). They agree there, contradicting relative ampleness of \(K+T\). If the base is a point, adjunction to the normalization of a surviving component of \(\Sigma\) gives the same comparison on a surface. The smaller different is at least \(1-u\) times the larger different, since the part of the different obtained with \(\Gamma=0\) is effective and the added restriction is linear in the coefficient of \(\Gamma\). The smaller adjoint is numerically trivial and the larger is ample, again a contradiction.

On a nef output \(K+T\) is both anti-nef and rationally equivalent to the nonzero effective sum of the surviving components. The effective sum has positive intersection with a very ample square, contradicting anti-nefness.
\end{proof}

\begin{corollary}\label{cor:cy-volume}
Fix a rational DCC set $I\subset[0,1]$. There is a constant $V(I)$ such that, over any algebraically closed field of characteristic $p>5$, every projective klt threefold pair $(X,B)$ with $X$ of Fano type, $\coeff B\subset I$ and $K_X+B\equiv0$ satisfies $\vol(-K_X)\le V(I)$.
\end{corollary}
\begin{proof}
By \cref{lem:numerical-torsion}, $K_X+B\simq0$. Take a small $\Q$-factorialization using a klt log Fano boundary. Its crepant transform is klt with nef and big negative log divisor, so the new model is of Fano type by Kodaira perturbation. The transform of $(X,B)$ is crepant. The anticanonical Weil section spaces agree on the two models by the description of sections via rational functions and prime divisors.

On the $\Q$-factorial model fix a rational $0<h<\min\{1,\eta(I)\}$ from \cref{lem:places} and set $D=hB\simq-hK_X$. For every $m$ with $mD$ Cartier and every $G\in|mD|$, the effective rational divisor $\Phi=(1-h)B+G/m$ satisfies $K_X+\Phi\simq0$. If $(X,\Phi)$ fails to be klt, interpolate from $B$ to $\Phi$ up to the first non-klt boundary. This boundary is lc and the parameter is rational by calculation on a common log resolution. It still dominates $(1-h)B$ and has rationally trivial log divisor. By \cref{lem:places} it has no lc place, a contradiction.

At a smooth closed point $x$ off $\Supp B$ the ordinary point blowup now gives $\operatorname{mult}_xG<3m$. The restriction of sections of $mD$ to $\cO_{X,x}/\mathfrak m_x^{3m}$, using a local trivialization, is injective. This quotient has length $\binom{3m+2}{3}$. Taking the asymptotic bound in the Cartier degrees gives $\vol(D)\le27$, so $\vol(-K_X)\le27h^{-3}$.
\end{proof}

\begin{theorem}[Discrepancies of exceptional completions]\label{th:exceptional-gap}
Fix a rational DCC set $I\subset[0,1]$. There is $e(I)>0$, uniform over algebraically closed fields of characteristic $p>5$, such that the following holds. Let $(X,B)$ be a projective exceptional threefold pair with $X$ of Fano type, $\coeff B\subset I$ and $L=-(K_X+B)$ nef. Then every effective rational $P\simq L$ satisfies
\[
 A_{X,B+P}(E)\geq e(I)
\]
for every prime divisor $E$ over $X$. The same bound holds for $(X,B)$.
\end{theorem}

\begin{proof}
If there is no uniform bound, choose a sequence of completions $B+P$ and primes $E$ for which $a=A_{X,B+P}(E)$ tends to zero. Each $a$ is positive by exceptionality. Write $u=1-a$, suppressing the sequence index. Pass to a small $\Q$-factorialization, which preserves completions and their discrepancies. Take $a$ smaller than $1/2$ and both constants of \cref{lem:uniform-acc-consequences} for $J=I\cup\{0,1\}$. If $E$ is exceptional, apply \cite[Lemma~7.7]{Bpos} to the klt pair $(X,B+P)$ to obtain a projective $\Q$-factorial extraction $f:Y\to X$ of just $E$. If $E$ is on $X$, set $f=\operatorname{id}_X$. Write
\[
 K_Y+\Psi=f^*(K_X+B+P),\qquad \Psi=C+uE+R,
\]
where $C$ is the strict transform of $B$ away from $E$ and $R\geq0$ contains no $E$. The inequality $K_Y+C+uE\geq f^*(K_X+B)$ holds since the difference is supported on $E$ with coefficient $\operatorname{coeff}_E(f^*P)$. A small positive contribution of the crepant log Fano witness from $X$, averaged with $\Psi$, gives an effective klt boundary with nef and big negative log divisor. At the extracted prime effectivity follows from $u>0$. Kodaira perturbation gives Fano type of $Y$.

Set $T=C+E$ and run a $-(K_Y+T)$-MMP. This is a klt MMP after realizing a small positive multiple of its divisor as a klt log divisor on $Y$. The models are $\Q$-factorial and of Fano type, the latter by pushing forward a big klt Calabi--Yau boundary. The pushforward of $K_Y+\Psi$ stays crepant by its principal relation. Thus on each model $(1-a)T\leq\Psi$ and the threshold part of \cref{lem:uniform-acc-consequences} gives log canonicity of $T$. The prime $E$ survives. In fact, if a divisorial contraction $g$ contracts $E$, then $K+T=g^*(K'+T')+bE$ with $b<0$, since $K+T$ is positive on the contracted ray and $E$ is negative there. Its crepant coefficient over the lc target would be $1-b>1$.

If the output $W$ has $-(K_W+T_W)$ nef, choose an effective rational $Q_W\simq-(K_W+T_W)$ by semiampleness. On a common resolution with maps $r$ to $Y$ and $s$ to $W$, negativity gives
\[
 s^*(K_W+T_W+Q_W)\geq r^*(K_Y+T)\geq(f\circ r)^*(K_X+B).
\]
Write $n(K_W+T_W+Q_W)=\divi(h')$ and set $n(K_X+B^\dagger)=\divi(h')$ on $X$. The two principal log divisors have equal pullbacks, so the inequality gives $B^\dagger\geq B$. The surviving $E$ has coefficient at least one on $W$. The same equality transports its nonpositive discrepancy to $(X,B^\dagger)$. Since $B^\dagger-B\simq L$, this contradicts exceptionality.

The output is consequently a Mori fibre space $W\to Z$ with $K_W+C_W+E_W$ relatively ample. On its ray, $K_W+C_W+uE_W$ has nonpositive degree because the effective difference $\Psi_W-C_W-uE_W$ has nonnegative degree on a curve moving in the general fibre. Interpolation gives a rational $t\in[u,1)$ with $K_W+C_W+tE_W\equiv_Z0$. This pair is klt by convexity, since the pair with coefficient $u$ lies below $\Psi_W$ and the pair with coefficient one is lc. For $\dim Z>0$, restrict to the normal generic fibre over $k(Z)$ \cite[Lemmas~3.4--3.5]{DasVol}. The numerical part of \cref{lem:uniform-acc-consequences} makes the two restricted boundaries with coefficients $t$ and one equal, since the former dominates $u$ times the latter. This contradicts relative ampleness of the rounded log divisor.

For the resulting sequence of point bases, $\rho(W)=1$ and $t\to1$. Take a subsequence with $t$ strictly increasing. The coefficients of $C_W+tE_W$ then lie in a single DCC set $J'$ obtained by adjoining these values of $t$ to $I$. By \cref{cor:cy-volume}, $\vol(-K_W)\leq V(J')$. Also $E_W$ is ample and $-K_W\equiv C_W+tE_W$, so $t^3E_W^3\leq V(J')$. Adjunction to the normalization $S$ of $E_W$ gives an lc pair $(S,\Theta)$ with coefficients in $J^{\mathrm{adj}}$ and
\[
 K_S+\Theta\equiv (1-t)E_W|_S,\qquad
 (K_S+\Theta)^2=(1-t)^2E_W^3.
\]
The restriction is ample. With $m=m(J^{\mathrm{adj}})$ from \cref{lem:adjoint-surfaces}, we obtain
\[
 0<m^{-2}\leq (K_S+\Theta)^2\leq (1-t)^2t^{-3}V(J')\longrightarrow0,
\]
a contradiction. Finally $L$ is semiample, so a completion exists and deleting its effective increment cannot decrease discrepancies.
\end{proof}

For $I=\Phi(R)$ the discrepancy bound makes the coefficients of an exceptional boundary $B$ range over a finite set. Indeed, $1-r/j\le1-e(I)$ gives $j\le r/e(I)$ whenever $r>0$.

For $B\ne0$ we increase $t$ to the pseudo-effective threshold of $-K_X-tB$. At the threshold the MMP gives either a fibration to which the fibre complement theorem applies or a Fano threefold of Picard number one.

\begin{proposition}[Reduction by the pseudo-effective threshold]\label{pr:fib}
Let $(X,B)$ be an exceptional pair as in Theorem~\ref{th:intcomp} with $B\ne0$ and $L=-(K_X+B)$ big. In characteristic $p>5$ define the threshold by
\[
      \lambda=\sup\{t\ge1\mid -K_X-tB\text{ is pseudo-effective}\}
\]
on a small $\Q$-factorialization. It is rational and exceeds one. There is a birational model $T$ and a Mori fibre space $T\to Z'$ such that $K_T+\lambda B_T$ is $\Q$-linearly trivial over $Z'$ and $-(K_T+\lambda B_T)$ is nef. A complement of $(T,\lambda B_T)$ transfers to $(X,B)$ after clearing the rounding denominator.

If $\dim Z'>0$, these pairs $(X,B)$ have a common complementary index in large characteristic. If $Z'$ is a point, $T$ is a $\Q$-factorial Fano threefold of Picard number one and $(T,\lambda B_T)$ is $e(R)$-lc with $K_T+\lambda B_T\simq0$.
\end{proposition}
\begin{proof}
Work on the small $\Q$-factorialization. Bigness at $t=1$ gives $\lambda>1$, and intersection with an ample square bounds $\lambda$ above. To see rationality, choose rational numbers on either side of $\lambda$. Choose a small rational $a>0$. At each endpoint write $a(-K_X-tB)$ as a klt adjoint up to rational linear equivalence. For this, start with a log Fano boundary, reserve a small common ample part and use a general perturbation in the remaining ample class \cite[Theorem~1]{Tan}. Interpolation is linear in $t$. The pseudo-effective part of this segment is rational polyhedral by \cite[Setting~6.1 and Corollary~6.5]{Sato}.

Run the divisor MMP for $-(K_X+\lambda B)$ to a nef model $Y$ and let $f:Y\to Z$ be its semiample contraction. The divisor lies on the boundary of the big cone, so $\dim Z<3$. We have $-(K_Y+\lambda B_Y)\simq f^*A$ with $A$ ample on $Z$, taking $A=0$ if $Z$ is a point. All models are of Fano type, since the pushforward of a big effective klt Calabi--Yau boundary gives one on each model. On a common resolution with maps $u$ to $X$ and $v$ to $Y$, write
\[
   u^*(-K_X-\lambda B)=v^*(-K_Y-\lambda B_Y)+F,\qquad F\ge0 .
\]
For $P_Y\ge0$ with $P_Y\simq-(K_Y+\lambda B_Y)$, the boundary $\lambda B+u_*(F+v^*P_Y)$ completes $(X,B)$. Its log divisor and $K_Y+\lambda B_Y+P_Y$ are crepant. In fact a rational function trivializing a multiple on one model trivializes the same multiple on the other. Theorem~\ref{th:exceptional-gap} therefore gives the gap $e(R)=e(\Phi(R))$ on $Y$.

The transform of $L$ is big and is numerically $f^*A+(\lambda-1)B_Y$. Thus $B_Y$ is big over $Z$. For small rational $u_0>0$, the MMP for $K_Y+(\lambda-u_0)B_Y$ over $Z$ is negative for $-B_Y$. It ends in a Mori fibre space $T\to Z'$ over $Z$. Every step is crepant for $K_Y+\lambda B_Y$, which comes from $Z$. Furthermore $B_T$ is ample over $Z'$ and the preceding transfer of the principal log divisor still applies. For any resulting index divisible by $q$, integrality and \eqref{eq:round} give the required rounding on $X$.

Suppose $\dim Z'>0$. Write $F_\eta$ for the generic fibre over $k(Z')$. It has ample $-K_{F_\eta}$ and $K_{F_\eta}+\lambda B_{F_\eta}\equiv0$. The discrepancies restrict by localizing a birational model displaying each divisor over the generic fibre \cite[Lemmas~3.4--3.5]{DasVol}. They give an $e(R)$ lower bound also for $F_\eta$ with zero boundary. Choose $q_0$ clearing the finitely many possible coefficients of $B$ given by Theorem~\ref{th:exceptional-gap} and the hyperstandard formula. Localization makes $q_0B_{F_\eta}$ integral.

For a curve fibre $F_\eta$, regularity and $H^0(F_\eta,\cO)=k(Z')$ give $\deg K_{F_\eta}=2g-2$ over $k(Z')$. Its negative degree forces $g=0$. Hence $\lambda=2q_0/j$ for a positive integer $j\le2q_0$. For a surface fibre choose $0<e_0<e(R)$ rational. The boundedness theorem for del Pezzo surfaces over their fields of constants gives integers $r,V_2$ depending only on $e_0$ such that $H=-rK_{F_\eta}$ is ample Cartier and $K_{F_\eta}^2\le V_2$. Geometric integrality and tameness in $p>5$ are included in the application of \cite[Theorems~5.6, 5.12 and Corollary~5.8]{BM}. Intersecting with $H$ gives
\[
 \lambda=\frac{q_0H^2}{rj},\qquad
 1\le H^2\le r^2V_2,\qquad
 1\le j=(q_0 B_{F_\eta})H\le q_0rV_2 .
\]
Thus $\lambda$ takes finitely many values in either case. The coefficients of $\lambda B_T$ belong to a fixed hyperstandard set after adjoining the numbers $1-b$ for these finitely many coefficients $b$ to $R$. Apply \cite[Theorem~1.1]{Fib} to $T\to Z'$ and transfer the resulting principal relation. If $Z'$ is a point, then so is $Z$, and $-K_T\simq\lambda B_T$ is ample. This gives the stated second outcome.
\end{proof}

\subsection{Birational models in bounded families}\label{s:arithmetic}
For a projective family we need control of Cartier indices, divisors of small log discrepancy and birational contractions. On each component of the reduced closure of the selected geometric points, we spread the construction over a dense open as in Section~\ref{s:prelim}, then apply Noetherian induction to its complement. A finite surjective cover of the open reaches every geometric point there.

\begin{lemma}\label{lem:index}
Let a projective family have normal threefold fibres and separately marked prime divisors with fixed rational coefficients. Suppose the selected positive-characteristic fibres \((Y_s,D_s)\) are klt with \(K_{Y_s}+D_s\) \(\Q\)-Cartier. Above the absolute characteristic bound in \cite[Corollary~4.4]{Bern}, the family can be stratified as above so that a relative line bundle gives a fixed Cartier multiple of \(K_{Y_s}+D_s\) on each stratum. The coefficients of the crepant boundary on a relative marked log resolution are constant, including on the geometric generic fibre. If \(-(K_{Y_s}+D_s)\) is ample on the selected fibres, the same holds on the geometric generic fibre.
\end{lemma}

\begin{proof}
Take a log resolution of the geometric generic fibre and spread it to \(h:W\to Y\), with \(W\) smooth projective over the base, the exceptional primes \(E_i\) and strict marks relatively SNC, and \(Y\) flat with geometrically normal fibres. Spread the open isomorphism of \(h\). Shrinking keeps the dimensions of the components of its complement and of their images fixed, so the exceptional primes on the fibres are precisely the \(E_i\).

The classes of the \(E_i\) are independent relative to \(h\) by the negativity lemma applied to both signs of a numerically trivial exceptional divisor. Choose contracted curves on the geometric generic fibre giving an invertible matrix of intersections with the \(E_i\). Spread the curves and their images. On selected fibres the constants in
\begin{equation}\label{eq:matrix}
 K_{W_s}+\widetilde D_s=h_s^*(K_{Y_s}+D_s)+\sum_i a_i E_{i,s}
\end{equation}
are therefore the solution of a single rational linear system. Clear the denominators of the marks and the \(a_i\) by \(r\). The line bundle \(M=\cO_W(r(K_{W/S}+\widetilde D-\sum_i a_i E_i))\) is numerically trivial over \(Y_s\) on each selected fibre. Bernasconi's descent theorem makes it the pullback of a line bundle on \(Y_s\), namely \(h_{s,*}M_s\). The theorem applies over the algebraically closed field of each selected fibre, with the same characteristic bound.

After shrinking, \(h_*M\) commutes with passage to geometric fibres. Indeed, over affine opens of \(Y\) compute the pushforward by finite affine \v Cech complexes. Their terms are flat over the base, and generic flatness for their coherent cohomology sheaves gives the required base change. Fibrewise local freeness and the evaluation \(h^*h_*M\to M\) then give an invertible pushforward and an isomorphism on an open containing the generic fibre, by generic flatness, Nakayama's lemma and properness. In codimension one on each fibre the pushforward represents \(r(K_Y+D)\). The pullback identity gives the generic crepant trace, so the relative SNC calculation proves the discrepancy assertion. Ampleness of the inverse line bundle is open. Repeat on the closed remainder of the selected locus.
\end{proof}

\begin{lemma}\label{lem:toroidal}
Let \(W\) be smooth, \(\Sigma=\sum_i\Sigma_i\) reduced SNC and \(D_W=\sum_i d_i\Sigma_i\) with \(d_i\le 1-e\), \(e>0\). Any prime exceptional over \(W\) of log discrepancy at most one for \(D_W\) is obtained by blowing up strata of \(\Sigma\). The sum of its positive orders along the original components is at most \(1/e\). On a relative SNC family finitely many such modifications suffice, and on each stratum a single spread smooth model can dominate them.
\end{lemma}

\begin{proof}
For a divisorial valuation \(v\), the SNC discrepancy formula gives
\[
 A_{W,D_W}(v)=A_{W,\Sigma}(v)+\sum_i(1-d_i)v(\Sigma_i).
\]
The first term is a nonnegative integer. An exceptional centre outside the boundary has discrepancy at least two, as is seen by blowing up its smooth generic locus. If the centre is in the boundary, the sum is positive and the assumed bound forces \(A_{W,\Sigma}(v)=0\) and \(\sum_i v(\Sigma_i)\le 1/e\). The centre of a zero-discrepancy valuation is a stratum. Indeed, a centre of codimension \(r\) in only \(s<r\) components would contribute at least \((r-s)v(E)>0\) on blowing up its generic locus with exceptional divisor \(E\) and adding \(E\) to the reduced boundary. Blowing up the intersection of two components through the centre replaces the two positive orders \(w_i,w_j\) by \(\min(w_i,w_j)\) and \(|w_i-w_j|\), with zero entries omitted. This is crepant for the reduced boundary and decreases the sum of the orders. The process ends when the centre itself is a divisor. The same argument, without the order bound, describes all lc places of a reduced SNC pair as monomial valuations on its strata. For the bounded list, split the strata into geometric components and resolve the common graph on the geometric generic threefold before spreading it.
\end{proof}

\begin{lemma}\label{lem:chambers}
Let \(T\to S\) be projective over an integral stratum, whose geometric generic fibre \(T_\eta\) is a \(\Q\)-factorial threefold of Fano type. Assume the generic characteristic is zero or greater than five. After a finite cover and shrinking, all projective normal birational contractions of the geometric fibres occur in finitely many projective families. The families admit rational log Fano boundaries on every remaining geometric fibre.
\end{lemma}

\begin{proof}
\emph{Divisor classes.} Spread a smooth projective resolution \(g:W_\eta\to T_\eta\) with its exceptional primes. In characteristic zero, rationality of klt singularities and log Fano vanishing give \(H^1(W_\eta,\cO)=H^2(W_\eta,\cO)=0\). Comparison over \(\C\), after descent to a finitely generated field, and Lefschetz \((1,1)\) show that divisors span the second \(\ell\)-adic cohomology. In positive characteristic any two general points on \(W_\eta\) can be joined by a chain of rational curves. This rational chain connectedness follows from \cite[Theorem~4.8 and Proposition~4.11]{GNT}. The cycle class theorem \cite[Theorem~1.2]{GJ} gives the same spanning assertion, and the Albanese is trivial because every rational chain maps to a point on an abelian variety. Thus in both cases the first Betti number vanishes.

Choose a rational basis of \(\Cl(T_\eta)_\Q\) with Cartier representatives \(D_i\). Their pullbacks and the exceptional primes of \(g\) give a basis of \(\Pic(W_\eta)_\Q\) by pushforward and negativity. Spread these divisors. For one prime \(\ell\) invertible on the base, smooth proper base change keeps their classes a basis in degree two and keeps the first Betti number zero. The reduced Picard variety of \(W_s\) is therefore zero, so the Kummer injection makes the spread divisors a basis of \(\Pic(W_s)_\Q\). Pushing down and using the exact description of the kernel by exceptional divisors shows that the \(D_{i,s}\) form a rational basis. In particular \(T_s\) is \(\Q\)-factorial.

\smallskip\noindent\emph{The generic models.} In characteristic zero \(T_\eta\) is a Mori dream space by \cite[Corollary~1.3.2]{BCHM}. Its rational polyhedral Mori chambers express each effective rational class on a common resolution as the pullback of a semiample class on a birational contraction plus an effective exceptional divisor \cite[Propositions~2.8 and~2.13]{Okawa}. Here and below the positive and exceptional parts vary rationally linearly on each cone.

For positive characteristic use an effective klt log Fano boundary \(\Theta\) on \(T_\eta\). Represent the vertices of a rational box around zero in the coordinates \(D_i\) by divisors \(D_v\). Choose small \(t>0\) rational and an effective ample rational \(A\) such that \((T_\eta,\Theta+A)\) is klt and \(tD_v-(K_{T_\eta}+\Theta+A)\) is ample at every vertex. Semiample perturbation gives klt boundaries \(\Delta_v\ge A\) with \(K_{T_\eta}+\Delta_v\simq tD_v\). The convex hull of the \(\Delta_v\) has the common big part \(A\). The finiteness of minimal models for this polytope \cite[Theorem~1.4]{BW} applies to its pseudo-effective adjoints. Scaling a class into the box shows that the resulting finitely many maps \(T_\eta\dashrightarrow T_j\) give minimal models for all pseudo-effective rational classes. They extract no primes because the pairs are klt.

Pushing forward a big klt Calabi--Yau boundary and its principal multiple shows that each \(T_j\) is of Fano type. The cone theorem and base point free theorem make its nef cone rational polyhedral and its rational nef classes semiample.

On a common resolution with maps \(u\) to \(T_\eta\) and \(v\) to \(T_j\), impose that the transform \(D_j\) is nef and that \(u^*D-v^*D_j\) is effective and \(v\)-exceptional. These are rational linear inequalities and equalities, since the difference is computed in a finite divisor basis. The cones cover the rational pseudo-effective classes by the minimal-model discrepancy inequalities and negativity. Conversely, for any class satisfying the conditions, the exceptional difference preserves the sections in divisible degrees. This gives the same finite description by semiample models in positive characteristic.

\smallskip\noindent\emph{Passage to the fibres.} We use curves to lift big classes from special fibres to the generic fibre. For each facet of the generic pseudo-effective cone, choose a rational class in its relative interior and use one of its semiample models. The semiample morphism has positive relative dimension. Its generic fibre has an integral reduced component after finite extension of the function field of its base. In using this component we allow a nonreduced geometric generic fibre. The component dominates the actual generic fibre.

On the semiample model the loci where the inverse map fails to be an isomorphism and the images of the resolution-exceptional divisors have codimension at least two. A general fibre avoids these loci if the relative dimension is one. In higher relative dimension take general complete-intersection curves on the reduced component, avoiding their intersections with the component. Spreading these curves gives a covering family on the smooth open of \(T_\eta\), with degree zero on the chosen facet class and positive degree on an ample divisor. Its degree functional is nonnegative on effective classes since the curves cover an open. The functional thus detects that facet. In rank one take ample complete-intersection curves for the facet at zero.

Spread the test families as proper flat families of geometrically integral curves with evaluation maps. After shrinking, their evaluations cover an open on each fibre and all the degrees against the \(D_i\) and a fixed relative very ample divisor stay constant. A big class on a special fibre has strictly positive degree on each test curve. To see this, subtract a small ample class, choose an effective representative of the remainder, and take a curve in the family avoiding containment in its support. The special class therefore lifts into the interior of the generic pseudo-effective cone.

Let \(T_s\dashrightarrow X\) be a birational contraction and let \(D_s\) be the transform of a very ample Cartier divisor \(H_X\). Since all primes of \(X\) occur on \(T_s\), testing effectivity of rational sections on a common graph gives
\begin{equation}\label{eq:traceample}
 H^0(X,jH_X)=H^0(T_s,jD_s)\qquad(j\ge0).
\end{equation}
Thus \(D_s\) is big and lifts into one of the generic cones. Spread the diagrams on common resolutions, the semiample and exceptional divisors on each rational ray, and the principal functions for the finitely many ray identities. Freeness of the chosen multiples and the exceptional image dimensions persist on an open. For a positive combination of a subset of ray generators, take the Stein factor of the product of the corresponding semiample maps. Their positive combination pulls back an ample class from this factor, by finiteness over the product image and ampleness on that product. Whenever the lifted class is big the factor is birational. Spread the finitely many such birational factors and morphisms, keeping the factors normal and finite over the product images. The same positivity then holds on the fibres for every rational combination with that subset of positive coefficients.

On the common resolution, an effective exceptional divisor preserves the sections of a pulled-back Cartier divisor in divisible degrees. Together with the ray identities, \eqref{eq:traceample} recovers \(X\) as the corresponding spread ample model. Each generic birational contraction has an effective big klt Calabi--Yau boundary by pushforward of the boundary and its principal relation. Bigness follows by comparing on the common resolution with the pullback of an ample divisor on the target. If the pushed boundary is \(\Xi\simq A_0+G\), \(A_0\) ample and \(G\ge0\), then \((1-h)\Xi+hG\) is klt with ample negative log divisor for small \(h>0\). Spread a rational log Fano boundary, its log resolution and a Cartier multiple of its log divisor with ample inverse. These data give witnesses on the surviving fibres.
\end{proof}

\subsection{Klt complements from bounded birational systems}\label{s:partner}\label{s:places}\label{s:roots}
A marking is a list of prime divisors on a variety. A birational linear system maps the threefold and its uncontracted marked divisors into projective space. We bound the degrees of these images to place them in families of finite type. If \(\Gamma\) is a reduced divisor on \(X\) and \(|G|\) is birational, the bounds on \(\vol(G)\) and \(\vol(K_X+\Gamma+14G)\) in the proof of \cite[Theorem~5.3]{Sato} bound the degrees of the image and its marked divisors by \(\vol(G)\) and \(8\vol(K_X+\Gamma+14G)\), respectively. On the normalized image the marks record the transforms of \(\Gamma\) and the primes exceptional over \(X\). Nondegeneracy bounds the ambient dimension in terms of the image degree.

These estimates give finite-type arithmetic parameter spaces, also over any one fixed prime. An integral subvariety of degree \(d\) is cut out set-theoretically by forms of degree at most \(d\). To separate a point outside it, project to a space of dimension one greater than that of the subvariety, choosing the fibre through the point to miss the subvariety. The restriction is finite since the pullback of \(\cO(1)\) is ample. The equation of the image hypersurface, pulled back, separates the point. The equations for the image and its marked components therefore lie in finitely many Grassmannians.

On the geometric generic fibres take the reduced schemes and spread them. The two defining ideals have the same radical on the spread fibres by the generic power containments. Geometric integrality and normality on dense opens spread the normalizations as finite birational maps to these reductions. Resolve their geometric generic fibres, marking the divisorial components over the marked component of the images and the resolution-exceptional primes. By stratification these constructions give finitely many smooth projective families with relative SNC markings. We use this construction both over horizontal strata and over fixed primes above the bound for Cartier descent.

\begin{theorem}[Bounded volume of the reduced adjoint]\label{th:partner}
Fix \(0<e<1\) and \(U>0\). Let \(X\) be a Fano-type threefold with an effective rational \(e\)-lc boundary \(\Omega\) and \(K_X+\Omega\simq0\). With \(D_X=(\Supp\Omega)_{\rm red}\), assume
\begin{equation}\label{eq:reducedadj}
                    \vol(K_X+D_X)\le U .
\end{equation}
This holds with \(U=(\delta^{-1}-1)^3V\) when
\begin{equation}\label{eq:partner}
 \operatorname{coeff}^{>0}\Omega\ge\delta>0,\qquad
 0<\vol(-K_X)\le V.
\end{equation}
In sufficiently large characteristic the underlying varieties are arithmetically bounded and have klt complements of a common index, depending only on \(e,U\). If a rational boundary \(B\le\Omega\) with coefficients in a fixed finite set is prescribed, the complement can dominate \(B\) with index also depending on that set. The same conclusions hold over any fixed prime above the bound for Cartier descent, with the family and index allowed to depend on that prime.
\end{theorem}

\begin{proof}
Take a small \(\Q\)-factorialization. It preserves the section spaces of strict Weil transforms. The effective divisor \(D_X-\Omega\simq K_X+D_X\) lies between \(eD_X\) and \((\delta^{-1}-1)\Omega\) when \eqref{eq:partner} holds, so \eqref{eq:reducedadj} follows. On a resolution add its reduced exceptional divisor to the strict support of \(\Omega\), obtaining \(D_W\). The crepant trace \(\Omega_W\) satisfies \(K_W+D_W\simq D_W-\Omega_W\ge eD_W\), with the same support on the right. This divisor is big since its support contains the support of the big pullback of \(\Omega\). Lemma~\ref{lem:floor}, followed by pushforward, makes \(G=b(K_X+D_X)\) birational for a fixed \(b=b(e)\). The volumes of \(G\) and \(K_X+D_X+14G\) are bounded above by \((1+14b)^3U\). The marked image construction provides a bounded smooth model \(W\) on which \(\Omega_W\) is supported in the marks. Indeed an unmarked prime corresponds to a prime on \(X\) outside \(\Supp\Omega\), with crepant coefficient zero.

Any prime of \(X\) missing on \(W\) has log discrepancy at most one for \(\Omega\). Use Lemma~\ref{lem:toroidal} and enlarge the finite marked diagram so that \(W\dashrightarrow X\) is a birational contraction. Among the marks record the primes \(E_i\) contracted on \(X\) and the strict transforms \(S_j\) of \(\Supp\Omega\). Record the prescribed coefficients \(b_j\) of \(B\) on the \(S_j\), including zeros. This leaves finitely many choices of markings and coefficients.

Fix such a choice, take its closure, and let \(H=\sum_i E_i+\sum_j S_j\). Choose a relative very ample divisor \(A\) and a rational \(0<e'<e\). We select new coefficients \(\alpha\) supported in \(H\) by imposing
\[
 \alpha_{E_i}\le 1-e',\qquad 0<\alpha_{S_j}\le 1-e',
 \qquad \alpha_{S_j}\ge b_j .
\]
For \(U_W=K_{W/S}+\alpha\) impose the linear equations
\[
 U_W A^2=U_W A K_{W/S}=U_W A E_i=U_W A S_j=0
 \quad\text{for all }i,j.
\]
The intersection numbers are constant. The crepant boundary of any one selected fibre gives a rational solution, which we now use throughout the stratum. Then \(U_W A^2=U_W^2 A=0\), and the Hodge index theorem for divisors gives \(U_W\equiv0\) on each fibre \cite[Proposition~2.9]{Hu}. The coefficients of \(\alpha\) make the smooth marked pair sub-klt. On a selected fibre the difference between the pullbacks of \(K_{W_s}+\alpha_s\) and \(K_X+\alpha_X\), with \(\alpha_X\) the pushforward, is exceptional and numerically trivial over \(X\). Negativity gives equality. Thus \(\alpha_X\) is effective klt and \(K_X+\alpha_X\equiv0\). The base point free theorem on a Fano-type variety makes a numerically trivial rational Cartier divisor rationally trivial. Hence \(U_{W,s}\simq0\) on each selected fibre.

To find one index, we construct a generic contraction carrying exactly the primes of \(X\). On each selected fibre \(H_s\) supports the big effective divisor obtained by pulling \(\Omega\) to a common resolution and pushing it to \(W_s\). Thus \(K_{W_s}+H_s\simq H_s-\alpha_s\ge e'H_s\) is big with the same support on the right. Lemma~\ref{lem:floor} gives a common birational degree \(b'\). Four sections whose three ratios are algebraically independent give at least \(\binom{j+3}{3}\) sections in degree \(jb'\). Semicontinuity gives this lower bound generically for each \(j\), hence bigness on the geometric generic fibre. Since \(H_\eta-\alpha_\eta\equiv K_{W_\eta}+H_\eta\) has full support \(H_\eta\), every positive rational combination with this support is big.

Set \(\Gamma=\sum_j\alpha_{S_j}S_j+(1-e'/2)\sum_i E_i\). This is a big klt boundary on the geometric generic fibre, and \(F=\Gamma-\alpha\) is effective with support exactly the \(E_i\). On each selected fibre, \(h^0(W_s,m(K_{W_s}+\Gamma_s))\le1\) in each degree clearing the coefficients. Indeed after a further multiple the divisor is linearly equivalent to a multiple of \(F_s\). An effective divisor in this latter linear system pushes to an effective principal divisor on \(X\), so its defining function is constant. Two independent sections in degree \(m\) would have independent powers in the further degree. Semicontinuity thus bounds the generic spaces by one in each of the divisible degrees.

Run the generic MMP for \(K_W+\Gamma\), using the big-boundary minimal-model theorem \cite[Theorem~1.2 and Corollary~1.4.2]{BCHM} or \cite[Theorems~1.4--1.6]{BW} according to characteristic. Its final nef log divisor is semiample by the big-boundary base point free theorem (in positive characteristic \cite[Theorem~1.4(3)]{HNT}). The section bound forces Iitaka dimension zero, so the log divisor is rationally trivial. Throughout the MMP the pushforward of \(U_{W,\eta}\) stays numerically trivial by negativity. A negative divisorial step therefore contracts a component of the current \(F\), because the contracted curves cover the prime. On the final model \(T_\eta\) the pushforward of \(F\) is effective and numerically zero, hence zero. Exactly the \(E_i\) have been contracted. The pushforward of \(\Gamma\) is a big klt Calabi--Yau boundary, making \(T_\eta\) of Fano type.

The reduced Picard variety of \(W_\eta\) is trivial by log Fano vanishing in characteristic zero and rational chain connectedness of a smooth common resolution in positive characteristic \cite[Theorem~4.8 and Proposition~4.11]{GNT}. Consequently \(U_{W,\eta}\simq0\). Choose an integral principal multiple on the geometric generic fibre. The rational function giving that multiple descends after finite extension and spreads, after shrinking, to a trivialization on the family. On each selected \(X\) push forward this principal relation. It gives a klt complement of the stated index with the required dominance.

Spread also the generic MMP output and a common resolution with \(W\). Keeping the divisorial image dimensions fixed ensures that \(W_s\dashrightarrow T_s\) contracts exactly the \(E_i\) and extracts no prime. Then \(T_s\dashrightarrow X\) is small, and Lemma~\ref{lem:chambers} bounds \(X\), including the small contraction back from the chosen \(\Q\)-factorialization. Stratifying the remainders gives finitely many such constructions. A common multiple of their indices finishes the proof.
\end{proof}

By \cref{lem:places}, to perturb an lc complement it suffices to increase the discrepancies at its lc places. We do this with sections of a single birational linear system.

\begin{corollary}\label{cor:move}
In Lemma~\ref{lem:places} let \(X\) be a klt Fano threefold. If \(P\ge0\), \(P\simq-K_X\) and \(A_{X,P}(v)>0\) at every lc place \(v\) of \(\Delta\), then \((1-h)\Delta+hP\) is klt for any rational \(0<h<\eta\).
\end{corollary}

\begin{proof}
Choose an effective klt complement \(\Xi\) of individual index by semiample perturbation. Interpolate from the klt boundary \((1-h)\Delta+h\Xi\) to \((1-h)\Delta+hP\). Every boundary in this segment dominates \((1-h)\Delta\) and has positive discrepancy at the old lc places. If the endpoint fails to be klt, the first non-klt boundary is lc and occurs at a rational parameter by calculation on a common log resolution. Its lc place contradicts Lemma~\ref{lem:places}.
\end{proof}

\begin{proposition}\label{pr:allplaces}
Let \(\Delta\) be an lc \(N\)-complement on a klt Fano threefold in characteristic \(p>5\). Suppose for each lc place \(v\) of \(\Delta\) there is \(G_v\in|-mK_X|\) with \(v(G_v)<mA_{X,0}(v)\), for one \(m\). Then a klt complement exists of index depending only on \(m,N\).
\end{proposition}

\begin{proof}
Resolve the complement and the fractional base ideal of \(|-mK_X|\) on a common smooth model. The pullbacks of its section divisors are \(F+M_G\), with \(F\) a fixed effective rational divisor and \(M_G\) a free system. Every old lc place has centre a stratum of the coefficient-one components of the crepant SNC boundary, by the discrepancy calculation of Lemma~\ref{lem:toroidal}. Choose a member of the free system avoiding containment of all these finitely many strata. Its order at every old lc place is zero. Thus the corresponding section divisor attains the minimum order in the original system at every lc place of the original complement at once. Apply Corollary~\ref{cor:move} to \(P=G/m\) with a fixed rational \(h<\eta(\{0,1/N,\ldots,1\})\). The mixture is a klt complement at index \(mN\operatorname{den}(h)\), by multiplying and adding the two principal relations.
\end{proof}

\begin{lemma}\label{lem:lipschitz}
Let \((Z,\Sigma)\) range over a bounded smooth projective family of threefolds with reduced relative SNC boundary and a relative very ample divisor \(H\). For an effective rational divisor \(G\), let \(v_w(G)\) be its monomial order along a stratum with nonnegative weights \(w\) in the local boundary parameters. There is a family constant \(C\) such that
\[
 |v_w(G)-v_{w'}(G)|\le C(H^2G)\|w-w'\|_\infty.
\]
If \(H^2B\le C_0\), \(B\ge0\) and \((Z,B)\) is \(e\)-lc, there is \(\eta>0\) depending only on the family, \(C_0,e\), such that at each lc place \(v\) of \((Z,\Sigma)\),
\[
                        A_{Z,B}(v)\ge\eta A_{Z,0}(v).
\]
\end{lemma}

\begin{proof}
Choose \(\ell\) so that \(\ell H-\Sigma_i\) is globally generated for all marked components after stratification. For a stratum of codimension \(k\), choose a finite projection \(Z\to\PP^3\) with pullback \(\cO(\ell H)\) of \(\cO(1)\), in which the first \(k\) affine coordinates cut out the components transversely at the generic point. For the prescribed coordinates take the defining sections of the components times general sections of \(\ell H-\Sigma_i\). Their common zero locus has codimension at least \(k\) by SNC and successive general base point free cuts. The remaining \(4-k\) sections make the projection base point free, hence finite by ampleness. The unrestricted coordinate ratios can be chosen with independent differentials along the stratum by very ampleness. The projection is therefore \'etale at its generic point.

For prime \(G\), the pushed hypersurface cycle has degree \((\ell H)^2G\). Its equation pulls back to a multiple of a local equation of \(G\). In the completed local ring of the stratum take the transverse projection coordinates as parameters. The residue extension along the stratum is separable, so one may choose the coefficient field to contain the unrestricted ratios. The equation of the pushed cycle is then polynomial in the transverse variables with total degree at most \((\ell H)^2G\). For positive rationally independent weights the first monomial of each factor is unique. Their product bounds the total exponent of the first monomial of \(G\) by the same degree. It follows by continuity on the positive cone that the order is the minimum of linear forms with that exponent bound. On a face the variables of weight zero belong to the coefficient ring for the remaining parameters and give the same minimum. This proves the Lipschitz bound for primes and hence for effective rational sums.

The function \(F(w)=\sum_iw_i-v_w(B)\) has a uniform Lipschitz constant \(M\) and is homogeneous. For any nonzero nonnegative integral vector \(u\), discrepancy gives \(F(u)\ge e\), using homogeneity if \(u\) is nonprimitive. Normalize \(\|w\|_\infty=1\), and take an integer \(Q>2M/e+2\). By simultaneous approximation there are \(1\le j\le Q^k\) and an integral vector \(u\ge0\) with \(\|jw-u\|_\infty\le1/Q\). A coordinate of \(u\) is \(j\). Hence \(jF(w)\ge e/2\), giving the assertion with \(\eta=e/(6Q^3)\) since \(\sum_iw_i\le3\).
\end{proof}

\begin{lemma}\label{lem:root}
Let \(Z\) be a smooth projective threefold in characteristic \(p>0\), \(H\) very ample and \(D\) integral Cartier. Let \(V_j\subset H^0(Z,jD)\) be nonzero linear systems for multiples \(j>0\) of \(r\), with base ideals \(\mathfrak b_j\). Fix an integer \(t>0\), and for \(q=p^e\) take \(j\) least divisible by \(r\) with \(j\ge qt\). For all sufficiently large \(q\) there is an ideal \(I\) such that \(I\otimes\cO_Z(K_Z+tD+4H)\) is globally generated and, for each divisorial valuation \(v\),
\[
 \frac{v(\mathfrak b_j)}q-(1-1/q)A_{Z,0}(v)
          \le v(I)\le\frac{v(\mathfrak b_j)}q .
\]
\end{lemma}

\begin{proof}
Let \(I\omega_Z\) be the trace image of \(F^e_*(\mathfrak b_j\omega_Z)\). In a local Frobenius basis the coefficient projections are obtained by multiplying in the source and taking trace. Thus \(\mathfrak b_j\subset I^{[q]}\), yielding the upper bound. For the lower bound display \(v\) on a smooth model. A local canonical generator on \(Z\) pulls back with order \(A_{Z,0}(v)-1\). In function coordinates the trace therefore applies the regular trace on the valuation ring to a source multiplied by a factor of order \(-(q-1)(A_{Z,0}(v)-1)\). The regular trace takes an element of order \(s\) to order at least \(\lfloor s/q\rfloor\), by extracting powers of a uniformizer and applying regularity to the remainder. Using \(\lfloor s/q\rfloor\ge(s-(q-1))/q\) gives the lower bound.

By the projection formula and the evaluation of \(V_j\), the twisted trace image is a quotient of copies of
\[
             F^e_*\big(\omega_Z((qt-j)D+4qH)\big).
\]
These sheaves are \(0\)-regular for \(H\) because the \(i\)-th cohomology after twisting by \(-iH\), \(1\le i\le3\), vanishes by Serre vanishing since \(qt-j\) takes only finitely many values and \(4-i>0\). This proves global generation. The required \(q\) may depend on \(Z,D,r\), while the output twist uses the given \(t\) and \(H\).
\end{proof}

\begin{theorem}[A klt complement from bounded scaled volume]\label{th:birvol}
Fix \(\eps>0\), \(N>0\) integral and \(U>0\). Let \(X\) be an \(\eps\)-lc Fano threefold with an lc \(N\)-complement \(\Delta\). Suppose some positive integer \(b=b_X\) satisfies
\begin{equation}\label{eq:birvol}
                 |-bK_X|\ \text{birational},\qquad b^3(-K_X)^3\le U .
\end{equation}
In sufficiently large characteristic such \(X\) are arithmetically bounded and have one klt complementary index depending only on \(\eps,N,U\). The same construction gives boundedness and an index over each fixed prime above the bound for Cartier descent, allowing dependence on the prime.
\end{theorem}

\begin{proof}
Write \(L=-K_X\). A klt index suffices for the boundedness conclusion by Theorem~\ref{th:partner}, since \(L^3\le U\). We may assume \(\Delta\) has an lc place. Put \(G=bN\Delta\sim bNL\). The system \(|G|\) is birational because it contains products of sections of \(|bL|\). Moreover
\[
  \vol(G)\le N^3U,\qquad
  \vol(K_X+(\Supp\Delta)_{\rm red}+14G)\le(15N)^3U .
\]
The latter bound follows from \((\Supp\Delta)_{\rm red}\le N\Delta\), using a small \(\Q\)-factorialization to compute volumes. The construction in Subsection~\ref{s:partner} gives a bounded smooth model \(Z\) over the normalized image, with reduced SNC divisor \(\Sigma\) of marked primes. Write \(A\) for the pullback of the projective image's hyperplane and \(H\) for a relative very ample divisor on \(Z\). There is a common \(c\) with \(cA-H\) effective up to linear equivalence. Indeed the image of a general \(H_0\in|H|\) is a surface of degree at most \(A^2H\). Project the projective image finitely to \(\PP^3\) and pull back the equation of the image of that surface. The resulting member of a bounded multiple of \(|A|\) contains \(H_0\).

\emph{Bounds on the marked model.} On a common resolution with maps \(f\) to \(X\) and \(g\) to \(Z\), put \(B_0=g_*f^*\Delta\). Choose \(\divi_X(u)=N(K_X+\Delta)\), and define \(K_Z+\Delta_Z=\divi_Z(u)/N\) using compatible canonical divisors. Both \(B_0\) and \(\Delta_Z\) are supported in \(\Sigma\). Also \(\Delta_Z\le\Sigma\) by lc. Write \(M=bNf^*\Delta=g^*A_\Delta+F\), where \(A_\Delta\) is the hyperplane member corresponding to \(bN\Delta\) and \(F\) is the effective fixed part. The three divisors \(M,g^*A,g^*H\) are nef, so successively replacing classes by larger ones with effective differences gives
\[
 bN H^2B_0=(g^*H)^2M\le c^2(g^*A)^2M\le c^2 M^3\le c^2N^3U .
\]
There is a prime of coefficient one in \(\Delta_Z\), since a crepant SNC boundary with every coefficient below one is sub-klt. At that prime the coefficient in \(B_0\) is the log discrepancy for \((X,0)\), at least \(\eps\). We obtain a bound on \(b\). The negative coefficients of \(\Delta_Z\) are bounded as well, since \(H^2\Delta_Z^-=H^2(\Delta_Z^++K_Z)\) and \(\Delta_Z^+\le\Sigma\). Fix a nonnegative integer \(C\) such that \(B_0-\Delta_Z\le C\Sigma\) on the bounded families.

Choose an effective rational \(\eps/2\)-lc Calabi--Yau boundary \(\Xi_X\), by taking half of a klt representative of \(2L\) \cite[Theorem~1]{Tan}. Its crepant trace \(\Xi_Z\) is \(\Delta_Z-B_0+g_*f^*\Xi_X\). Thus \(\Xi_Z\ge-C\Sigma\), and with \(\alpha=1/(1+C)\) the boundary \(B=\alpha\Xi_Z+(1-\alpha)\Sigma\) is effective and \(\alpha\eps/2\)-lc. Its \(H^2\)-degree is bounded by its numerical class. Every old lc place \(v\) is an lc place of \((Z,\Sigma)\). Applying Lemma~\ref{lem:lipschitz} gives a uniform comparison at all of them,
\begin{equation}\label{eq:compareold}
 \alpha A_{X,0}(v)\ge \alpha A_{X,\Xi_X}(v)=A_{Z,B}(v)
                     \ge \eta A_{Z,0}(v).
\end{equation}

\emph{Transfer of sections in a fixed degree.} Set \(D=N(\Delta_Z+C\Sigma)\), an effective integral divisor with bounded coefficients, and \(T=g^*D-Nf^*\Delta\). The difference \(g^*B_0-f^*\Delta\) is effective by negativity, since it is \(g\)-exceptional and \(g\)-anti-nef. Hence \(T\ge0\) and \(v(T)\le v(D)\) for every \(v\). For a rational function \(s\) the section divisors satisfy
\begin{equation}\label{eq:transferroot}
 g^*(\divi_Z(s)+kD)=f^*(\divi_X(s)+kN\Delta)+kT .
\end{equation}
The original linear system \(V_k=H^0(X,kN\Delta)\) is consequently a subsystem of \(H^0(Z,kD)\). Choose an individual \(r\) with \(rN\Delta\) base point free Cartier. For \(k\) divisible by \(r\), the order of the base ideal \(\mathfrak b_k\) along each \(v\) is \(k v(T)\), since a section of the base point free system can avoid the centre of \(v\).

The divisor \(R=bD-A_\Delta\) is effective because \(g^*R=bT+F\). There is a common \(\ell\) and an effective \(J_0\sim\ell A-(K_Z+4H)-2\Sigma\), using the effective representative of \(cA-H\) and relative Serre vanishing for a high \(H\)-twist. With \(m_2=b\ell\) put
\[
       J=\ell R+2\Sigma+J_0\sim m_2D-(K_Z+4H).
\]
Then \(v(J)\ge m_2v(T)+2v(\Sigma)\) for every valuation. Its \(H^2\)-degree is uniformly bounded from its linear class. At the old lc places Lemma~\ref{lem:lipschitz} gives \(v(J)\le C_1 A_{Z,0}(v)\) for a uniform \(C_1\). Choose an integer \(t\) with \(Nt>(C_1+1)\alpha/\eta\).

For each old lc place \(v\), Lemma~\ref{lem:root} gives a section divisor \(S_v\in|K_Z+tD+4H|\) from the generated ideal with \(v(S_v)=v(I)\). Choose \(q\) large enough also that \((r/q)\max\operatorname{coeff}D\le1\). Write \(S_v+J=\divi_Z(s)+(t+m_2)D\) and set \(G_v=\divi_X(s)+(t+m_2)N\Delta\). To see that \(G_v\ge0\), test each prime \(P\) of \(X\) using \eqref{eq:transferroot}. If its centre lies in \(\Sigma\), then \(P(\Sigma)\ge1\) and
\[
 A_{Z,0}(P)-P(\Sigma)=A_{Z,\Sigma}(P)\le A_{X,\Delta}(P)\le1 .
\]
The lower trace bound gives \(P(S_v)\ge tP(T)-A_{Z,0}(P)\). The term \(2P(\Sigma)\) in the lower bound for \(P(J)\) cancels this possible loss, so \(P(S_v+J)\ge(t+m_2)P(T)\). A prime with centre generically off \(\Sigma\) is already on \(Z\), since an exceptional valuation there would have \(A_{Z,\Sigma}\ge2\). For such a prime both \(P(D)\) and \(P(T)\) vanish, and effectivity follows.

At the chosen old place use the upper trace bound and \(qt\le j<qt+r\) to obtain
\[
 v(G_v)\le v(J)-m_2v(T)+(r/q)v(T)
     \le(C_1+1)A_{Z,0}(v)<N(t+m_2)A_{X,0}(v).
\]
Here \(v(T)\le v(D)\le(\max\operatorname{coeff}D)v(\Sigma)\) and \(v(\Sigma)=A_{Z,0}(v)\). We have constructed an effective section divisor of degree \(N(t+m_2)\) for each old place. The degree ranges over a finite set since \(b\) is bounded. Taking powers puts all sections in one common degree and preserves the strict inequalities. Proposition~\ref{pr:allplaces} gives the klt index, and Theorem~\ref{th:partner} gives arithmetic boundedness. On vertical strata the same order and trace estimates apply to the spread smooth marked models.
\end{proof}

\section{Anticanonical volume and birationality}\label{s:degrees}

\subsection{The first anticanonical degree and the rank-one reduction}\label{s:scaled}\label{s:volume}
Let $X$ be an exceptional Fano threefold with zero boundary and set $L=-K_X$. Since $L$ is ample, the integer $d=\min\{j>0\mid H^0(X,jL)\ne0\}$ exists. In fact $d\ge2$. An effective member of $|L|$ would have a component of coefficient at least one and would give a non-klt completion.

\begin{lemma}\label{lem:window}
Let $X$ be a klt Fano threefold in characteristic $p>5$ and put $L=-K_X$. If $|bL|$ defines a generically finite map, then
\[
      H^1(X,-jL)=0\qquad\text{whenever }j\ge1
                                \text{ and }1+(p-1)j>5b.
\]
For exceptional $X$ there is an absolute $b_0$ such that $|b_0(d-1)L|$ is birational. In particular, for $p>10b_0+1$,
\[
                         \chi(X,(d-1)L)\le0 .
\]
\end{lemma}
\begin{proof}
Suppose a class in $H^1(X,-jL)$ is nonzero. Its iterated Frobenius pullbacks eventually vanish by \cite[Lemma~4.9]{XZ}, applied separately to the residues modulo a very ample Cartier multiple of $L$. Let the last nonzero class have twist $-D_0=-p^v jL$. It remains nonzero on $U=X_{\rm reg}$ since a reflexive sheaf on the normal variety $X$ has depth at least two along $X\setminus U$.

On $U$ the divisor $-D_0$ is Cartier. Write $\mathcal L=\cO_U(-D_0)$. A cocycle for the class defines an affine-line bundle with transitions $y_i=g_{ij}y_j+c_{ij}$ in local frames of $\mathcal L$. Since Frobenius kills the class, there are functions $h_i$ with $h_i-g_{ij}^p h_j=c_{ij}^p$. The equations $y_i^p=h_i$ therefore define a finite flat cover of $U$ in the bundle. A root in $k(X)$ would be regular in each normal affine chart by integrality. The roots would glue to a section of the affine-line bundle, killing the class. The cover therefore has function field of degree $p$ over $k(X)$. The relative canonical bundle of the affine-line bundle is the pullback of $\mathcal L^{-1}$, and the equations give a Cartier divisor of class $\mathcal L^{p}$. Adjunction gives the pullback of $K_U-(p-1)D_0$ as its canonical class. This is the local construction of \cite[Theorem~4.5]{XZ}.

Let $g:T\to X$ be the normalization in this field extension. Over $U$, normalization contributes the effective conductor divisor $C$ by finite duality for the hypersurface cover \cite[proof of Corollary~4.6]{XZ}. The complement of $g^{-1}(U)$ has codimension at least two, so taking closures gives a Weil divisor relation
\[
     K_T+C\simq g^*\bigl(K_X-(p-1)p^v jL\bigr).
\]
In particular $K_T+C$ is $\Q$-Cartier.

On a smooth resolution of the graph of $|bL|$ its free part $M$ is nef and big Cartier, and the pullback of $bL$ dominates $M$ up to rational linear equivalence. Take a smooth birational model $W$ of $T$ mapping to that graph resolution, with $M_W$ the pullback of $M$. Write $h:W\to T$ for the induced morphism. The divisor $K_W+5M_W$ is pseudo-effective \cite[Theorem~3.4]{KT}, so it has nonnegative intersection with the square of the nef divisor $h^*g^*L$. Add the effective remainder of $b\,h^*g^*L-M_W$ with coefficient five and the strict transform of $C$. The projection formula gives
\[
 0\le (K_T+C+5b\,g^*L)(g^*L)^2
    =\bigl(5b-1-(p-1)p^v j\bigr)(g^*L)^3 ,
\]
which is incompatible with the stated inequality.

Now let $G\in|dL|$ and assume $X$ exceptional. The pair $(X,G/d)$ is $e$-lc for one absolute $e>0$ by Theorem~\ref{th:exceptional-gap} with zero prescribed boundary. On a log resolution $f:W\to X$ take $H$ to be the union of the transformed support of $G$ and all exceptional primes. The difference between $H$ and the crepant boundary of $G/d$ is supported on $H$ and dominates $eH$. It is big because it dominates a positive multiple of $f^*G$. Lemma~\ref{lem:floor}, followed by pushforward, gives birationality of $|b_0(K_X+\Supp G)|$. Adding the fixed effective divisor $b_0(G-\Supp G)$ embeds this system in $|b_0(d-1)L|$. Apply the first assertion with $j=d$. Cohen--Macaulay duality for divisorial sheaves \cite[Theorems~3 and~19]{BK} identifies $H^2(X,(d-1)L)$ with $H^1(X,-dL)^\vee$. The former therefore vanishes, as do $H^0$ by the choice of $d$ and $H^3$ by duality and ampleness. This proves the assertion about $\chi$.
\end{proof}

\begin{proposition}\label{th:scaled}
For exceptional Fano threefolds with zero boundary, a uniform upper bound on $d^3(-K_X)^3$ gives a klt complement of uniformly bounded index and arithmetic boundedness in sufficiently large characteristic.
\end{proposition}
\begin{proof}
For $G\in|-dK_X|$ put $\Omega=G/d$. Theorem~\ref{th:exceptional-gap} gives an absolute positive lower bound for the log discrepancies of $(X,\Omega)$. Since $\Supp\Omega\le G$,
\[
       \vol(K_X+\Supp\Omega)\le(d-1)^3(-K_X)^3 .
\]
Apply Theorem~\ref{th:partner}. Thus a section whose degree $d$ varies with $X$ suffices when the volume in that degree is bounded.
\end{proof}

\begin{corollary}\label{cor:bare}
Exceptional Fano threefolds having a section in one fixed anticanonical degree form an arithmetically bounded family in sufficiently large characteristic.
\end{corollary}
\begin{proof}
If the degree is $n$, the resulting klt boundary $\Omega=G/n$ has all positive coefficients and all log discrepancies at least $1/n$, since $n(K_X+\Omega)$ is principal. The volume bound \eqref{eq:exceptional-volume} and Theorem~\ref{th:partner} apply.
\end{proof}

\begin{corollary}\label{cor:allexc}
Let $(X,B)$ satisfy the assumptions of Theorem~\ref{th:intcomp}. If $(X,B)$ is exceptional and $L=-(K_X+B)$ is not big, then it admits a klt complement of bounded index.
\end{corollary}
\begin{proof}
A suitable multiple of $L$ defines a contraction. If its image has dimension one or two, apply \cref{cor:nonbig}. If its image is a point, $K_X+B\simq0$. The coefficients of $B$ belong to a finite set by the discrepancy gap and its positive coefficients are at least $\delta_R$. On a small $\Q$-factorialization we have $\Supp B\le B/\delta_R$ and $K_X+B\simq0$. Hence $\vol(K_X+\Supp B)\le(\delta_R^{-1}-1)^3\vol(-K_X)$, and smallness preserves both volumes. The right side is bounded by \cref{cor:cy-volume}. This gives the volume hypothesis of Theorem~\ref{th:partner}. Apply it with $B$ prescribed. Any effective increment over $B$ is numerically trivial, hence zero by intersection with a very ample square. Thus a fixed multiple of $K_X+B$ is principal. Taking a common multiple with $q$ gives \eqref{eq:definition}.
\end{proof}

Large anticanonical volume gives an effective rational divisor of small anticanonical degree and high multiplicity at a smooth point. We scale its coefficient until a prescribed log discrepancy $e$ is attained and take a model carrying a prime with that discrepancy. In the resulting Mori fibre space, degree estimates on curves and surfaces force the base to be a point.

\begin{theorem}\label{th:rankone}
Fix $0<\epsilon\le1$. There is a rational $0<e<\epsilon$, depending only on $\epsilon$, with the following property in characteristic $p>5$. If the volumes of the $\epsilon$-lc Fano threefolds are unbounded, there are $\Q$-factorial $e$-lc Fano threefolds $T_i$ of Picard number one and prime divisors $D_i$ such that
\[
 -K_{T_i}\simq s_iD_i,\qquad
 s_i\longrightarrow\infty,\qquad (-K_{T_i})^3\longrightarrow\infty .
\]
Each $(T_i,D_i)$ is lc, and $T_i$ admits an effective $e$-lc rational boundary $\Omega_i\ge(1-e)D_i$ with $K_{T_i}+\Omega_i\simq0$.
\end{theorem}

\begin{proof}
We adapt the extraction in \cite[Section~9, proof of Theorem~1.6]{Bcomp}. By the ACC for log canonical thresholds, we may choose rational $e<\min(\epsilon,1/2)$ such that $(T,(1-e)D)$ lc implies $(T,D)$ lc for a prime $\Q$-Cartier divisor on a klt threefold \cite[Theorem~1.5]{Wimp}.

Set $L=-K_X$. If $0<a<1$ is rational and $a^3L^3>27$, compare $h^0(X,maL)$ with the length of the jets of order $3m$ at a smooth point. For large divisible $m$ this gives $G\ge0$, $G\simq aL$, of multiplicity greater than $3$. On a log resolution, increase $t$ from zero until $(X,tG)$ is $e$-lc with a prime $E$ of log discrepancy exactly $e$. This first parameter is rational, is less than $1$, and is computed among the components of the resolution boundary. Write $\Delta=tG$ and $\alpha=ta$.

To complete $\Delta$ with the same discrepancy bound, on a log resolution $h:W\to X$ divide the positive part of its crepant trace by $1-e$. This gives an SNC boundary $C$ with coefficients at most one. Semiample perturbation gives $H_W\simq h^*((1-\alpha)L/(1-e))$ with $C+H_W$ lc \cite[Theorem~1]{Tan}. Sections of a Cartier multiple of this class are pulled back from $X$, so $(1-e)H_W=h^*P$ for $P\ge0$ with $P\simq(1-\alpha)L$. The crepant trace of $\Delta+P$ is at most $(1-e)(C+H_W)$, which is $e$-lc by convexity with $(W,0)$.

Extract just $E$ by a $\Q$-factorial modification $\mu:Y\to X$, taking a small $\Q$-factorialization if $E$ is on $X$ \cite[Lemma~7.7]{Bpos}. Write $K_Y+\Delta_Y=\mu^*(K_X+\Delta)$ and $P_Y=\mu^*\Delta$. Then $\Delta_Y\ge0$ has coefficient $1-e$ on $E$, and
\[
 P_Y\ge(\epsilon-e)E
\]
since $\ord_E\mu^*\Delta=A_{X,0}(E)-A_{X,\Delta}(E)\ge\epsilon-e$. The crepant completion $\Omega_Y=\Delta_Y+\mu^*P$ is effective and $e$-lc. For rational $0<\tau<1-e$, run the MMP for $K_Y+\Omega_Y-\tau E\simq-\tau E$ to a Mori fibre space $T\to Z$ \cite[Theorem~1.7]{BW}. The transform of $E$ has positive degree on each contracted ray, so it survives and its final transform $D$ is relatively ample. The principal relation for $K_Y+\Omega_Y$ descends crepantly. Thus $\Omega_T\ge(1-e)D$, and the choice of $e$ gives $(T,D)$ lc.

Put $b=(1-\alpha)/\alpha$. The transforms on $T$ satisfy
\[
 -K_T\simq\Delta_T+bP_T,\qquad P_T\ge(\epsilon-e)D .
\]
If $\dim Z>0$, let $F$ be the generic fibre over $k(Z)$. Relative rank one and $\Omega_T\ge(1-e)D$ make $-K_F$ ample. Localization gives the discrepancy bound on the fibre with zero boundary as in Proposition~\ref{pr:fib}. The restriction $D_F$ is an integral effective divisor over this field. In dimension one, $\deg(-K_F)=2$, whence $b(\epsilon-e)\le2$. In dimension two, choose rational $0<e_0<e/2$. There are $r,V_2$ depending only on $e_0$ with $-rK_F$ ample Cartier and $K_F^2\le V_2$, over the same field \cite[Theorems~5.6, 5.12 and Corollary~5.8]{BM}. Intersecting with $-rK_F$ gives $b(\epsilon-e)\le rV_2$, since $(-rK_F)D_F$ is a positive integer.

Along the unbounded sequence choose $a\to0$ with $a^3L^3>27$. Then $b\to\infty$, so eventually $Z$ is a point. Rank one gives $-K_T\equiv sD$ with $s\ge b(\epsilon-e)$. The big klt boundary $\Omega_T$ makes $T$ of Fano type, so numerical proportionality lifts to rational linear equivalence by base point freeness. Finally the pushforward from $Y$ to $T$ preserves every section of a divisible multiple of $bP_Y$. Hence
\[
 (-K_T)^3\ge\vol(bP_T)\ge\vol(bP_Y)=(1-\alpha)^3L^3 .
\]
\end{proof}

\begin{lemma}\label{cor:thresholdpoint}
Let $T$ range over an arithmetically bounded family of klt Fano threefolds in sufficiently large characteristic. Suppose $-K_T\equiv tD_T$ for an effective nonzero rational divisor $D_T$, where $q_0D_T$ is integral for a fixed integer $q_0$ and $t\ge a>0$. Then $t$ ranges over a finite set depending on the family, $q_0$ and $a$.
\end{lemma}
\begin{proof}
Choose a relatively very ample divisor $H$ on each of the finitely many strata. Lemma~\ref{lem:index} gives a uniform $r$ for which $-rK_T$ is represented by a relative line bundle. Thus $h=(-rK_T)H^2$ belongs to a finite set of positive integers. Writing $j=(q_0D_T)H^2$ gives
\[
             t=\frac{q_0h}{rj},\qquad 1\le j\le\frac{q_0h}{ra}.
\]
Both integers in the fraction have only finitely many possibilities.
\end{proof}

\begin{proposition}\label{pr:noncan}
Let $X$ be an exceptional Fano threefold with zero boundary and a prime valuation $E$ of log discrepancy $a\le1-\gamma$ for fixed $\gamma>0$. Extract $E$ and apply the threshold construction of Proposition~\ref{pr:fib} to its coefficient. If the final Mori base has positive dimension, $X$ has a uniform klt complement depending on $\gamma$. The same holds at a point base when $(-K_T)^3$ is uniformly bounded on the resulting Fano threefold $T$. If $X$ is noncanonical at the generic point of a curve, one can take $\gamma=1/3$.
\end{proposition}

\begin{proof}
On the extraction $f:Y\to X$ set $b=1-a$, so $K_Y+bE=f^*K_X$ and $1>b\ge\gamma$. The effective klt weak log Fano pair makes $Y$ of Fano type. Define
\[
 v=\sup\{t\ge b:-(K_Y+tE)\ \text{is pseudo-effective}\}.
\]
Bigness and intersection with an ample square give $b<v<\infty$. Apply the rationality and two MMPs of Proposition~\ref{pr:fib}. On the first semiample model the pushforward of $-(K_Y+bE)$ is big, so the transform of $E$ is relatively big. The final Mori fibre space $T$ carries the relatively ample prime $E_T$, and $-(K_T+vE_T)$ comes from the first semiample base. An effective completion of $vE_T$ transfers by the same principal function to an effective completion of $bE$ and hence to a completion of $(X,0)$. It has the uniform discrepancy bound, giving $v<1$ independently of the denominator of $b$.

On a generic curve fibre $v=2/j$, where $j=\deg E_F$ is an integer between $1$ and $2/\gamma$. On a surface fibre take $H=-rK_F$ as in Proposition~\ref{pr:fib}. Then
\[
 v=\frac{H^2}{rj},\qquad
 1\le H^2\le r^2V_2,\qquad 1\le j=HE_F\le rV_2/\gamma.
\]
These computations over the generic field give finitely many possible coefficients. The fibre theorem supplies a principal relation of bounded order. At a point $(T,vE_T)$ is a klt Calabi--Yau pair with $v\ge\gamma$. The volume bound and Theorem~\ref{th:partner} bound $T$ arithmetically. The same ample-square argument as in Lemma~\ref{cor:thresholdpoint}, using $v\ge\gamma$, gives finitely many $v$. Prescribing $vE_T$ in Theorem~\ref{th:partner} bounds the order of $K_T+vE_T$. In either case the MMP comparison transfers the relation to $X$.

For the transverse assertion, localize at the generic curve point, or take the corresponding closed point on the generic surface of a parameter separating on the curve. The local ring and log discrepancies agree. On the minimal regular resolution of this klt surface singularity write $K_V+\sum b_iE_i=f^*K$. The coefficients are nonnegative by the minimal resolution calculation for rational klt surface singularities \cite[Proposition~2.6 and proof of Lemma~5.9]{BM}. Noncanonicity gives a nonzero sum. Negative definiteness provides an $E_i$ with $K_VE_i>0$. Write $E_i^2=-d_i j_i$, where $d_i$ is the degree of its constant field and $j_i$ is a positive integer. Adjunction gives $K_VE_i/d_i\ge j_i-2$, and this degree is a positive integer. Intersections with the other exceptional curves are nonnegative, so
\[
 b_i\ge\frac{K_VE_i}{d_i j_i}\ge
 \frac{\max(1,j_i-2)}{j_i}\ge\frac13.
\]
Taking the closure of $E_i$ over the curve gives the claimed valuation of $X$.
\end{proof}

\subsection{Volume bounds from moving divisors}
A pencil with prime members has a geometrically integral generic fibre. Connectedness on its normalized graph forces a lower bound for the fixed part along every horizontal exceptional prime. On a Mori fibre of a resolution of the pencil, that fixed part can be measured by a curve or surface degree.

\begin{proposition}\label{pr:moving}
Let $X$ be a $\Q$-factorial $\epsilon$-lc Fano threefold of Picard number one over an algebraically closed field of characteristic $p>5$, and let $S\equiv d(-K_X)$ be prime. If $h^0(X,\mathcal O_X(S))\ge2$, then $d\ge\lambda(\epsilon)>0$, uniformly in $p$. This applies to members of a nonconstant flat family of integral prime divisors over an integral curve. More precisely, fix rational $0<e<\min(1,\epsilon)$ and let $r,V_2$ be the Cartier multiple and volume bound for the $e/2$-lc del Pezzo fibres in the proof. One may take
\[
 \lambda(\epsilon)=\frac{e}{4\max\{2,rV_2\}+2e}.
\]
\end{proposition}

\begin{proof}
A smooth projective resolution $U$ of $X$ is rationally chain connected \cite[Theorems~1.5 and~4.8]{GNT}, so $\Pic^0(U)(k)=0$. Over the normalization of a parameter curve, the closure of the strict transform of the family is a Cartier divisor on the smooth product. Its fibre divisors push to the original members, as can be checked where $U\to X$ is an isomorphism in codimension one. Their line bundle classes on $U$ coincide by algebraic equivalence. Thus the prime members are Weil linearly equivalent on $X$. This argument also works without the Picard number assumption.

Put $L=-K_X$ and choose a nonconstant function $u\in H^0(X,\mathcal O_X(S))$. It has simple pole divisor $S$. The normalized graph $\pi:W\to X$ carries $f:W\to\PP^1$. Its fibres map finitely to $X$, and over infinity its unique prime fibre component is the strict transform of $S$, with multiplicity one. Thus the finite map in the Stein factorization has degree one. The simple pole also shows $u\notin k(X)^p$, so $k(X)/k(u)$ is separably generated. Consequently the generic fibre is geometrically integral. General closed fibres map finitely and birationally to prime divisors $S_t\sim S\simq dL$.

Suppose $d<1/2$. For distinct general parameters put $Q=S_{t_1}+S_{t_2}$ and $K_W+B_W^0=\pi^*(K_X+Q)$. If a $\pi$-exceptional prime $E$ dominates $\PP^1$, then $E\to\pi(E)\times\PP^1$ is finite onto its image. Dimension gives $\dim\pi(E)=1$ and the image is the product. In particular for general $x\in\pi(E)$ both graph fibres meet $E$ over $x$, and a point of the first meeting can be chosen off the other exceptional primes. Both graph fibres have coefficient one in $B_W^0$. Birational connectedness \cite[Lemma~3.4]{FW} makes the non-klt subset of the fibre over $x$ connected. Every contracted curve lies in the union of exceptional primes. Indeed, subtracting the pullback of the pushforward from a relatively ample Cartier divisor gives a nonpositive exceptional divisor by negativity, and this difference is positive on every contracted curve. The chosen point must therefore lie on a positive-dimensional non-klt component inside $E$. Varying $x$ shows that $E$ itself has coefficient at least one. Equivalently,
\begin{equation}\label{eq:pencil_horiz}
  \ord_E\pi^*Q\ge A_{X,0}(E).
\end{equation}

Write $K_W+\Gamma_W-R_W=\pi^*K_X$ with effective disjoint parts. After adding a pullback from $\PP^1$ we have $B_W\ge\Gamma_W$ by \eqref{eq:pencil_horiz}. On a log resolution $h:Z\to W$ use boundary
\[
 \Theta=\Gamma_W^{\rm str}+(1-e/2)\Exc(h)_{\mathrm{red}}.
\]
It is $e/2$-lc. Run its adjoint MMP over $W$ to a nef model $Y$ \cite[Theorems~1.6--1.7]{BW}. If $B_Y$ is the crepant trace of $B_W$, then $\Theta_Y-B_Y$ is nef over $W$ and pushes to a nonpositive divisor. Negativity gives $\Theta_Y\le B_Y$. A surviving horizontal $h$-exceptional prime therefore has order on $Q$ at least $A_{X,0}(E)-e/2$.

Let $H,Q_Y$ be the pullbacks of $L,Q$ to $Y$, and write $K_Y+\Gamma-R=-H$ with effective disjoint parts. Since $\Gamma\le\Theta_Y$, $(Y,\Gamma)$ is $e/2$-lc. For horizontal primes over $\PP^1$ the preceding estimates give
\[
 Q_Y\simq2dH,\qquad
 Q_Y\ge R^{\rm hor},\qquad Q_Y\ge(e/2)E
\]
for each horizontal prime $E$ exceptional over $X$. On the generic surface of the pencil, $H$ is nef and big and $(K_Y+\Gamma)H\le(2d-1)H^2<0$. Thus $K_Y+\Gamma$ is not relatively pseudo-effective. Its MMP over $\PP^1$ gives a Mori fibre space $\psi:T\to V$ with $1\le\dim V<3$.

Some prime $E_T$ exceptional over $X$ dominates $V$. In fact, on a common graph take the pullback of a general very ample divisor $A_V$ and push it to $X$. Its pushforward is ample by rank one. Its pullback, traced on $T$, is big and is $\psi^*A_V$ plus an effective divisor exceptional over $X$. If every prime in the latter were vertical, one could dominate the entire divisor by a pullback from $V$, contradicting bigness.

Let $G$ be the generic fibre of $\psi$ over $k(V)$. Localization preserves log discrepancies \cite[Lemma~2.1 and Corollary~2.2]{DW}, and relative rank one makes $-K_G$ ample because $\Gamma_G$ is effective. Restricting the inequalities on $Y$ and their transforms gives
\[
 -K_G\simq H_G+\Gamma_G-R_G
       \simq \frac{e(1-2d)}{4d} E_G+G_0,\qquad G_0\ge0 .
\]
Indeed, expand $H_G\simq (1/(2d)-1)Q_G+Q_G$ and use $Q_G\ge R_G$. The divisor $E_G$ is integral over $k(V)$. If $G$ is a curve, it has genus zero and $\deg(-K_G)=2$. If $G$ is a surface, apply the results over the generic field of \cite[Theorem~5.6, Corollary~5.8 and proof of Theorem~5.12]{BM} with strict gap $e/4$. They give $-rK_G$ ample Cartier and $K_G^2\le V_2$. The positive integer $(-rK_G)E_G$ then gives $rV_2\ge e(1-2d)/(4d)$. In the curve case $2$ serves in place of $rV_2$. Solving for $d$ gives the stated bound.
\end{proof}

A birational linear system supplies three coordinates for differentiating sections. Each derivative raises the degree by a fixed amount. Differentiating with respect to a parameter of the family preserves the degree.

\begin{lemma}\label{lem:jetderivatives}
Suppose $L=-K_X$ and $|bL|$ is birational. There are separating coordinates $t_1,t_2,t_3$ and an effective $D\sim(7b-1)L$, cut out by a rational section $u$, for which
\[
 u^{2a}\partial_t^{[I]}H^0(X,mL)\subseteq H^0(X,(m+(14b-2)a)L),
 \qquad |I|=a .
\]
The same estimate holds with parameters from a second factor, differentiating only in the first.
\end{lemma}

\begin{proof}
Choose $t_i=s_i/s_0$ among ratios of sections in degree $b$ whose differentials span the cotangent space of the function field. Put $H=\divi(s_0)+bL$, choose $K_X=\divi(dt_1\wedge dt_2\wedge dt_3)$, and set $J=K_X+4H$, $D=J+3H$. If $f_0$ cuts out $H$ along a prime and $f_i=f_0t_i$, then $f_0^4dt_1\wedge dt_2\wedge dt_3$ is the alternating sum of the regular forms $f_i\,df_0\wedge\cdots\widehat{df_i}\cdots\wedge df_3$. Hence $J\ge0$.

At a prime $P$, use regular \'etale coordinates $x$. Write $h=H_P$, $j=J_P$, and $k_P=(K_X)_P$. The $x$-Taylor coefficients of $t_i$ of order $c$ have order at least $-(c+1)h$. Indeed, for $h>0$, expanding the denominator loses at most $h+c$, and if $h=0$ the coefficients are integral. The inverse Jacobian has entries of order at least $-j$, since its determinant before inversion has order $k_P$ and its cofactors have order at least $-4h$.

The formal inverse change $w=x(t+y)-x(t)$ exists because the coordinates are separating. Its coefficient of degree $a\ge1$ has order at least $-(2a-1)(j+3h)$. This is the Jacobian estimate in degree one. A term of Taylor degree $c\ge2$ in the recursion for degree $a$ uses lower-order coefficients whose orders sum to $a$. Induction bounds its pole by $(2a-c)(j+3h)+(c+1)h$ before multiplication by the inverse Jacobian. Use $(c-1)(j+3h)\ge j+(c+1)h$ to conclude. The $x$-Taylor coefficient of order $c$ of a section $s$ of $mL$ has order at least $mk_P-c$, by expanding a uniformizer to its integral order times a unit. Substituting $w$ shows $\ord_P(\partial_t^{[I]}s)\ge mk_P-2aD_P$. If $D_P=0$, then $h=k_P=0$ and the expansions are integral. Multiplication by $u^{2a}$ clears these poles at every prime, proving the assertion.
\end{proof}

To remove the prime factors through the moving point, work at the generic diagonal of $X\times X$ and let $s$ be a section in the first factor, of diagonal order $M$. A prime $P$ of its divisor through the diagonal with coefficient $a_P$ can be removed by a Hasse derivative of order $a_P$ in a coordinate direction whose first derivative takes a local equation of $P$ to a unit at its generic point. In the Taylor coefficient of that order, only the product of the first derivatives of the $a_P$ local factors can survive modulo $P$. The derivative loses at most $a_P$ in diagonal order. There is such a direction among the two sets of separating coordinates because the reduced prime in the smooth product is generically smooth over $k$.

If none of the directions in the parameter factor works, all those derivations preserve the equation ideal. In the completed regular local ring along the diagonal, with the first projection as coefficient field, the leading homogeneous part of the equation then has all its partial derivatives zero. Thus its degree is at least $p$, and $a_P\le M/p$ in this case. Every prime through the diagonal dominates both factors, so the clearing multipliers are units there. After putting the sections in one degree, a general linear combination of the derivatives avoids all the finitely many primes of $s$ through the diagonal.

We will also use the following observation. If $a_P>M/2$, the equation of $P$ has order one along the diagonal. On the generic fibre over the parameter, this gives a smooth rational point of the corresponding prime. Its regular local ring contains the algebraic constants of the function field by normality, and its residue field is the base field. The constants are therefore trivial, and smoothness at the point gives separable generation. The prime is geometrically integral. After spreading, its members contain the moving point and form a nonconstant covering family.

\begin{lemma}\label{lem:removeprime}
Let $X$ be a $\Q$-factorial $\epsilon$-lc Fano threefold of Picard number one in characteristic $p>5$, and suppose $|bL|$ is birational. With $\lambda=\lambda(\epsilon)$ as in Proposition~\ref{pr:moving}, take rational $v>2/\lambda$ with $v^3<L^3$. If
\[
  \frac{6}{v}\left(1+\frac{14b-2}{\lambda}\right)<\alpha<\frac12 ,
\]
there are positive integers $n,\mu$ with $3n<\alpha\mu$ and a family of two members of $|nL|$ parametrized by a dense smooth open of $X$. Both have multiplicity at least $\mu$ at the parameter point, and they have no common prime through that point.
\end{lemma}

\begin{proof}
Counting sections and jets at the rational point of $X_{k(X)}$ given by the diagonal gives, in a large divisible degree $m$, a section of diagonal order $M>vm$. For a prime $P$ through the diagonal, choose a reduced geometric component of its generic fibre through the geometric diagonal point. Descend the component and incidence to a finite extension of the parameter field and spread to a covering family of integral primes. Its numerical coefficient along $L$ is at most $m/a_P$ since the section has coefficient at least $a_P$ on the component, even if the geometric fibre is nonreduced. Proposition~\ref{pr:moving} gives $a_P\le m/\lambda$. The preceding removals and Lemma~\ref{lem:jetderivatives} thus give the required pair in a common degree
\[
 n\le m\left(1+\frac{14b-2}{\lambda}\right).
\]
Each section has order at least $\mu=\ceil{M/2}$ because $m/\lambda<M/2$. Membership in the diagonal ideal power and the dimension bound for the base locus at the diagonal spread to a dense open.
\end{proof}

Two high-order members have a non-klt centre of codimension at least two at the moving point. Connectedness first excludes point centres, then bounds the degree of the remaining curve. Applying connectedness over the parameter surface then gives the volume bound.

\begin{lemma}\label{lem:smallcurves}
Under the hypotheses of Lemma~\ref{lem:removeprime}, there is a flat projective family of distinct geometrically integral curves $C_t$ that covers a dense subset of $X$, with $t$ on a smooth integral surface. Each curve satisfies
\[
 0<L C_t\le (1-2\alpha)^{-1}.
\]
For each closed $t$ and rational $\theta>0$ there is an effective $\Delta_t\simq\alpha L$ non-klt at the generic point of $C_t$, with all coefficients $<\theta$. After a finite purely inseparable extension of the parameter field and restriction to an open, the normalizations form a smooth projective family.
\end{lemma}

\begin{proof}
Fix one parameter $z_0$ and work outside the first member at $z_0$. For every point $x$ in this open, the general first member avoids $x$. After scaling by $3/\mu$, it is non-klt at its own parameter and its class is $cL$ with $c<\alpha$. Thus an effective $P\simq\beta L$ with $\beta<1-\alpha$ cannot have $x$ as an isolated non-klt point. Otherwise addition of the scaled divisor gives a disconnected non-klt locus for a log Fano pair, contrary to \cite[Theorem~4.12]{NT}. Nor can $x$ itself be an lc centre if $(X,P)$ is lc near $x$. To see this, take a member $H_x$ of a very ample system through $x$, avoiding the generic points of the positive-dimensional lc centres through $x$. On a resolution displaying a point-centred lc place, let $h$ be the minimum of $A_{X,0}(E)/\ord_E(H_x)$ over the displayed lc places centred at $x$. For small rational $\delta>0$, the boundary $(1-\delta)P+\delta hH_x$ has a point-centred lc place, all displayed discrepancies over the other centres through $x$ are positive, and the class is still $< (1-\alpha)L$. The SNC calculation makes $x$ an isolated non-klt point, which was excluded.

Resolve the pencil ideal at the geometric generic point of the parameter space. Descend the marked log resolution to a finite extension of the parameter field and spread to the normalization over a dense open. The crepant coefficients, orders of the ideal and list of centres through the parameter point are constant after shrinking. Choose a marked prime $E$ computing the threshold $\tau$ at that point. The smooth point blowup gives $0<\tau\le3/\mu$. Its centre $C$ has dimension at most one because no common prime passes through the point.

To represent this threshold with arbitrarily small coefficients, work in a closed fibre at $z$. Subtract the common fixed Weil divisor $Q$, which avoids $z$. On a further log resolution $\pi:W\to X$ the residual pencil pulls back to a free Cartier pencil $M$ and an effective exceptional fixed divisor $F$. If $K_W+B_W=\pi^*K_X$, then $B_W+\tau F$ is lc over a neighbourhood of $z$ and has coefficient one along the chosen prime. The positive part capped at one is an SNC lc boundary globally. Apply semiample perturbation to this boundary and $a\tau M$ for an arbitrarily large integer $a$. The effective representative $J'$ so obtained has coefficients at most one and avoids the coefficient-one prime. The pushforward $P=\pi_*(J'/a)$ satisfies
\[
 \pi^*P=J'/a+\tau F,\qquad P\simq \tau(nL-Q).
\]
These are identities of the pulled-back principal relation after clearing denominators. Thus $(X,P)$ is lc near $z$ and retains $E$ as an lc place. By rank one $P\simq\beta L$ with $\beta\le\tau n<\alpha$. Add a small-coefficient representative of $(\alpha-\beta)L$ supported away from $z$ to obtain $\Delta$. Taking $a$ large makes all its coefficients as small as desired. The point-centre exclusion shows that $C$ is a curve.

Set $H=(1-2\alpha)L$. If $HC>1$, take a general smooth point $x\in C$ in the open used for the point-centre exclusion and where $(X,\Delta)$ is lc. On a log resolution choose a coefficient-one prime $E'$ above $C$. We may choose a component of its fibre over $x$ whose generic point lies on no other marked divisor, by avoiding the finitely many images of vertical components of the pairwise intersections. For large divisible $m$, Riemann--Roch gives a section of $mH|_C$ of order $m'>m$ at $x$. Use the subspace of $H^0(X,mH)$ restricting to its span. By Serre vanishing and generation of $\mathcal I_C(mH)$ its base scheme is supported at the finitely many zeros on $C$. Resolve the subspace, writing its pullback as free part plus fixed part. The same capped-boundary perturbation gives an effective $B\simq H$ with $(X,\Delta+B)$ lc on a punctured neighbourhood of $x$ and with pullback containing the fixed part divided by $m$.

At the chosen component on $E'$ the base ideal is contained in $(u,w^{m'})$, where $u$ cuts out $E'$ and $w$ is a local parameter along the fibre on $E'$. The monomial valuation of weights $(m',1)$ has log discrepancy one for the crepant trace of $(X,\Delta)$ and order at least $m'$ on this ideal. Hence $(X,\Delta+B)$ is non-lc at $x$. The old lc places through $x$ have positive-dimensional centres by the point-centre exclusion, so they have order zero on $B$ by lc on the punctured neighbourhood. On a common resolution the first local threshold for adding $B$ is therefore in $(0,1)$ and gives a new lc place centred at $x$. Its boundary class is strictly below $(\alpha+1-2\alpha)L$, again excluded. This proves $HC\le1$.

The marked images form a projective family, with geometrically integral generic image since the marked prime is geometrically integral and scheme-theoretic image commutes with flat field extension. Shrink to obtain flatness and geometrically integral fibres. The Hilbert image has dimension two or three since the curves cover $X$ and the original point parameter has dimension three. In dimension three, a general hyperplane of a compactification meets a projective closure of a generic evaluation-fibre curve off its finite boundary. Thus a surface component still gives a covering family. Restrict to a smooth open in the actual Hilbert image, so the curves are distinct and each closed parameter retains the constructions above. To normalize them simultaneously, normalize the generic curve over the perfection of the parameter field, where the normalization is smooth. Descend the finite map and common isomorphism open to a finite purely inseparable extension and spread. The parameter change is radicial and preserves distinctness.
\end{proof}

\begin{lemma}\label{lem:graphlct}
Let $g:W\to Z$ be a dominant morphism of smooth varieties of dimension at most three over an algebraically closed field. The local log canonical threshold of $g^*H$ at a closed point $w$ is at most that of an effective rational divisor $H$ at $g(w)$. An effective rational divisor of multiplicity $<1$ at a smooth point is klt there.
\end{lemma}

\begin{proof}
Clear the coefficients to work with a Cartier ideal. The threshold of the pullback to $W\times Z$ equals the original one by the product resolution. Near $(w,g(w))$ the graph is cut out by the successive smooth equations $z_i-g^*z_i=0$ in \'etale coordinates on $Z$. Apply the restriction inequality for log canonical thresholds along this flag \cite[Corollary~4.5]{ZhuJet}. The same inequality along a general smooth curve through a point bounds the threshold below by the reciprocal multiplicity.
\end{proof}

\begin{lemma}\label{lem:curvedisjoint}
Let $X$ be a projective $\Q$-factorial $e$-lc threefold in characteristic $p>5$, where $0<e\le1$. Let $f:W\to T$ be a smooth projective family of geometrically integral curves over a smooth integral surface, with a generically finite dominant map $g:W\to X$. Each $F_t$ maps finitely and birationally onto $C_t$, and the $C_t$ are distinct. Suppose for each closed $t$ and rational $\theta>0$ there is $H_t\ge0$ with coefficients $<\theta$, non-klt at the generic point of $C_t$, and $H_tC_t\le a<e/2$. After shrinking, the curves are disjoint, lie in the smooth locus, and have integral anticanonical degree. The hypotheses are impossible if in every dense parameter open two distinct closed parameters admit such $H_s,H_t$ with $-(K_X+H_s+H_t)$ nef and big for arbitrarily small $\theta$.
\end{lemma}

\begin{proof}
Shrink $T$ to make every contracted prime of $g$ in $W$ horizontal over $T$. There are finitely many, since $g$ is generically finite. Write $d=[k(W):k(X)]$. For a fixed $H\ge0$ with coefficients $<\theta$ decompose $g^*H=H_h+H_v$ over $T$. A vertical prime is the reduced pullback of a prime curve on $T$, by smoothness and integral fibres. It maps onto a prime of $X$ and has coefficient $<d\theta$. The degree of $H_h$ on each fibre is $HC_t$, constant in the family for this fixed $H$. At a general point of a subvariety whose parameter image is positive-dimensional, the fibre is not contained in $H_h$, since a horizontal divisor cannot contain all the fibres above a curve. The multiplicity of $H_h$ there is at most its degree on the fibre. At most one vertical prime passes through such a point, and it is smooth there.

For general $t$ take $H=H_t$ and $a+d\theta<1$. The multiplicity bounds and Lemma~\ref{lem:graphlct} show that $\operatorname{Nklt}(W,g^*H_t)$ maps to finitely many parameters. The generic point of $C_t$ is smooth on $X$ as the family covers. At general closed points the threshold inequality puts every component over $C_t$ dominating $C_t$ in the non-klt locus. Each is a fibre, and distinctness forces it to be $F_t$. Since $F_t\to C_t$ is birational, a general point of a general $C_t$ has only one preimage. Choosing it in the open where the number of geometric points of a fibre of $g$ is the separable degree shows that $k(W)/k(X)$ is purely inseparable.

Choose a power $q_0$ of $p$ with $k(T)^{q_0}\subset k(X)$. Resolve the graph of the induced rational map from $X$ to a normal projective compactification of the Frobenius twist of $T$. Write $\pi:Y\to X$ for the smooth model with SNC exceptional divisor, and $h$ for the map to the compactification. The normalization in $k(W)$ is finite radicial over $Y$. It agrees over a dense base open with $W$, since over the generic point both are normal proper curves with the same function field. Denote the resulting map by $q:W\to Y^\circ$. Shrink the base so that each exceptional divisor of $\pi$ meeting $Y^\circ$ is horizontal and $\pi$ is an isomorphism at the general point of each $C_t$.

Write $K_Y+\Gamma=\pi^*K_X$. The positive part of $\Gamma$ has SNC support and coefficients at most $1-e$. For distinct $s,t$, choose $H=H_s+H_t$ with $2a+2d\theta<e$. By the same multiplicity bounds, $g^*H/e$ has multiplicity $<1$ at a general point of any subvariety with positive-dimensional parameter image. The graph inequality for $q$ gives that $(Y,\pi^*H/e)$ is klt there. For $e<1$ take the convex combination with the lc pair $(Y,\Gamma^+/(1-e))$, of weights $e,1-e$, to see that $(Y,\Gamma^++\pi^*H)$ is klt there. For $e=1$ the positive part is zero. Therefore $Z=\operatorname{Nklt}(Y,\Gamma+\pi^*H)$ has finite parameter image over the open. It contains the strict transforms of $C_s,C_t$. If they met on $X$ at $x$, birational connectedness would give a connected set $Z\cap\pi^{-1}(x)$ \cite[Lemma~3.4]{FW}. Its image under $h$ would contain two distinct isolated points, which is impossible. If $-(K_X+H)$ is nef and big, $\operatorname{Nklt}(X,H)$ is connected \cite[Theorem~4.12]{NT}, and $Z$ is connected by properness and the connected fibres of its map to that locus. Its image under $h$ gives the same contradiction.

Disjointness and fibrewise finiteness now make $g$ quasi-finite. The finite surjection $q$ then makes $\pi|_{Y^\circ}$ quasi-finite. Zariski's main theorem identifies it with an open immersion into $X$. Thus the curves lie in a smooth open and have integral anticanonical degree.
\end{proof}

\begin{theorem}[Volume in a fixed birational degree]\label{th:fixedbvolume}
Fix $\epsilon>0$ and an integer $b\ge1$. In characteristic $p>5$, the volumes $(-K_X)^3$ of $\Q$-factorial $\epsilon$-lc Fano threefolds of Picard number one with $|-bK_X|$ birational are uniformly bounded.
\end{theorem}
\begin{proof}
Take rational $0<e<\min(\epsilon,1)$ and $\alpha=e/16$. Fix an integer $v$ large enough for the inequalities in Lemma~\ref{lem:removeprime}, with $\lambda=\lambda(e)$. If $(-K_X)^3>v^3$, Lemmas~\ref{lem:removeprime} and \ref{lem:smallcurves} give a covering family with $H_t=\Delta_t$ and $H_t C_t\le\alpha/(1-2\alpha)<e/2$. On the family of normalizations the evaluation map is generically finite since source and target have dimension three. Moreover $-(K_X+H_s+H_t)\simq(1-2\alpha)(-K_X)$ is ample. This contradicts Lemma~\ref{lem:curvedisjoint}.
\end{proof}

\begin{theorem}[Complements with nonzero boundary]\label{th:nonempty}
Fix a finite rational set $R\subset[0,1]$. In sufficiently large characteristic depending only on $R$, every exceptional pair as in Theorem~\ref{th:intcomp} with $B\ne0$ admits a monotonic klt complement of a common index depending only on $R$.

An exceptional Fano threefold with zero boundary also has a uniform klt complement if it has a prime valuation of log discrepancy at most $1-\gamma$ for fixed $\gamma>0$.
\end{theorem}
\begin{proof}
Propositions~\ref{pr:fib} and \ref{pr:noncan} treat the positive-dimensional bases. At a point base they give a $\Q$-factorial rank-one Fano $T$ with an $e$-lc completion $P$ whose positive coefficients are at least a fixed $\delta>0$. Let $D=\Supp P$. On a log resolution the reduced marking of $D$ and the exceptional divisors exceeds the crepant trace by at least $e$ times the marking. This difference is big, since it dominates a positive multiple of the pullback of $P$. Lemma~\ref{lem:floor} gives a birational system $|J|=|m(K_T+D)|$ with $m=m(e)$. For $a_0>0$ and $a_0\ge m(\delta^{-1}-1)$ we have
\[
 a_0(-K_T)-J\simq(a_0+m)P-mD\ge0.
\]
At a general smooth point the free hyperplane test gives $\varepsilon(-K_T;x)\ge1/a_0$. Lemma~\ref{lem:trace} makes $|-bK_T|$ birational for $b+1>6a_0$. Theorem~\ref{th:fixedbvolume} bounds the volume of $T$.

The ample-square computation of Lemma~\ref{cor:thresholdpoint} now puts $\lambda$ in a finite set when $P=\lambda B_T$. If $P=vE_T$, the same computation uses $v\ge\gamma$ to bound the denominator. Applying Theorem~\ref{th:partner} with $P$ prescribed gives a bounded principal relation on $T$, and Propositions~\ref{pr:fib} and \ref{pr:noncan} transfer it to $X$. Corollary~\ref{cor:allexc} supplies the lower-dimensional cases.
\end{proof}

\subsection{Bounding the first birational degree}\label{s:movingvolume}
We now bound the volume as the birational degree varies. A section vanishing to high order at a general point can have a large fixed prime through that point. Moving the point then gives a family of primes of small degree. If there is no such prime, derivatives give two sections without a common divisor near the point, and their intersection yields a curve to which the preceding degree estimates apply. Write $R$ for the polarization and set $V=R^3$. For rational $v\ge16$, let $l$ be the smallest positive integer greater than $vV^{-1/3}$.

\begin{theorem}\label{th:jetsurface}
Consider either of the following classes of projective $\Q$-factorial threefolds over algebraically closed fields.
\begin{enumerate}[label=\textup{(\roman*)}]
\item $(Y,B)$ is dlt, $N(K_Y+B)\sim0$ for fixed $N$, and $R$ is a nef, big, semiample rational divisor with $R\ge B$ and $\Supp R=\Supp B$. Let $d$ be the first birational degree for the systems $|\floor{dR}|$.
\item $Y$ is $e$-lc for fixed $e>0$, $R=-K_Y$ is nef, big and semiample, $d\ge2$ is the first effective degree, and $|bR|$ is birational for some $b\le Jd$, with $J$ fixed.
\end{enumerate}
There are constants $P,K$ depending only on $N$ in \textup{(i)} and on $J,e$ in \textup{(ii)} such that, for every $v\ge16$, if $p>P$ and $d>Kl$, there is a nonconstant family of integral prime divisors covering $Y$ whose general member satisfies $R^2 S/V<2l/v$.
\end{theorem}

\begin{proof}
Throughout we may increase $P,K$ by constants depending only on the indicated data. In (i) set $b=d$. Take the rational point $\xi$ of $Y_{k(Y)}$ given by the diagonal, and use a bar for geometric base change. There is an integral curve $C$ through $\bar\xi$ with
\[
 0<\gamma=\frac{RC}{\operatorname{mult}_{\bar\xi}C}<D_1/d,
 \qquad b\gamma<D_2 ,
\]
for fixed $D_1,D_2$. In (ii) a larger Seshadri constant would give a section in degree $d-1$ by Lemma~\ref{lem:trace}, so one can take $D_1=4$. In (i) put $j=\floor{(d-1)/N}-N$. A Seshadri constant $>6/j$ would give a birational system in $K_Y+\ceil{jR}$ by trace with rounding boundary. Its $N$-fold products preserve the ratios and, by the principal relation for $N(K_Y+B)$ and
\[
 N\ceil{jR}-NB\le N(j+N)R ,
\]
give a birational system below $d$. The inequality uses that each positive coefficient of $R$ is at least $1/N$. For large $d$, $6/j<20N/d$. In each use of trace with a nef and big divisor, write it as ample plus effective, and add a small multiple of the effective part to the trace boundary. This gives an ample positive class and retains the strict Seshadri bound at the diagonal point. Sections and birationality descend under the ground-field extension. The degree on any curve through $\bar\xi$ is positive since the semiample map is an isomorphism on a dense open.

\emph{Sections with high diagonal order.} Write $\Lambda_j=H^0(Y,\floor{jR})$. After multiplication by a fixed rational function to a suitable power, a Hasse derivative of order $a$ in the first factor takes $\Lambda_j$ into $\Lambda_{j+D_0a}$, with $D_0\le c_0d$. The multiplier is a unit along subvarieties dominating the first factor. In (ii) this is Lemma~\ref{lem:jetderivatives}, with $D_0=14b-2$. In (i) take the separating coordinates $t_i\in\Lambda_b$, using the section $1$. Set $H=\floor{(b+N)R}$, $K_Y=\divi(dt_1\wedge dt_2\wedge dt_3)$ and $D'=K_Y+7H\ge0$. For $s\in\Lambda_j$, the Taylor calculation of Lemma~\ref{lem:jetderivatives}, starting with $\ord_P(s)\ge-\floor{jR_P}$, gives
\[
 \ord_P(\partial_t^{[I]}s)\ge-\floor{jR_P}-2aD'_P,\qquad a=|I|.
\]
If $D'_P=0$ then $H_P=R_P=0$, so the expansions are integral there. Replace $2aD'$ by $2aND'$ and multiply by $u^{2a}$ with $\divi(u)=N(K_Y+B)$ to obtain $D_0=14N(b+N)$. Derivatives in the parameter factor cost no degree.

Since $l^3V>v^3$, jet counting yields a section of $nlR$ of diagonal order $M>vn$ for large divisible $n$. If a prime through the diagonal has coefficient $>M/2$, it is geometrically integral over the parameter field by the smooth rational point argument preceding Lemma~\ref{lem:removeprime}. Spreading it gives $R^2S/V\le nl/\operatorname{coeff}(S)<2l/v$, as desired. Otherwise use the prime removals preceding that lemma. Every coefficient is at most $M/2$, and those needing first-factor differentiation are at most $M/p$. We obtain a pair in a common degree $n'\le nl+D_0M/p$ with diagonal order at least $M/2$ and no common prime there. Let
\[
 h=\ceil{l+16D_0/p}.
\]
Taking $h$th powers gives two members in degree $n'hR$ of order at least $4n'$, with the same common-prime exclusion. Pass to a multiple $n_1$ of $n'$ clearing a Cartier degree for $hR$.

\emph{Recovering the degree of the curve.} Let $V_c(x)$ be the space of sections of $n_1hR$ of order at least $\ceil{c n_1}$ at $x$. Let $\Upsilon$ be the reduced image of $C$ over the parameter field and $\mathcal U\subset Y\times Y$ its closure, with first coordinate the curve point. Its generic point $\zeta$ dominates both factors. Choose a fixed rational $\tau>0$ with $\tau(1+D_0\gamma)<1$. This is possible uniformly since $D_0\gamma<c_0D_1$. We have
\[
 h\gamma\le(l+1)D_1/d+16c_0D_1/p<e_0\tau/8 ,
\]
by increasing $K,P$, where $e_0=1/N$ in (i) and $e_0=\min(e,1)$ in (ii). A member of $V_\tau(\xi)$ must contain $\Upsilon$, since a proper intersection with $C$ at $\bar\xi$ would give $\tau\le h\gamma$.

Every member $s$ of $V_4(\xi)$ has order at least $\ceil{\tau n_1}$ at $\zeta$. If its order there were $a<\tau n_1$, a joint Hasse derivative of order $a$ would be nonzero modulo the ideal of $\mathcal U$. Indeed the regular local ring at $\zeta$ has dimension two and its conormal space injects into the cotangent space of the smooth product at the generic point of $\mathcal U$. Substituting the independent linear Taylor terms of two parameters into the initial form of $s$ gives a nonzero coefficient. The corrected derivative is nonzero on $C$, has degree at most $n_1h+D_0a$ and diagonal order at least $4n_1-a$. This contradicts
\[
 4n_1-a\le(n_1h+D_0a)\gamma
\]
by the choices of $h\gamma$ and $\tau$.

To remove the multiplicity from $\gamma$, represent $\zeta$ as $(y,x)$ over an algebraically closed extension of large transcendence degree. Each geometric component of $\Upsilon_x$ contains $x$, by conjugacy and the rational diagonal point. Only one component passes through $y$, since $y$ is transcendental of degree one over $k(x)$ and distinct components meet at finitely many points. It is also the unique component of the base locus of $V_4(x)$ through $y$, since the base locus contains $\Upsilon_x$ and has no divisorial component through $x$, by the constructed pair. These assertions at the generic point persist under base change. The order bound at $\zeta$ gives $V_4(x)\subseteq V_\tau(y)$, by specializing the generic jet kernel and restricting the ideal power to the fibre at $x$. Each component of $\Upsilon_y$ lies in this base locus and contains $y$. Thus all of them equal the unique component through $y$. Since $y$ also represents the generic point of $Y$, the reduced curve is geometrically irreducible and $(\Upsilon_x)_{\rm red}=(\Upsilon_y)_{\rm red}$. The singular points of the curve on the right are algebraic over $k(y)$. The point $x$ has transcendence degree one over $k(y)$, so it is smooth on the curve. Consequently $\operatorname{mult}_{\bar\xi}C=1$ and $RC=\gamma$.

\emph{The contradiction.} Descend the reduced geometric curve to a finite extension of the parameter field and spread it with the jet sections. The local ideal membership at $\zeta$ spreads, so the sections have order at least $\ceil{\tau n_1}$ at the generic curve point of a general closed fibre. The common-prime exclusion there also spreads, by including a power of the curve ideal in the two-generated ideal at $\zeta$. The ordinary curve blowup bounds the threshold of this pencil by $2/(\tau n_1)$. Apply the residual-pencil perturbation from Lemma~\ref{lem:smallcurves} at the generic point of the curve. The common fixed divisor misses the generic point, so subtracting it lowers the degree on the curve. This gives $H_t$ with arbitrarily small coefficients, non-klt at the generic point of $C_t$, and
\[
 H_tC_t\le 2h\gamma/\tau<e_0/2.
\]
As in Lemma~\ref{lem:smallcurves}, take a covering parameter surface in the Hilbert image and normalize after a purely inseparable parameter extension. The degree remains $RC_t=\gamma$. In (i) the pair $(Y,0)$ is $1/N$-lc, since positive completed log discrepancies are at least $1/N$ and lc places lie over the generically smooth dlt strata. Lemma~\ref{lem:curvedisjoint} directly gives $-K_Y C_t=B C_t\in\Z$. This is positive because $B$ is big, and at most $RC_t$. In (ii) that lemma makes $RC_t$ itself integral. Either contradicts $0<\gamma<1$ for large $d$.
\end{proof}

\begin{proposition}\label{pr:prepare}
Let $X$ be an $\epsilon$-lc Fano threefold with an lc $N$-complement $\Delta$. On a $\Q$-factorial dlt modification $f:Y\to X$ extracting only lc places, set
\[
 K_Y+B=f^*(K_X+\Delta),\qquad
 R=q f^*\Delta,\qquad q=\ceil{1/\epsilon}.
\]
Then $Y$ is of Fano type and the data satisfy Theorem~\ref{th:jetsurface}\textup{(i)}. Birationality of $|\floor{mR}|$ implies birationality of $|-NmqK_X|$.
\end{proposition}
\begin{proof}
An extracted prime $E$ has coefficient one in $B$ and $\ord_E(f^*\Delta)=A_{X,0}(E)\ge\epsilon$. Thus $R\ge B$ and their supports agree. Mixing $B$ with a small proportion of the crepant boundary of $(X,0)$ on $Y$ gives an effective klt weak log Fano pair, and Kodaira perturbation gives Fano type. For a rational function, the effective-divisor inequality in degree $m$ of $R$ holds exactly when it holds in degree $mq$ of $\Delta$, by pushing down and pulling back the rational Cartier inequality. Taking $N$-fold products retains the ratios, and $N\Delta\sim-NK_X$ gives the stated degree.
\end{proof}

\begin{theorem}[Uniform anticanonical volume]\label{th:fullupper}
There is an absolute $P_0$ such that for each $\epsilon>0$ all $\epsilon$-lc Fano threefolds in characteristic $p>P_0$ satisfy $(-K_X)^3\le U(\epsilon)$.
\end{theorem}
\begin{proof}
Take $P_0$ above $5$, the threshold in \cref{th:nonexceptional} for zero boundary, and the characteristic bound of Theorem~\ref{th:jetsurface}\textup{(i)} for that index $N$. If the volumes are unbounded, Theorem~\ref{th:rankone} gives $\Q$-factorial $e$-lc Fano threefolds $T$ of Picard number one with unbounded volumes and $-K_T\simq sD$ for $s\to\infty$. They lie in the same range of characteristics, which can include a fixed characteristic. The effective divisor $sD$ shows that $(T,0)$ is non-exceptional, so \cref{th:nonexceptional} gives lc $N$-complements on them.

Apply Proposition~\ref{pr:prepare} with $q=\ceil{1/e}$. Fix an integer $v\ge16$ with $2q/v<\lambda(e)$ of Proposition~\ref{pr:moving}. Since $R^3=q^3(-K_T)^3\to\infty$, eventually $l=1$. If the first birational degree exceeds $K$ in Theorem~\ref{th:jetsurface}, the resulting prime $S$ pushes forward to a moving prime of numerical coefficient
\[
 \frac{(-K_T)^2 f_*S}{(-K_T)^3}=q\frac{R^2 S}{R^3}<2q/v .
\]
A general member avoids the exceptional primes, and distinct general members have distinct images. Proposition~\ref{pr:moving} excludes this degree. Thus the prepared models have bounded birational degree. By Proposition~\ref{pr:prepare} and taking a common multiple, so do the anticanonical systems on the $T$. Theorem~\ref{th:fixedbvolume} contradicts their unbounded volumes.
\end{proof}

The degree of a pencil member controls the mixed intersection numbers on its generic surface. We seek a generically finite system on that surface whose degree against the restricted polarization is bounded in the same scale. Hyperplanes in the system bound the degrees of curves inside the fibres. The pencil bounds the degrees of transverse curves.

\begin{theorem}\label{th:smallpencil}
Let $k$ be algebraically closed of characteristic $p>\max(5,N)$, let $(X,B)$ be a projective $\Q$-factorial dlt threefold, and suppose $N(K_X+B)\sim0$. Let $R$ be a nef, big and semiample $\Q$-divisor with $R\ge B$ and $\Supp R=\Supp B$. Suppose $S_0,S_\infty$ are distinct linearly equivalent primes with $R^2 S_\infty\le\sigma R^3$. Then $|\floor{mR}|$ is birational for some $m\le C_N(1+\sigma)$.
\end{theorem}

\begin{proof}
\emph{Sections from intersection bounds.} Let $V$ be a normal projective surface over an infinite $F$-finite field $K$. At a regular closed point $x$ the denominator in the Seshadri ratio is $[\kappa(x):K]\operatorname{mult}_x C$. Suppose $M$ is integral, $\Psi\ge0$ is rational, and $D=M-(K_V+\Psi)$ is nef, big and semiample with $\varepsilon(D;x)>3$. If the boundary at $x$ is supported on at most one regular curve of coefficient $c\le1$, then $|M|$ is generically finite. The blowup and ample-model arguments of \cref{lem:trace} give jets modulo $\mathfrak m_x^{\ceil{q\beta}}$ for $3<\beta<\varepsilon(D;x)$.

To compute their trace over $\kappa(x)$, put $R=\cO_{V,x}$ and choose parameters $z_1,z_2$, taking $z_1=0$ as the boundary equation when $c>0$. Finite duality identifies the unrestricted trace as a generator of $\operatorname{Hom}_R(F^e_*R,R)$ over the source ring \cite[\S2]{SchwedeAdj}. Kunz's theorem makes $F^e_*R$ free \cite[Tag~0EC0]{Stacks}. Since $w=(z_1z_2)^{q-1}\notin\mathfrak m_x^{[q]}$, some functional takes $F^e_*w$ to a unit. Thus the trace takes $aw$ to a unit for some $a\in R$. This element belongs to the restricted source because $w$ is divisible by $z_1^{\ceil{(q-1)c}}$. Prescribe its jet and those of $z_i^qaw$. Errors in $\mathfrak m_x^{3q-1}$ trace into $\mathfrak m_x^2$ by the two-variable monomial calculation. The resulting section ratios generate $\mathfrak m_x$ by Nakayama's lemma, so $x$ is isolated in its fibre. This argument also applies when $\kappa(x)/K$ is inseparable.

If $P$ is merely nef and big, write $P\simq H+E$ with $H$ ample and $E\ge0$. Off $E$, a Seshadri bound $r$ for $P$ gives a bound $r/2$ for $P-E/2$. Apply trace to $M=K_V+\ceil{aP}$ and $\Psi=\ceil{aP}-aP+aE/2$ with $a>6/r$. For a geometrically integral regular surface this calculation can be made at the rational diagonal over $K(V)$, a regular extension of $K$. The resulting generically finite system descends by base change of sections.

Resolve the pencil on $h:W\to X$, with $f:W\to\PP^1$, generic fibre $F$, and $P=h^*R|_F$. Write $\widetilde R=h^*R$, $J=f^*\mathcal O(1)$ as a divisor. Suppose a model of $F$ has a free system $A$ defining a generically finite map, with $PA\le cP^2$. Then $tP-A$ is big for an integer $t>2c$, and $n\widetilde R-J$ is big for $n>3\sigma$, by the nef Morse inequality \cite[Proposition~5.4]{Cut}. Spread the fibre system and an effective representative of $tP-A$ over a base open. At a general closed point away from these effective differences, free hyperplanes give degree-to-multiplicity ratio at least one for $t\widetilde R$ on vertical curves and for $n\widetilde R$ on transverse curves. Therefore $\varepsilon(R;x)\ge1/(t+n)$ on a dense open of $X$. Trace gives a birational system in $K_X+\ceil{7(t+n)R}$. Taking $N$-fold products and using $N(K_X+B)\sim0$ puts the ratios in degree $N(7(t+n)+N)$ of $R$, since each positive coefficient of $R$ is at least $1/N$. It suffices to construct $A$ with $c=O_N(1+\sigma)$.

\emph{Adjunction on the generic surface.} Let $W_0$ be the normalized graph and choose $W$ so that $\Sigma=\Exc(h)_{\mathrm{red}}+\Supp B^{\rm str}$ is SNC. The graph argument of Proposition~\ref{pr:moving} applies to this pencil. Its function has a simple pole along $S_\infty$, and the fibre at infinity on $W_0$ has one reduced prime because it is finite over $X$. The Stein factor is trivial and the function is separating. Thus $F$ is regular and geometrically integral over $K=k(\PP^1)$. Every horizontal prime $E$ of $W_0$ exceptional over $X$ satisfies
\begin{equation}\label{eq:pfhor}
 2\ord_E S_\infty\ge A_{X,0}(E)\ge 1/N .
\end{equation}
Indeed the first inequality is \eqref{eq:pencil_horiz}, applied to two general pencil members whose pullbacks have the same horizontal part. Its connectedness proof uses only $\Q$-factoriality of $X$ and the prime pencil. For the second inequality use $N(K_X+B)\sim0$ and the SNC locus at generic points of dlt strata.

Put $\Theta=\Sigma|_F$, write $K_F+B_F=h^*(K_X+B)|_F$, and set $L_F=h^*S_\infty|_F$. The last divisor is the horizontal fixed part on $F$. The projection formula and Hodge index give
\[
 P^2=R^2S_\infty>0,\qquad PL_F=RS_\infty^2\le\sigma P^2.
\]
The components of $\Theta$ with positive $P$-degree arise from $B$, with coefficients in $P$ at least $1/N$, or from horizontal exceptional primes on $W_0$, which satisfy \eqref{eq:pfhor}. Thus $P\Theta\le N(1+2\sigma)P^2$. On the generic fibre of $W_0$ the trace of $B_F+2L_F$ is effective by \eqref{eq:pfhor}. On its minimal regular resolution the crepant trace stays effective, by relative nefness of the canonical divisor and negative definiteness. Since $F$ dominates that resolution, $PK_F\le2PL_F$. In particular
\begin{equation}\label{eq:pfbudget}
 N(K_F+B_F)\sim0,\qquad E=\Theta-B_F\ge0,\qquad
 PE\le\lambda P^2,\qquad \lambda=N(1+2\sigma)+2\sigma .
\end{equation}
Positive coefficients of $E$ are at least $1/N$. A coefficient-one component of $B_F$ is an old lc place of the dlt pair. Over its generic centre the log discrepancy for the zero boundary is an integer and equals the order of $B$, so $P$ has coefficient at least one there.

Run the $(K_W+\Sigma)$-MMP over $\PP^1$ \cite[Corollary~1.8]{Wlc}. The adjoint is rationally equivalent to an effective divisor by \eqref{eq:pfbudget} and its counterpart on $W$. On the generic fibre the flips are isomorphisms, and we obtain a morphism $F\to T$ to a dlt surface with $D_T=K_T+\Theta_T$ semiample by \cite[Theorem~4.10]{Tsurf}. Pullbacks of the free systems from its semiample model have mixed degree bounded by the corresponding multiple of $PE$, by MMP negativity.

Adjunction reduces the positive-dimensional cases to a curve calculation. If a regular projective curve $C$ over $K$ has $a=\deg(K_C+\Delta)>0$ with standard coefficients in $\Delta$, then $a\ge c[H^0(C,\mathcal O_C):K]$ for an absolute $c>0$. In genus zero this follows from the DCC for sums of at most six nonzero standard coefficients, counted with their degrees over the constant field. In higher genus it follows from the canonical degree and the fact that positive coefficients are at least $1/2$. One can choose a closed point of degree $O(a)$ off a finite subset. Use a member of $|-K_C|$ in genus zero, $|3K_C|$ in genus at least two, or three times a point of $\Delta$ in genus one. These systems are free by Riemann--Roch. Adjunction to the normalization of a reduced component of an lc reduced boundary gives standard coefficients \cite[Proposition~2.7]{Wimp}.

If $\kappa(D_T)=2$, pass to the ample log canonical model. Take a component $C$ of $E_T$, with coefficient at least $1/N$. Adjunction gives the positive degree $a=D_TC$ on its normalization. Choose a regular point of $C$ of degree $O(a)$ at which it is the sole boundary component. Test curves other than $C$ against $E_T$, and $C$ against its degree $a$. This gives a Seshadri lower bound for $D_T$ depending only on $N$. Surface trace for $M=mD_T$ and $\Psi=\Theta_T$ supplies a generically finite system for fixed $m=O(N)$, of mixed degree at most $m\lambda P^2$ on a resolution over $F$.

If $\kappa(D_T)=1$, let $g:T\to C$ be the semiample contraction and $D_T\simq g^*D_C$. The divisor $E_T$ vanishes on the generic curve $G$, so $\Theta_G=B_G$. The polarization on $G$ is positive, supported on this boundary, and has coefficients at least one there. Consequently $G$ is a regular genus-zero conic, smooth for $p>2$ \cite[Lemma~2.17]{BT}, and $P_G$ has degree at least two. A horizontal boundary component has normalization $S$ finite separable of degree $\ell\le2$ over $C$. Adjunction gives $\deg(K_S+\mathrm{Diff})=\ell\deg D_C$, so $\deg D_C\ge c/2$ by the curve calculation. Riemann--Hurwitz gives $\deg K_C\le\deg D_C$. For a free $H_C$ of degree $O(\deg D_C)$, use $3K_C$ in genus at least two or $-K_C$ in genus zero. In genus one use three times the pushforward of a canonical member on $S$ when $\deg K_S>0$, or of a point from the different of degree at most $2\ell\deg D_C$ when $\deg K_S=0$. The constant field of $C$ is $K$ because $F$ is geometrically integral.

The pullback $J_C$ of $H_C$ on $F$ has $PJ_C=O(\lambda)P^2$. At the generic diagonal, a vertical curve is the smooth conic, while transverse curves are tested by $J_C$. Hence $\varepsilon(P+J_C;x)\ge1$. The ample perturbation trace for $P+J_C$ with fixed $a>6$ gives a generically finite system in $K_F+\ceil{a(P+J_C)}$. Its mixed degree is $O(\lambda)P^2$ by \eqref{eq:pfbudget}, the bound on $PK_F$, and the rounding error supported on $\Theta+\Supp J_C$.

\emph{Iitaka dimension zero.} Here $E_T=0$, and its support on $F$ consists exactly of the contracted curves. Indeed a curve contracted at any step must be negative against the current effective difference $\Theta-B$, and every component of $E$ vanishes on $T$. Thus $\Theta_T=B_T$ is reduced dlt and $N(K_T+\Theta_T)\sim0$. The pushforward $P_T$ is nef and big by the projection formula and increase of self-intersection under birational surface pushforward. Its support is exactly $\Theta_T$ with coefficients at least one.

Choose compatible canonical divisors on $T$ and $F$. Normalize both surfaces in the field extension obtained by adjoining an $N$th root of the function with divisor $N(K_T+\Theta_T)$. This is cyclic separable of degree dividing $N$. At a prime of ramification index $r$ above coefficient $b$, the new crepant coefficient is the integer $1-r(1-b)$ by the Kummer and tame Hurwitz formulas \cite[Sections~2.41--2.43]{Kol13}. The cover of $T$ is quasi-\'etale and is \'etale near the SNC locus of every lc centre, where the function has divisor $N$ times a locally Cartier divisor, so the residual Kummer datum is a unit. Elsewhere the tame formula preserves klt. The covered pair is dlt with principal log canonical divisor.

Let $V_0$ be the minimal regular resolution of the cover of $T$. Its crepant exceptional coefficients are nonnegative by relative nefness of $K_{V_0}$ and negativity, integral and less than one by dlt, hence zero. The strict boundary $D_0$ is disjoint from the exceptional set, since $K_{V_0}+D_0\sim0$ and both terms have nonnegative degree on exceptional curves. The pullback $P_0$ of $P_T$ thus has support $D_0$ with coefficients at least one. The reduced boundary $D_0$ is SNC, and for small $\eta>0$ the effective klt pair $(V_0,D_0-\eta P_0)$ has nef and big negative adjoint. Hence $V_0$ is of del Pezzo type.

Let $Y$ be the minimal regular resolution of the normalization of $F$ in the extension, with $\pi:Y\to F$. There is a morphism $\rho:Y\to V_0$. Indeed, on a common regular resolution a last exceptional curve for the map towards $Y$ also contracts towards $V_0$. Otherwise its negative canonical degree would contradict relative nefness of $K_{V_0}$ and the nonnegative intersection with the effective exceptional discrepancy divisor. Descend through the point blowups one at a time. This uses factorization of birational morphisms of regular surfaces into blowups of closed points \cite[Tag~0C5R]{Stacks}.

Over an imperfect field, a contracted prime with $KC<0$ and $C^2<0$ has $KC=C^2=-[H^0(C,\mathcal O_C):K]$ by adjunction and integrality over its constants. Its arithmetic genus is zero and $\mathcal O_C(-C)$ has degree one over those constants. Riemann--Roch gives a degree-one map to the projective line and hence the regular point contraction \cite[Tags~0C2L and~0C2M]{Stacks}. This also proves relative nefness of the canonical divisor on the minimal regular resolutions. Intersection calculations upstairs use the actual constant field of $Y$, a finite separable extension of $K$ because $F$ is geometrically integral.

\emph{Controlling the point blowups.} Write $P'=\pi^*P$, $Q=P'^2$, and let $\Gamma$ be the reduced union of the strict transform $D$ of $D_0$ and the $\rho$-exceptional curves. The divisor $\Gamma$ is SNC at every stage, since blowing up a closed point of the regular surface with this reduced total boundary preserves normal crossings. The crepant boundary on $Y$ has the form $B_Y=D-D^-$ with $D^-\ge0$ exceptional. Indeed its coefficients are integral by tame Hurwitz. A coefficient-one curve lying over a curve of $F$ must map onto a coefficient-one component of $B_F$. Such a component has coefficient zero in $\Theta-B_F$ and survives on $T$, so this curve is not $\rho$-exceptional. A coefficient-one curve exceptional over $F$ would lie over a zero-dimensional lc centre of $(F,B_F)$, a crossing of two coefficient-one components. There $B_F$ is integral and Cartier, the Kummer datum is a unit after removing an $N$th power, and the cover is \'etale. Its minimal regular resolution introduces no curve there.

The blowup sequence also shows that the negative part is disjoint from $D$. A point on $D$ meeting only nonnegative coefficients produces a nonnegative coefficient, while a blowup off $D$ with nonpositive coefficients produces a negative one. A crossing of two positive curves cannot have been blown up as that would introduce a coefficient-one exceptional prime on $Y$.

We have $P'\ge D$ and
\begin{equation}\label{eq:upperbudget}
 P'(\Gamma-B_Y)\le\lambda Q .
\end{equation}
For curves mapping onto curves of $F$, this follows from the tame coefficient formula and \eqref{eq:pfbudget}, summing with the residue degrees. Their images belong to $\Theta$ because they are either boundary or contracted to $T$. Other curves have zero $P'$-degree. Denote by $E_i^*$ the total transforms of the individual exceptional curves of the blowups from $V_0$ to $Y$, and by $n_i$ their blowup residue degrees. The pushforward $\rho_*P'$ equals $P_0$ because the difference is effective, supported on curves exceptional over the cover of $T$ and disjoint from $D_0$, and relatively nef, so it is zero by negativity. Write
\[
 P'=\rho^*P_0-\sum_i b_i E_i^*,\qquad b_i\ge0,\qquad n_i b_i\le\lambda Q .
\]
Indeed $E_i^*\le\Gamma-B_Y$. Its coefficient at creation is at least one, and subsequent blowups add the effective excess $(2-s)E_j$ to the pullback of $K+\Gamma$ when $s\le2$ is the number of reduced components through the centre. Orthogonality and $K_{V_0}+D_0\sim0$ give
\begin{equation}\label{eq:numerical}
 P_0^2-Q=\sum_i n_i b_i^2,\qquad
 \sum_i n_i b_i=K_YP'+D_0P_0\le\lambda Q+P_0^2.
\end{equation}

If $\lambda Q\le1/4$, these inequalities give $P_0^2\le2\lambda Q$. Run a $K+2P_0$-negative program on $V_0$, pushing forward the polarization at each step, by the surface cone and contraction theorems \cite[Theorems~2.14 and~4.4]{Texc}. The negative rays are $K$-negative and the adjoint remains big, since $K+2P\simq2P-D\ge P$. Thus the steps are birational between regular surfaces. Write $l:V_0\to V$ for the output, with $K_V+2P_V$ nef and $K_V+D_V\sim0$. We have
\[
 2(l^*P_V-P_0)\le K_{V_0}-l^*K_V=l^*D_V-D_0\le l^*P_V,
 \qquad l^*P_V\le2P_0.
\]
The reduced boundary $D_V$ is lc and $P_V$ is semiample by the surface base point free theorem applied to a Cartier multiple and the klt boundary $D_V-\eta P_V$ \cite[Theorem~4.2]{Texc}. Some component $C$ satisfies $D_V C>0$. If all row sums of the boundary intersection matrix were nonpositive, its off-diagonal nonnegativity and $2a_i a_j\le a_i^2+a_j^2$ would force $P_V^2\le0$. Also $P_VC\ge D_VC/2$. The normalization of $C$ has genus at most one by lc adjunction, so $3D_V|_{C^\nu}$ is free by its degree. Choose a regular point of $C$ off the other components of degree at most $3D_VC$. It has $\varepsilon(P_V;x)\ge1/6$ by testing curves against $P_V$ and $C$ against its degree. Trace gives a generically finite system in $K_V+\ceil{aP_V}$ for fixed $a>18$. Because $\ceil{aP_V}-D_V\le\floor{aP_V}$, its free mixed degree after pullback to $Y$ is at most $2aP_0^2=O(\lambda)Q$.

Suppose $\lambda Q>1/4$. Contract a $K$-negative exceptional component $C$ of $D^-$ whenever its current normalized degree $a=PC/n$ is less than $1/(1+\lambda)$, where $n$ is the degree of its constant field. Blowing up a point on an exceptional curve raises its canonical degree by at least its initial residue degree. Thus a $K$-negative curve in the blowup sequence has no remaining later centre on it. Removing it increases $P^2$ by $na^2$, decreases $PD^-$ by $na$, and increases $P\Gamma$ by at most $na$. In fact $D^-_{\rm old}$ is the pullback of $D^-_{\rm new}$ plus $C$, while for $\Gamma$ the pullback difference has coefficient $s-1$ with $s\le2$. On the resulting $V_1$ we obtain
\[
 \begin{gathered}
 Q\le P_1^2\le2Q,\qquad P_1D_1^-\le\lambda Q,\qquad
 P_1\Gamma_1\le3\lambda Q,\\
 K_{V_1}P_1\le\lambda Q,\qquad -(D_1^-)^2\le\lambda(1+\lambda)Q .
 \end{gathered}
\]
The first line follows by summing the increments and \eqref{eq:upperbudget}. For the last inequality sum $-KC\le(1+\lambda)P_1C$ over the components of $D_1^-$ with their coefficients. The bound holds because $-KC=n$ for remaining $K$-negative curves, and $KC=D_1^-C$ on this support.

For a contraction from $V_1$ to a curve, a reduced integral fibre $G$ over a rational point satisfies $P_1G\ge1/(4\lambda)$. If a component of the transform of $D_0$ is horizontal, its coefficient gives $P_1G\ge1$. Otherwise take the last exceptional curve $E$ in the blowup sequence horizontal for the contraction. Such a curve exists since $P_1G>0$ and $P_1$ is supported on the boundary and exceptional curves. The contraction descends past the later blowups, whose exceptional curves are vertical. On the model creating $E$, let $P_i,G_i$ be the corresponding polarization and fibre. For $n$ the residue degree and $\delta=P_1G$ we have
\[
 E^2=-n,\quad G_i^2=0,\quad \ell=G_iE\ge n,\quad
 P_iE=nb_i\le\lambda Q,\quad P_iG_i=\delta,\quad P_i^2\ge Q.
\]
The Hodge inequality on the span of $E,G_i$ gives
\[
 Q\le 2\delta(P_iE)/\ell+n\delta^2/\ell^2
   \le2\lambda Q\delta+\delta^2.
\]
Since $\lambda Q>1/4$, this implies the asserted lower bound. The same computation holds after the regular ground-field extension to the generic diagonal.

The surface $V_0$ is of del Pezzo type in characteristic at least seven, hence geometrically integral with $H^1(\mathcal O)=0$ \cite[Corollary~5.5 and Theorem~5.7]{BT}. These properties persist under the point blowups and under separably generated ground-field extension. Over the function field of $V_1$, consider any integral curve $C$ through the rational diagonal point $x$. If $h^0(\mathcal O(C))=1$, the exact sequence for $C$ gives $H^0(\mathcal O_C(C))=0$. Its Hilbert point has zero relative tangent space and therefore a finite separable field of definition over the original field. The function field extension is regular, so the curve descends and cannot contain $x$. Thus $C$ moves and $C^2\ge0$. If $C^2=0$, it and a linearly equivalent member define a free pencil, because they are disjoint. The Stein factor is trivial since one fibre is reduced integral and contains a rational point. The preceding fibre estimate applies.

By Hodge index perpendicular to $P_1$ we have
\[
 K_{V_1}C\le D_1^-C
 \le(\lambda+\sqrt{2\lambda^2+\lambda})P_1C\le3\lambda P_1C .
\]
If $j=\operatorname{mult}_x C$, the strict transform on the point blowup has constant field degree at most $j$, by its intersection with the exceptional divisor. Adjunction there gives $C^2+KC\ge j(j-1)-2j$. If $P_1 C/j<1/(100\lambda)$, then $C^2\le(P_1C)^2/P_1^2$, so the two inequalities force $j\le3$ and $C^2=0$. This contradicts the fibre estimate. Therefore $\varepsilon(P_1;x)\ge1/(100\lambda)$. Use the ample perturbation trace with $a=\ceil{600\lambda}+1$ to obtain a generically finite system in $K_{V_1}+\ceil{aP_1}$. The displayed bounds give mixed degree $O(\lambda)Q$. Pullback to $Y$ preserves this upper bound since $P'\le g^*P_1$ for the blowdown $g$.

\emph{Descent of the system.} In both regimes the tame cover has a generically finite free system on a birational model, of mixed degree $O(\lambda)Q$. Push its divisors and sections down by the norm to $F$. For any generating ratio $z_i=s_i/s_0$, the norms of $s_i-a s_0$ for finitely many constants $a\in K$ recover the coefficients of $\operatorname{Nm}(z_i-T)$ by interpolation. Thus the upstairs ratios are algebraic over the field of downstairs ratios and the system downstairs is generically finite. By the projection formula and the degree of the cover, its free part has mixed degree $O(N\lambda)P^2$. The initial Seshadri argument now yields the required $m$.
\end{proof}

\section{Exceptional complements}\label{s:exceptional}

\subsection{Pencils on canonical threefolds with a two-form}\label{s:twoforms}
For the canonical Fano threefolds in this subsection, a nonzero regular two-form on a smooth resolution forces the generic surface of the pencil to have pseudo-effective canonical divisor. We use exceptional curves on this surface to find a bounded Cartier multiple of a big adjoint divisor.

\begin{lemma}\label{lem:twoformpencil}
Let $X$ be a canonical Fano threefold, $L=-K_X$, $V=L^3$, and $d$ the first effective anticanonical degree. Suppose $p\ge11$ and a smooth projective resolution of $X$ has a nonzero regular two-form. There is an absolute $C$ such that a nonconstant flat family of prime divisors of $L^2$-degree $cV$ implies $d\le C(c+1)$.
\end{lemma}

\begin{proof}
\emph{The canonical divisor of the generic surface.} Take a smooth projective birational model $g:U\to X$ with a contraction $f:U\to C_0$ to a smooth curve. Let $P=g^*L$ and $F$ be its generic surface. If $K_F$ is not pseudo-effective, the relative $K_U$-MMP ends in a terminal Mori fibre space \cite[Theorems~1.1, 1.3 and~1.4]{HNT}. The regular two-form is preserved on smooth projective models by properness in codimension one and extension across codimension two. It is regular on the smooth locus of the terminal output by the same argument.

If the Mori base is a curve, its generic fibre is a regular del Pezzo surface with the base function field as constants, since terminal threefold singularities are isolated. It is smooth by \cite[Proposition~5.2]{BT}. On a smooth closed fibre $G$ the cotangent sequence filters the restriction of the two-form by $\Omega^1_G$ and $K_G$. Both spaces of global sections vanish for the smooth rational del Pezzo surface. This forces the form to vanish on a dense open of the threefold, a contradiction.

If the Mori base is a surface, the generic fibre is a regular conic, smooth in this characteristic \cite[Lemma~2.17]{BT}. Resolve the base to $Z$ and the graph to a smooth $g':U'\to X$ with contraction $h:U'\to Z$. The actual generic fibre $G$ is unchanged. The cotangent sequence on $G$ and $H^0(G,\Omega^1_G)=0$ identify the two-form with the pullback of a rational two-form on $Z$. A point of degree at most two on $G$ gives a multisection of separable degree at most two. The form pulls back regularly on its resolution. At any prime of $Z$ the ramification above it in that field is tame, and its differential correction $r-1$ cannot cancel a pole of order $r$ or more. Thus we can take $K_Z\ge0$.

Put $P'=-(g')^*K_X$, so $K_{U'}=-P'+E$ for $E\ge0$ exceptional. Since $P'G\ge2$, set $\lambda=2/(P'G)\le1$. For a very ample $H$ on $Z$, take $D_0\ge0$ in $|mr\lambda P'-h^*H|$ for large $m$, where $r\lambda P'$ comes from a very ample Cartier divisor on $X$. For small rational $a>0$ and large $N$ use general $D_i\in|r\lambda P'|$ to form
\[
 \Delta'=aD_0+\frac{1-amr}{rN}\sum_{i=1}^N D_i
       \simq \lambda P'-ah^*H .
\]
Multiplicities in the fixed system are uniformly bounded by intersection with two very ample divisors through the point. Choose $D_i$ with no four meeting, so small $a$ and large $N$ make the multiplicity of $\Delta'$ less than one at every point. It is klt by Lemma~\ref{lem:graphlct}. On the geometric generic conic the $D_i$ can be chosen reduced and disjoint, avoiding $D_0$ and the exceptional locus there. Indeed the map given by the very ample system on $X$ is immersive along the generic conic away from finitely many points, and the system has no base points. The open conditions over the geometric generic field can be met over the infinite field $k$. Reducing $a$ also makes the geometric generic boundary klt.

The relative $(K_{U'}+\Delta')$-MMP and semiampleness with relatively big boundary express the resulting log divisor as pulled back from a surface birational to $Z$ \cite[Theorem~1.4]{HNT}. The geometric generic pair is unchanged since the log degree on the conic is zero. The curve canonical bundle formula makes the base divisor rationally equivalent to the canonical divisor plus an effective divisor \cite[Proposition~3.2]{WitCBF}. Its canonical divisor is effective by the two-form on $Z$. Negativity on a common resolution now gives $P'^2(K_{U'}+\Delta')\ge0$. Yet
\[
  P'^2(K_{U'}+\Delta')=(\lambda-1)V-aP'^2h^*H<0 ,
\]
since a general pullback of $H$ is not exceptional over $X$. This rules out the conic case, and so $K_F$ is pseudo-effective.

\emph{A big adjoint with at most ten marked components.} Resolve so the reduced exceptional divisor $\Gamma$ over $X$ is SNC. The effective exceptional divisor $K_U+P$ is big on $F$ since $P|_F$ is big and $K_F$ is pseudo-effective. Thus $\Gamma|_F$ is big. We claim that at most ten of its primes suffice for $K_F+\Gamma_0|_F/3$ to be big, where $\Gamma_0$ denotes their closures on $U$.

The relative $K_U$-MMP gives a morphism $\rho:F\to F_0$ to a regular surface with $K_{F_0}$ nef, since flips are isomorphisms on the generic fibre and the terminal output is regular there. Moreover $K_F-\rho^*K_{F_0}$ is positive along each contracted curve. Apply surface abundance over the generic field \cite[Theorem~1.1]{Tsurf}. If $\kappa(K_{F_0})=2$, no primes are needed. If $\kappa(K_{F_0})=1$, some component of $\rho_*(\Gamma|_F)$ dominates the Iitaka base, since the pushforward is big. Adding this component gives a big divisor.

If $\kappa(K_{F_0})=0$, then $K_{F_0}\equiv0$. Write $D=\rho_*(\Gamma|_F)$ with components $D_i$, and let $a_i$ be the constant-field degree of $D_i$. Adjunction gives $D_i^2/a_i=2p_a(D_i)-2$ over its constant field. The intersection form on the support of $D$ has a positive direction, since the free moving part of a large multiple pushes down to a nef subdivisor of positive square. A positive diagonal supplies one curve, and a zero diagonal with nonzero off-diagonal entry supplies two. Otherwise discard isolated zeros. The remaining integral matrix $(-D_iD_j/a_i)$ is a symmetrizable generalized Cartan matrix. A minimal connected principal submatrix with indefinite symmetrization is hyperbolic and has rank at most ten \cite[Section~3]{Carbone}. The intersection form on that support has a positive direction with nonnegative coordinates, since replacing coordinates by absolute values can only increase the square. An effective rational divisor on these primes with positive square is big, so their reduced sum is big. Lift the selected primes to $F$ and subtract a small positive multiple of the pullback of the big adjoint on $F_0$, using the positive exceptional difference of canonical divisors and the pseudo-effectivity of $K_{F_0}$.

\emph{The Cartier multiple on the actual fibre.} Run the $(K_U+\Gamma_0/3)$-MMP over $C_0$ to $(Y,\Delta_Y)$. Put $J_Y=K_Y+\Delta_Y$ and $S=Y_\eta$. We show that $I J_Y|_S$ is Cartier for $I=3(2^{10})!$. In a generic surface step with canonical source the minimal regular resolution is crepant by relative canonical nefness and negativity. For a curve $D$ there, with strict transform $D'$ on that resolution, $D^2\ge(D')^2$ and $KD=K_{\rm res}D'$. If $KD>0$, let $a$ be the constant-field degree of $D'$ and $j=-(D')^2/a$. Adjunction and integrality give $K_{\rm res}D'/a\ge\max(1,j-2)$. In particular $KD+\min(0,D^2)/3\ge0$. Intersections with other boundary components are nonnegative, so a log negative contraction is $K$-nonpositive. It preserves canonicity, and all $K$-trivial contracted curves are boundary components. Thus at most ten such curves are contracted. A $K$-negative contracted curve has positive discrepancy afterwards. The minimal resolution of $S$ is crepant, is dominated by $F$, and consequently has at most ten exceptional curves.

At a point of $S$ the matrix $(-E_iE_j/a_i)$ of these curves, where $a_i$ is their constant-field degree, is integral with diagonal $2$ by adjunction. It is similar to a positive definite symmetric matrix with the same diagonal, so its determinant is at most $2^{10}$. The numerical pullback of a Weil divisor is determined by this matrix with integral right-hand side. Its determinant clears the exceptional coefficients. An integral numerical pullback descends at a rational surface singularity \cite[Section~14 and Proposition~17.1]{Lipman}. Canonical surface singularities over the generic field are rational \cite[Proposition~2.6]{BM}. Thus $(2^{10})!$ clears their Weil indices, proving the claim for $I J_Y|_S$.

\emph{Returning sections to $X$.} Take $H_0$ free Cartier of positive degree on $C_0$, set $M=f_Y^*H_0$, and put $N=J_Y+7M$. The ray length bound $6$ gives $N-M$ nef, and $N$ is nef and big since $J_Y|_S$ is nef and big. The base field $k(C_0)$ is relatively algebraically closed in $k(Y)$, and its single $p$-basis element remains outside $k(Y)^p$. Relative algebraic closedness and separable generation give geometric integrality of $S$. Over $\overline{k(C_0)}$ minimally resolve the normalization of the base change. On this resolution $R$ let $A$ be the pullback of $I J_Y|_S$. It is nef and big Cartier with $K_R A\le A^2/I$. Indeed the conductor formula gives the effective difference between the base-changed canonical divisor and that of the normalization \cite[Theorem~1.1]{PW}. Minimal resolution and negativity retain the inequality after intersecting with $A$. The pullback of the boundary is effective.

Descend the resolution, its Cartier polarization and its common birational isomorphism open to a finite extension of the parameter field, and spread to an open of that base curve. The generic Cartier multiple extends over an open by properness. On the resulting smooth surfaces its pullbacks are nef and big, with the same intersection inequality. Lemma~\ref{lem:surface} gives a Seshadri lower bound $1/4$ for $A$ on a dense open of such a surface. At a corresponding smooth point of $Y$, the vertical curves have Seshadri ratio at least $1/(4I)$ for $N$. The horizontal curves have ratio at least one, since $N-M$ is nef and the fibre divisor through the point is numerically $M/\deg H_0$. The point can be chosen off any prescribed proper closed subset.

Take $u=48I+2$. Write $uN\simq H'+E'$ with $H'$ ample and $E'\ge0$. The positive class $uN-E'/2$ has Seshadri constant at least $u/(8I)>6$ at a suitable point. Apply Lemma~\ref{lem:trace} with Weil target $K_Y+\ceil{uN}$ and boundary $\ceil{uN}-uN+E'/2$. Since $u\equiv-1\bmod 3$, the target is
\[
 K_Y+\ceil{uN}=(u+1)J_Y+7uM .
\]
Suppose $P^2 f^*H_0\le cV$. On a common smooth graph $W$, the free moving part $B_W$ of the birational system satisfies $P_W^2B_W\le7ucV$. Indeed the MMP comparison bounds the contribution from $J_Y$ by the degree of $K_U+\Gamma_0/3$ against $P^2$, which is $-V$. For $j+1>84u(c+1)$ the nef Morse inequality
\[
 \vol((j+1)P_W-4B_W)\ge (j+1)^3V-12(j+1)^2P_W^2B_W
\]
gives bigness. At a general smooth point off an effective representative of the difference, the hyperplanes of $B_W$ give $(j+1)\varepsilon(L;x)\ge4$. The trace criterion then gives a section of $jL$.

Finally lift general members of the given family to a smooth resolution. On a smooth parameter open their strict transforms form a Cartier family and are linearly equivalent by the same rational-chain-connectedness and Picard argument as in Proposition~\ref{pr:moving}. Resolve a pencil of two distinct members and take its Stein factorization over $C_0$. For the pullback $H_0$ of a point of $\PP^1$, we have $P^2 f^*H_0\le cV$ since removing the fixed part can only lower the degree. Taking $j$ as above proves the assertion, for example with $C=84u+2$.
\end{proof}

\begin{theorem}\label{th:twoformscaled}
The exceptional canonical Fano threefolds whose smooth resolutions admit a nonzero regular two-form have a common klt complementary index and are arithmetically bounded in sufficiently large characteristic.
\end{theorem}
\begin{proof}
Let $d$ be the first effective anticanonical degree and $V=(-K_X)^3$. Lemma~\ref{lem:window} gives a birational system in degree $b\le b_0d$ with $b_0$ fixed. For $d\ge2$ apply Theorem~\ref{th:jetsurface}\textup{(ii)} with $v=16$ on a small $\Q$-factorialization. If $d>Kl$, the family of primes descends to $X$ with normalized degree $<l/8$. Then Lemma~\ref{lem:twoformpencil} gives $d\le C(1+l/8)$. Thus $d=O(1+V^{-1/3})$, and the exceptional bound $V\le27$ gives a bound on $d^3V$. Proposition~\ref{th:scaled} gives the index and arithmetic boundedness.
\end{proof}

\subsection{Anticanonical nonvanishing and local Riemann--Roch}\label{sec:canslope}\label{s:terminalcorr}
We work on a crepant terminalization, a projective birational morphism $f:Y\to X$ with $Y$ terminal and $K_Y=f^*K_X$. Terminal threefolds have isolated singularities, by regularity of terminal surface germs \cite[Corollary~2.30(1)]{Kol13}. At each singular point the local Riemann--Roch correction records how the rounded pullback of $mK_Y$ varies with $m$. These corrections detect instability of the tangent bundle on a resolution. In the unstable case we use a rational curve of small anticanonical degree to bound the volume by counting jets at a singular point.

\begin{lemma}\label{lem:resolutionrr}
Let $T$ be a projective terminal threefold in characteristic $p>5$ and $g:Z\to T$ a projective log resolution which is an isomorphism over the regular locus. Put $A=-g^*K_T$, $V=A^3$, $c=c_2(T_Z)\cdot A$ and $\chi_0=\chi(\cO_T)$. Then for every integer $m$,
\begin{equation}\label{eq:canslope-smooth-rr}
 \chi(T,-mK_T)=\chi_0+\frac{m(m+1)(2m+1)}{12}V
                          +\frac m{12}c+\sum_Q\gamma_Q(m).
\end{equation}
Here $\gamma_Q(m)$ is a local intersection term over the non-Gorenstein point $Q$. It vanishes when $mK_T$ is Cartier near $Q$. If $T\to X$ is a projective crepant birational morphism to a canonical threefold, then $\chi(T,-mK_T)=\chi(X,-mK_X)$.
\end{lemma}

\begin{proof}
For a log resolution $h:U\to W$ of a canonical threefold, use
\[
 D_m=K_U-\lfloor(m+1)h^*K_W\rfloor
     =-mh^*K_W+(K_U-h^*K_W)+\{(m+1)h^*K_W\}.
\]
Pullback of effective $\Q$-Cartier divisors and nonnegativity of the discrepancy divisor give $h_*\cO_U(D_m)=\cO_W(-mK_W)$. Indeed, a section of the latter is effective after pullback and addition of the two remaining terms, and effectivity upstairs pushes down. The fractional part in the display is exceptional with coefficients in $[0,1)$. Grauert--Riemenschneider vanishing thus gives $R^ih_*\cO_U(D_m)=0$ for $i>0$ \cite[Definition~1 and Theorem~3]{BK}. On a common log resolution the divisors $D_m$ for $T$ and $X$ agree by crepancy, so Leray gives the comparison.

Apply smooth Riemann--Roch to $D_m$ on $Z$, writing $D_m=mA+P_m$ with $P_m=K_Z+A+\{(m+1)g^*K_T\}$. Subtract the formula at $m=0$. Products involving $A$ and an exceptional divisor vanish, since a Cartier multiple of $A$ is trivial on each exceptional fibre. Exceptional supports over distinct points are disjoint. This gives \eqref{eq:canslope-smooth-rr}, with $\gamma_Q(m)$ the difference between the values at $m$ and at $0$ of
\[
 \frac{P^3}{6}-\frac{P^2K_Z}{4}
                +\frac{P(K_Z^2+c_2(T_Z))}{12},
 \qquad P=\text{the part of }P_m\text{ supported over }Q .
\]
If $mK_T$ is Cartier near $Q$, the two fractional parts agree there. In particular Gorenstein points contribute zero.
\end{proof}

We describe the locality used in computing these terms. If $T$ is $\Q$-factorial, a resolution over an isolated point is the blowup of a primary ideal. To see this, replace a relatively ample Cartier divisor by a multiple and subtract the pullback of its Cartier pushforward. The resulting divisor is exceptional. The sections of its nonnegative multiples form ideals, since the sections have no poles downstairs in codimension one. A Veronese subalgebra generated in degree one presents the morphism as the blowup of the degree-one ideal. This ideal is the unit ideal away from the point.

Such an ideal contains a power of the maximal ideal and can be transferred through an isomorphism of completed local rings. By Artin approximation with agreement modulo that power \cite[Tag~0CAV]{Stacks}, the two blowups agree on a common pointed \'etale neighbourhood. The exceptional fibres, their discrepancy coefficients and the intersection numbers with tangent Chern classes therefore agree. One can use this argument also when the comparison variety lacks $\Q$-factoriality, by transporting the ideal from the $\Q$-factorial side. For a fixed $\Q$-factorial $T$ the numbers $\gamma_Q(m)$ are independent of the resolution. Indeed, $c$ is determined by \eqref{eq:canslope-smooth-rr} at a positive Cartier multiple of $K_T$, and the primary ideals at distinct points can be chosen independently. Comparing resolutions that differ only at $Q$ then gives the assertion.

For relatively prime $a,r$ with $r\ge2$, a cyclic chart of type $\frac1r(1,a,-a)$ means the degree-zero subring of $k[[x,y,z]]$ with these weights modulo $r$. This definition includes $p\mid r$. Write $b\equiv a^{-1}\pmod r$ and let $[u]_r$ denote the residue in $\{0,\ldots,r-1\}$. Put
\[
 \sigma_r=r-r^{-1},\qquad
 \ell_{r,b}(m)=\sum_{j=1}^m \frac{[jb]_r(r-[jb]_r)}{2r}
 \quad(m\ge0).
\]
A cyclic basket will mean a finite multiset of such pairs $(r,b)$. A point can contribute several pairs, with their corrections added.

\begin{lemma}\label{lem:cyclicrr}
At a point with a cyclic chart as above, the canonical index is $r$, the multiplicity is at most $2r$, and the local term of \eqref{eq:canslope-smooth-rr} is
\[
 \gamma_Q(m)=\frac{m\sigma_r}{12}-\ell_{r,b}(m)\quad(m\ge0),
 \qquad \gamma_Q(-1)=-\frac{\sigma_r}{12}.
\]
\end{lemma}

\begin{proof}
The affine cyclic chart is the first orthant in the lattice $\Z^3+\Z(1,a,-a)/r$. All its proper faces are regular. The functional with value one on the three coordinate rays has denominator $r$ in the dual lattice, giving the index. Take a projective regular subdivision keeping those faces. Its morphism is the blowup of an ideal primary to the vertex, so the subdivision computes the local term of the completed chart by the preceding comparison.

To compare with characteristic zero, complete the cone to a projective simplex fan and resolve outside the cone, keeping the cone fixed. Then perform the chosen subdivision of the cone. The last fan defines a smooth projective toric scheme over $\Z$. The discrepancies and the terms of the Riemann--Roch polynomial supported over the vertex agree on all fibres, by the lattice description and specialization of the exceptional intersection numbers in this smooth scheme.

Over $\C$ the compactification before the last subdivision has one terminal cyclic point. Reid's formula for the correction of the divisor locally equivalent to $iK$ \cite[(10.1) and Theorem~10.2]{ReidYPG}, with $i=[-m]_r$, gives
\[
 -\frac{i\sigma_r}{12}
    +\sum_{j=1}^{i-1}\frac{[jb]_r(r-[jb]_r)}{2r}.
\]
The polynomial part of Reid's formula agrees with that of \eqref{eq:canslope-smooth-rr}, as is seen at Cartier multiples. The displayed term gives the stated expressions by symmetry and $\ell_{r,b}(r)=r\sigma_r/12$. Finally $(x^r,y^r+z^r,yz)$ is an invariant parameter ideal. Its extension to $k[[x,y,z]]$ has colength $2r^2$. Taking the degree-zero summand divides the length by $r$, since the factor $k[x]/(x^r)$ has one basis element of every weight. The invariant ring is Cohen--Macaulay, so its multiplicity is at most the colength $2r$ of this parameter ideal.
\end{proof}

\begin{lemma}\label{lem:chieuler}
If $X$ is a canonical ample Fano threefold in characteristic $p>5$, then $\chi(\cO_X)\ge1$. If a smooth projective resolution $Z$ of $X$ has $H^0(Z,\Omega_Z^2)=0$, then $\chi(\cO_X)=1$.
\end{lemma}

\begin{proof}
By \cite[Theorem~4.8 and Proposition~4.11]{GNT} the resolution is rationally chain connected. Every morphism to an abelian variety is constant on rational chains, so $\operatorname{Pic}^0(Z)$ has dimension zero. The Picard deformation inequality gives
\[
 \dim\operatorname{Pic}^0(Z)\ \ge\
 h^1(\cO_Z)-h^2(\cO_Z).
\]
Recall why this applies to a possibly nonreduced Picard scheme. Write its completed local ring as $S/I$, with $S$ a power series ring in $h^1(\cO_Z)$ variables and $I\subset\mathfrak m_S^2$. The obstruction to lifting the universal line bundle across a square-zero extension with kernel $I/\mathfrak m_S I$, obtained by truncating $S/\mathfrak m_S I\to S/I$ in sufficiently high order, lies in $H^2(\cO_Z)\otimes I/\mathfrak m_S I$. Its image in each one-dimensional quotient of the kernel is nonzero. Indeed, a lift of the class would lift the tautological map from $S/I$, and changing the lifts of the parameters by elements killed by $\mathfrak m_S$ leaves the images of the relations in $I\subset\mathfrak m_S^2$ unchanged. This bounds the number of equations by $h^2(\cO_Z)$.

A section of $K_Z$ would push down to a section of $K_X$, so $h^3(\cO_Z)=0$. The rationality of the singularities of $X$ \cite[Theorems~3 and~19]{BK} gives $\chi(\cO_X)=\chi(\cO_Z)\ge1$. For the upper bound we need only the following consequence of the Hodge--Witt Euler formula and its symmetry in dimension three \cite[Theorem~4.4 and \S4.2]{XZ}. The $h_W^{0,j}$ are nonnegative, and hence
\[
 \chi(\cO_Z)\le h_W^{0,0}(Z)+h_W^{2,0}(Z)
              \le 1+h^0(Z,\Omega_Z^2).
\]
\end{proof}

\begin{proposition}\label{pr:basketslope}
Let $X$ be an exceptional canonical Fano threefold with zero boundary with $\chi(\cO_X)\le H$, and let $f:Y\to X$ be a projective crepant $\Q$-factorial terminalization. Suppose each non-Gorenstein point $Q$ has a cyclic basket expressing $\gamma_Q$ as a sum of the terms of Lemma~\ref{lem:cyclicrr}. Suppose the basket contains the local index and that $e_Q(Y)\le M\sigma_Q$, where $\sigma_Q$ is the sum of $\sigma_r$ over that basket. Then in sufficiently large characteristic the first effective degree of $-K_X$ is bounded in terms of $H$ and $M$.
\end{proposition}

\begin{proof}
Set $L=-K_X$, $A_Y=f^*L$ and $A=g^*A_Y$, where $g$ is the resolution of Lemma~\ref{lem:resolutionrr}. Put $V=L^3>0$, $c=c_2(T_Z)\cdot A$, $\chi_0=\chi(\cO_X)$ and $\sigma=\sum_Q\sigma_Q$. Since $\chi(X,K_X)=-\chi_0$ by Cohen--Macaulay Serre duality, evaluating \eqref{eq:canslope-smooth-rr} at $m=-1$ gives $c=24\chi_0-\sigma$. For $m\ge0$ the formula becomes
\begin{equation}\label{eq:canslope-rr}
 \chi(X,mL)=(2m+1)\chi_0+\frac{m(m+1)(2m+1)}{12}V
                      -\sum_{(r,b)}\ell_{r,b}(m),
 \qquad c=24\chi_0-\sigma.
\end{equation}
The sum is over all the baskets.

Write $d$ for the first effective degree. We may assume $d\ge2$. Take $p$ large enough for Lemma~\ref{lem:window}, which gives $\chi(X,(d-1)L)\le0$. Each summand in $\ell_{r,b}(d-1)$ is at most $r/8\le\sigma_r/6$. Together with $\chi_0\ge1$ this yields
\begin{equation}\label{eq:canslope-sigma-lower}
 \sigma\ \ge\ \frac{d(2d-1)}2 V+12\chi_0
           \ \ge\ \frac34 d^2V .
\end{equation}
If $\sigma\le48H$, all local indices of $K_Y$ are at most $48H+1$. Their common multiple then gives a positive lower bound for $V$, since the cube of a Cartier multiple of $A_Y$ is a positive integer. This bounds $d$.

Suppose $\sigma>48H$, so $c<-\sigma/2$. Choose an ample Cartier divisor $H_Z$ and a small positive rational $\eta$, and set $B=A+\eta H_Z$. With
\[
 v=-K_ZB^2,\qquad u=K_Z^2B,\qquad w=B^3,\qquad
 c_B=c_2(T_Z)\cdot B
\]
we can arrange $u,v>0$, $c_B<-\sigma/2$, $w/v<8/7$ and $v^2/w<2V$. Indeed $u,v,w$ tend to $V$ and $c_B$ tends to $c$.

Use the slope $\mu(\mathcal E)=c_1(\mathcal E)B^2/\operatorname{rk} \mathcal E$ and put $\theta=\mu_{\max}(T_Z)/v$. The Harder--Narasimhan filtration has semistable successive quotients of decreasing slope. Its first and last slopes are denoted by $\mu_{\max}$ and $\mu_{\min}$. Since $\mu(T_Z)=v/3$, we have $\theta\ge1/3$ and $\mu_{\min}(T_Z)\ge-2\theta v$. Consequently $\mu_{\max}(\Omega_Z^1)\le2\theta v$. Denote by $L_{\max},L_{\min}$ the limits of the extreme slopes under $e$ Frobenius pullbacks divided by $p^e$. For a rank $s$ sheaf, Langer's slope deviation bound is $(s-1)[L_{\max}(\Omega_Z^1)]_+/p$ \cite[Corollary~6.2]{LangerSlopes}. Applying it first to $\Omega_Z^1$ gives $[L_{\max}(\Omega_Z^1)]_+\le2p\theta v/(p-2)$. The bound for $T_Z$ now gives $L_{\max}(T_Z)\le2\theta v$ and $L_{\min}(T_Z)\ge-3\theta v$ for $p\ge7$.

Langer's discriminant inequality \cite[Theorem~5.1(5.1.1)]{LangerSlopes} applies after clearing the Cartier denominator of $B$. It gives
\[
 w(6c_B-2u)\ge
 -9\bigl(L_{\max}(T_Z)-v/3\bigr)\bigl(v/3-L_{\min}(T_Z)\bigr)
 \ge -72\theta^2v^2.
\]
Using $u>0$ and $c_B<-\sigma/2$, we get
\begin{equation}\label{eq:canslope-sigma-upper}
 \sigma<24\theta^2v^2/w<48\theta^2 V .
\end{equation}
Combining the two bounds for $\sigma$ gives $d<8\theta$. We may therefore assume $\theta>8$.

The first Harder--Narasimhan subsheaf $\mathcal F\subset T_Z$ is saturated of rank one or two. It is a $1$-foliation, meaning that it is closed under Lie brackets and $p$-th powers of derivations. In rank one the bracket vanishes modulo $\mathcal F$ and $F_Z^*\mathcal F$ has slope $p\theta v$, greater than the maximal slope of the quotient. In rank two the quotient has slope $(1-2\theta)v<0$. The determinant of $\mathcal F$ has positive slope, and the slope deviation bound gives $L_{\min}(\mathcal F)\ge\theta v(1-2/(p-2))>0$. Thus the bracket map from $\bigwedge^2\mathcal F$ to the quotient vanishes, as does the map from $F_Z^*\mathcal F$ induced by $p$-th powers after involutivity.

Take a general complete intersection curve of two very ample multiples of $B$ through a general point, away from the locus where $\mathcal F$ fails to be a subbundle. Foliated bend-and-break \cite[Theorem~2.1]{LangerFol} gives a rational curve through that point of $B$-degree at most
\[
 \frac{6p B^3}{((p-1)c_1(\mathcal F)-K_Z)B^2}
 \le \frac{6p}{(p-1)\theta+1}\frac wv
 < \frac8\theta.
\]
The normalizations of these curves occur among maps from $\mathbb P^1$ of finitely many degrees with respect to a fixed very ample multiple of $B$. One irreducible family of maps therefore dominates $Z$. We restrict to maps meeting the isomorphism open of $Z\to Y\to X$. The degree $\delta$ of the pullback of $A_Y$ to $\mathbb P^1$ is constant on the family and satisfies $0<\delta<8/\theta<1$. A curve of this degree must meet a non-Gorenstein point of $Y$, because the degree of a line bundle on $\mathbb P^1$ is an integer. The loci of maps meeting each such point are closed by properness of $\mathbb P^1$. The family is irreducible and the list of points finite, so one point $Q$ lies on every image in $Y$.

For sufficiently divisible $j$, a section of $jA_Y$ is determined by its germ modulo $\mathfrak m_Q^{\lfloor j\delta\rfloor+1}$. Indeed, a nonzero section in the kernel restricts nontrivially to some covering curve and would vanish above $Q$ to order greater than the degree $j\delta$. Comparing the Hilbert polynomial for $jL$ with the Hilbert--Samuel polynomial at $Q$ gives
\[
 V\le e_Q(Y)\delta^3
    \le M\sigma\delta^3
    < \frac{48M\,8^3}{\theta} V .
\]
This bounds $\theta$, and hence $d$, in terms of $M$.
\end{proof}

We compute the baskets required in Proposition~\ref{pr:basketslope}. A local trivialization of $rK_X$ makes the first $r$ reflexive powers of the canonical sheaf into a cyclically graded algebra. The spectrum of this algebra is the \emph{canonical cover}. It is a $\mu_r$-torsor over the regular locus. The grading works over the given field even when the torsor is inseparable.

\begin{lemma}\label{lem:terminalhyp}
Let $(X,Q)$ be a terminal threefold germ in characteristic $p>5$ over an algebraically closed field, and let $r$ be its canonical index. For a suitable local trivialization of $rK_X$, the canonical cover is normal with a unique point over $Q$ and is regular elsewhere near that point. Its completion is either regular or a hypersurface in four variables with an equation of order two. The variables can be chosen homogeneous for the $\Z/r\Z$-grading and the equation semi-invariant.
\end{lemma}

\begin{proof}
The graded pieces of the canonical algebra are Cohen--Macaulay \cite[Theorems~3 and~19]{BK}. A homogeneous element of nonzero degree has $r$th power in the maximal ideal at $Q$, since a homogeneous unit would trivialize the corresponding canonical power and contradict minimality of $r$. The fibre over $Q$ thus has a single point with residue field $k$.

For $p\nmid r$ the cover is \'etale over the regular locus. For $p\mid r$, first multiply the trivialization by a general affine-linear function $u$ of an affine embedding, nonzero at $Q$. Locally on the regular locus the new equation is $T^r-au$ with $a$ a unit. At a fixed regular point, $d(au)=0$ imposes three independent linear conditions on $u$, since affine-linear functions generate first jets. The critical locus is at most finite for general $u$. Shrink about $Q$ to remove it and the zeros of $u$. The Cohen--Macaulay cover is regular off its closed point and hence normal.

To obtain the hypersurface presentation we keep track of the field of definition. Consider the full canonical algebra
\[
 X^\sharp=\Spec_X\bigoplus_{j\in\Z}\cO_X(jK_X).
\]
The chosen homogeneous unit $s$ of degree $r$ makes it finite over $X\times\mathbb G_m$. It is Cohen--Macaulay by the same decomposition. Over the regular locus it is a torus torsor, so it is also normal. Write $K_X=\divi(\omega)$ and let $t$ be the rational coordinate of degree one. On the torus torsor the form $\Omega=\omega\wedge(dt/t)/t$ trivializes the canonical sheaf. Reflexive extension makes $X^\sharp$ Gorenstein.

On a log resolution $\pi:U\to X$ the normalization of the main component of the pullback of $X^\sharp$ has graded algebra $\bigoplus_j\cO_U(\lfloor j\pi^*K_X\rfloor)t^j$. Indeed the algebra is saturated in local normal crossings coordinates, and the $r$th power of each homogeneous element is in the original pullback algebra. Thus the normalized pullback is proper over $X^\sharp$. In normal crossings parameters $x_i$, put $d_i=\operatorname{ord}_{(x_i=0)}(\pi^*K_X)$. On a regular toroidal subdivision a ray has nonnegative integers $v_i=\operatorname{ord}(x_i)$ and $\operatorname{ord}(t)=\sum_i d_i v_i$. Taking the order of $\Omega$ in logarithmic coordinates gives the log discrepancy
\begin{equation}\label{eq:torusdiscrepancy}
  \sum_i(k_i+1-d_i)v_i=\sum_i A_X((x_i=0))v_i ,
\end{equation}
where $k_i$ is the order of $\omega$ on $U$. If $\sum_i v_i\ge2$ this is at least two. A ray with $\sum_i v_i=1$ has integral slope $d_i$ and, if exceptional, comes from an exceptional divisor on $U$, of discrepancy greater than one. The normal crossings criterion then gives terminality. To principalize the ideal of the reduced locus over $Q$, include the coefficient divisors of a finite set of homogeneous generators in the normal crossings support and refine the subdivision to principalize the resulting monomial ideal.

The morphism $s:X^\sharp\to\mathbb G_m$ is flat by torsion-freeness over $k[s,s^{-1}]$. Its generic fibre over $F=k(s)$ is Gorenstein and terminal by localization of the dualizing sheaf and the discrepancy formula. It is regular off the possible singular point over $Q$, which is $F$-rational because the reduced locus over $Q$ is $\Spec k[s,s^{-1}]$. The preceding principalization gives a regular model principalizing that rational point.

We apply the rational-point argument of \cite[Theorem~6.4]{Tell}. A generic hyperplane through the point over a purely transcendental extension of $F$ preserves multiplicity and decreases embedding dimension by one. Its strict transform is regular on the regular model and normal on the normalized point blowup by generic Bertini on the free parts \cite[Theorem~6.1]{Tell}. Downstairs the generic hyperplane $H$ is regular off the marked point by the same theorem. As a Gorenstein Cartier cut it satisfies $S_2$, hence it is normal by Serre's criterion.

Adjunction for the regular transform gives $\mathfrak m_H\omega_H\subseteq g_*\omega_{\widetilde H} \subseteq\omega_H$. The first inclusion uses the nonnegative threefold discrepancy and the rational-point pushforward in \cite[Lemma~6.2]{Tell}. The Gorenstein surface is canonical or elliptic \cite[Proposition~5.26]{Tell}. In the canonical case it is a rational double point, with multiplicity two and embedding dimension three, by the complete-ideal calculation for the fundamental cycle of square $-2$ \cite[\S18]{Lipman}. For an elliptic point the same values hold if its fundamental cycle has negative square at most two \cite[Theorem~5.23]{Tell}.

In the remaining elliptic case the surface point blowup is the normal canonical model and its relative canonical divisor is the negative of its point-ideal divisor \cite[Theorem~5.19 and Corollary~5.20]{Tell}. This blowup is the Cartier strict transform of the generic hyperplane on the threefold point blowup. Indeed, on a chart with generators $x_j/x_i$ the universal transform has equation $c_i+\sum_{j\ne i}c_j x_j/x_i=0$. Eliminating $c_i$ identifies it with the closure of the transform away from the point and the identification persists over the generic coefficient field. The Gorenstein property lifts across the Cartier cut. Normality along it follows by Serre's criterion, since a codimension-one nonregular locus meeting the cut would also be nonregular on the cut. Adjunction there makes the restriction of the threefold relative canonical divisor linearly trivial on the surface. It is effective and nonzero by terminality, since the cut meets the exceptional divisor in a curve. This contradicts negativity for the surface blowup. The rational-point hypothesis is used throughout over the actual field of the hyperplane.

After adjoining an $r$th root of $s$ to $F$, scale each homogeneous element by the corresponding power of this root. The generic fibre then identifies with the chosen finite cover after field extension. At the rational point, embedding dimension and multiplicity are unchanged by extension of scalars. Thus a nonregular cover has embedding dimension four and multiplicity two. Homogeneous lifts of a cotangent basis give a power series presentation. In the nonregular case its kernel has height one and is generated by a semi-invariant equation of order two.
\end{proof}

We give the comparison of discrepancies used both to identify the graded equations and to treat the noncanonical case. Consider a semi-invariant hypersurface equation $f$ of an isolated canonical cover. A sufficiently high polynomial truncation gives the same completed graded quotient. In fact the Tjurina ideal contains $\mathfrak m^s$ for some $s$. An error of order $N>2s+2$ is a combination of $f$ and its partial derivatives with coefficients of order at least $N-s$ and with the matching characters. A homogeneous coordinate change and multiplication by a unit remove it up to order at least $2(N-s)-1>N$. Iterating gives the isomorphism. We can compare the algebraic germs using common pointed \'etale neighbourhoods \cite[Tag~0CAV]{Stacks}.

Let $w_i$ be the coordinate weights of this polynomial modulo $r$ and $c$ its character. Freeness on the punctured cover gives $\gcd(w_i,w_j,r)=1$ for each pair of distinct coordinates: a plane fixed by a nontrivial subgroup would meet the hypersurface along at least a curve. It also gives a nonzero restriction to any coordinate axis of nonunit weight. The canonical character $\sum_i w_i-c$ is primitive modulo $r$. Indeed, for the canonical algebra $\bigoplus_i M_i$, finite duality identifies its dualizing module with $\operatorname{Hom}(\bigoplus_i M_i,M_1)$. Multiplication and projection to degree one give the graded shift, since $\operatorname{Hom}(M_i,M_1)=M_{1-i}$ by reflexive extension. Hypersurface adjunction computes the same character from the weights and the equation, up to sign.

Over $\C$ take a generic polynomial $f_0$ with the same support. It is Newton nondegenerate, since the monomials on each face give a base point free system on the face torus. It has an isolated normal singularity at the origin. To see isolation, use a Tjurina quotient modulo $\mathfrak m^N$ of length less than $N$ for the original $f$. A rank minor of the defining matrix that is nonzero in characteristic $p$ is also a nonzero polynomial of the coefficients in characteristic zero. The generic jet quotient has length less than $N$, which excludes positive dimension. The pair and axis freeness tests persist, so the complex quotient is isolated. Adjunction and the primitive character identify $f_0=0$ as its canonical cover of index $r$.

Use a regular subdivision subordinate to the Newton faces in the lattice $\Z^4+\Z(w_1,\ldots,w_4)/r$. In characteristic zero the strict transform is smooth and transverse to the toric strata. Its exceptional primes over the origin appear on open ray orbits, with primitive vectors $v$ of positive coordinates whose minimum face has at least two terms. Their log discrepancy is $\sum v_i-\nu_v(f_0)$, where $\nu_v$ is the minimum of the weights of the monomials.

The same face in the original characteristic has a nonunit Laurent equation with a codimension-one component of its zero locus. Dividing by a minimizing monomial gives an invariant Cartier equation on the open ray chart. Downstairs the quotient hypersurface $D$ satisfies $rD=\divi(f^r)$ because on the torus the ratio by a minimizing monomial defines the reduced equation, by strong grading and exactness of the degree-zero summand. Toric adjunction and normalization along a prime $E$ above the component give
\begin{equation}\label{eq:terminalface}
 A_X(E)=
 1+m\bigl(\sum v_i-1-\nu_v(f)\bigr)-
             \operatorname{coeff}_E(\operatorname{Cond}),
\end{equation}
where $m\ge1$ is the multiplicity of the restriction of the ray divisor and $\operatorname{Cond}$ is the effective conductor divisor subtracted from the dualizing divisor of the strict transform. The bracket is the coefficient of that ray in $K_{\widetilde Z}+\widetilde D-\pi^*(K_Z+D)$ on the ambient chart. All coordinates of $v$ are positive, so $E$ is exceptional over the origin. This comparison permits a singular or nonreduced face in characteristic $p$.

\begin{lemma}\label{lem:terminalweights}
For a cover in Lemma~\ref{lem:terminalhyp} with $r>1$, the regular case has weights $(a,-a,1)$ up to a change of generator and permutation, with $(a,r)=1$. In the hypersurface case the weights and equation character are $(a,-a,1,c;c)$, with $(a,r)=1$ and either $c=0$ or $(r,c)=(4,2)$. The equation restricts nontrivially to the last coordinate axis.
\end{lemma}

\begin{proof}
In the regular case the toric discrepancy calculation is independent of the field. In the hypersurface case \eqref{eq:terminalface} and terminality force $\sum v_i-\nu_v(f)>1$ for every ray contributing over $\C$. The support is fixed, so the complex quotient is terminal by the normal crossings criterion on its resolution. Mori's classification \cite[Theorems~12, 23 and~25]{MoriTerm} gives the cotangent weights and equation character of the complex cover, and hence those of the original cover. The axis restriction is nonzero by freeness.
\end{proof}

\begin{lemma}\label{lem:localRRall}
At each non-Gorenstein point $Q$ of a $\Q$-factorial terminal threefold in characteristic $p>5$, the local correction in \eqref{eq:canslope-smooth-rr} is the sum for a cyclic basket. The basket contains the local index, and its sum $\sigma_Q$ of $r-r^{-1}$ satisfies $e_Q\le9\sigma_Q$.
\end{lemma}

\begin{proof}
The regular cover gives the single entry of Lemma~\ref{lem:cyclicrr}. In the hypersurface case use the weights of Lemma~\ref{lem:terminalweights} and put $n=\operatorname{ord}_t f(0,0,0,t)$ with $t$ the last variable. We compare two projective hypersurfaces of the same weighted degree. One has the given germ, and the other replaces it by cyclic points. All other singular points on both are cyclic.

For $c=0$ take weights $(r,r,a,r-a,1)$ with $0<a<r$ and a large degree $D$ divisible by all weights. Along the axis of the two weight-$r$ coordinates, the restriction is an ordinary binary form of degree $D/r$. On the special member prescribe the finite jet of the graded equation at a point of this axis. The restriction has a root of order $n$ there and $D/r-n$ simple roots elsewhere. A general member has $D/r$ simple roots. At each simple root the weight-zero coordinate eliminates on the cover and the quotient has weights $(a,-a,1)$ modulo $r$.

For $(r,c)=(4,2)$ take weights $(4,2,1,3,1)$ and a large degree $D\equiv2\pmod4$ divisible by three. On the $(4,2)$-axis with coordinates $(u,t)$, the restriction is $tP(u,t^2)$ with $P$ a binary form of degree $(D-2)/4$. The weight-four point is always present. On a general member a nonzero linear $t$ term there gives the cyclic chart of index four. The other roots are simple and give cyclic points of index two, since the remaining three coordinate weights are odd. On the special member prescribe the given jet at the weight-four point. This leaves $(n-1)/2$ fewer simple roots elsewhere.

In both constructions the finite jet is prescribed in a weight sheaf, and one can twist its kernel by sufficiently high very ample Cartier powers. This leaves free first jets off the prescribed point wherever the sheaf is invertible, and the axis restrictions are binary forms with the indicated jets. Thus one can impose the simple roots, avoid the other isolated coordinate points of nontrivial stabilizer, and obtain smoothness off the axis. The simple-root calculation uses elimination on the graded cover, including when $p\mid r$. Each hypersurface is integral, since distinct positive-degree factors meet in dimension at least two, and a repeated factor is singular in codimension zero.

The weighted adjunction formula gives the canonical sheaf in degree $D-\sum w_i$. This follows on the smooth locus and extends by reflexivity. The weighted sheaves on each hypersurface are reflexive. Indeed, in the affine quotient charts the character summands are Cohen--Macaulay and of rank one. The graded hypersurface sequence in weighted projective space gives the same Euler characteristics for each weight on both members. In particular all the canonical powers have the same Euler characteristics. A common Cartier multiple of the two canonical divisors gives equality of $K^3$ by its Hilbert polynomial. Evaluating \eqref{eq:canslope-smooth-rr} at a positive Cartier multiple also gives equality of $c_2(T_Z)\cdot(-g^*K)$, since every local term vanishes there.

Use the primary-ideal comparison preceding Lemma~\ref{lem:cyclicrr} to transport a resolution from the given $\Q$-factorial germ to the weighted special member, and choose toric resolutions at the other cyclic points. Subtracting their known local terms from \eqref{eq:canslope-smooth-rr} gives the basket of the prescribed point. It consists of $n$ cyclic entries of index $r$ when $c=0$. When $(r,c)=(4,2)$ it consists of one index-four entry and $(n-1)/2$ index-two entries.

For $c=0$ use the invariant parameter ideal $(x^r+y^r,z^r,t-\lambda xy)$. After $t=\lambda xy$ and $z=0$, the equation has order at most $2n$ on each line $x=\xi y$ in $x^r+y^r=0$ for general $\lambda$. Indeed the coefficient of $y^{2n}$ has a nonzero term of degree $n$ in $\lambda$ from $t^n$. Counting lines with multiplicity, the resulting plane intersection has length at most $2nr$. Replacing $z$ by $z^r$ multiplies the colength by $r$, since $z$ is a parameter in the one-dimensional complete intersection before imposing $z=0$. The colength on the cover is thus at most $2nr^2$. Multiplicity for the parameter ideal scales by the finite cover rank $r$, giving $e_Q\le2nr\le(8/3)\sigma_Q$. For index four, the invariant ideal of fourth powers of the four variables extends to a reduction of the fourth power of the maximal ideal of the cover. Its multiplicity upstairs is $2\cdot4^3$, since the cover has multiplicity two. Again dividing by the rank, we get $e_Q\le32\le9\sigma_Q$.
\end{proof}

\subsection{The discrepancy gap and bounded complements}\label{s:mldgap}
On a noncanonical klt threefold some exceptional divisor has log discrepancy less than one. For Proposition~\ref{pr:noncan} we need such a discrepancy uniformly bounded away from one. We compare the local equation with a complex terminal equation and treat the remaining cyclic hypersurface directly.

\begin{lemma}\label{lem:splitcar}
In any characteristic, a normal isolated quotient germ
\[
 \left(k[[x,y,z,t]]/(xy+g(z^r,t))\right)^{\mu_r(a,-a,1,0)},
 \qquad (a,r)=1 ,
\]
is terminal if the hypersurface is isolated and normal and $g(0,t)\ne0$.
\end{lemma}

\begin{proof}
Graded finite determinacy permits a polynomial truncation of $g$, preserving the form of the equation. Put $U=z^r$, $G=xy=-g(U,t)$, and let $b$ be the inverse of $a$ modulo $r$ (take $b=0$ if $r=1$). On the quotient chart where $x^r$ is invertible use coordinates $A=x^r$, $Z'=zx^{-b}$ and $t$. Then $U=A^b(Z')^r$. The rational $r$-canonical form
\[
 \zeta=U\bigl(d\log A\wedge d\log Z'\wedge dt\bigr)^{\otimes r}
\]
is an invertible generator there, since it equals $A^{b-r}(dA\wedge dZ'\wedge dt)^{\otimes r}$.

On the chart where $y^r$ is invertible use $B=y^r$ and $Z''=zy^b$. Then $AB=G^r$, $Z''=Z'G^b$ and $U=B^{-b}(Z'')^r$. Modulo $dt$ one has $d\log G=\lambda(b\,d\log A+r\,d\log Z')$ for a rational function $\lambda$. The change of the two logarithmic differentials has matrix
\[
 \begin{pmatrix}-1+rb\lambda&r^2\lambda\\
                 b^2\lambda&1+rb\lambda\end{pmatrix}
\]
of determinant $-1$. Thus the local generators agree up to sign on the overlap. The complement $x=y=0$ has codimension at least two, so $\zeta$ generates $\omega_X^{[r]}$ on the germ. All these differentials are on the quotient.

For a prime $E$ on a regular model over the origin, the wedge of the three logarithmic differentials of $A,Z',t$ has at worst a simple pole, since all polar parts are multiples of the same uniformizer differential. Using $dt=t\,d\log t$ gives
\[
 A_X(E)=1+\operatorname{ord}_E(\zeta)/r
 \ \ge\ \operatorname{ord}_E(t)+\operatorname{ord}_E(U)/r>1 .
\]
The strict inequality uses that both functions vanish at the origin. At each punctured point at least one of $x,y,z$ is invertible by $g(0,t)\ne0$. Its weight generates $\Z/r\Z$, so the cover there is a flat torsor. Regularity of the cover descends by faithful flatness.
\end{proof}

\begin{proposition}\label{pr:positivegap}
In characteristic $p>5$, if a klt threefold has an exceptional log discrepancy less than one, it has an exceptional log discrepancy at most $1-1/13$.
\end{proposition}

\begin{proof}
Put $\delta=1/13$ from the complex gap theorem \cite[Theorem~1.4]{LX}. Suppose all log discrepancies are greater than $1-\delta$. On a log resolution the coefficients of the crepant trace are less than $1/2$. At an exceptional centre of codimension at least two on the resolution, the normal crossings lower bound is the sum of the log discrepancies of the parameters through the centre, and exceeds one. Hence all the exceptional places of discrepancy one appear on this resolution. Extract exactly these places on a $\Q$-factorial model by \cite[Lemma~7.7]{Bpos}. This extraction is crepant and leaves every place of discrepancy less than one exceptional. There are no codimension-two singularities on the new model. At such a point the actual transverse surface would either be noncanonical, with discrepancy at most $2/3$ by the surface calculation of Proposition~\ref{pr:noncan}, or be canonical and nonregular, with a crepant exceptional divisor.

We now work near an isolated point on the extracted model. Its torus algebra is again terminal. In \eqref{eq:torusdiscrepancy}, if the sum of the $v_i$ is at least two, the discrepancy exceeds one because all base discrepancies exceed $1/2$. If the sum is one, the base discrepancy is integral and, for an exceptional ray, different from one. The remaining arguments of Lemma~\ref{lem:terminalhyp} give the graded cover over the actual field, with the same free punctured action. A regular cover gives a toric comparison directly.

For a singular cover make the same Newton comparison, preserving the rank of the quadratic part. One can first pair complementary character spaces and diagonalize self-complementary blocks, since $p\ne2$. Varying the coefficients of the resulting fixed quadratic blocks preserves their rank. The genericity conditions for $f_0$ still hold on this support. On the complex toric resolution let $A_0=\sum v_i-\nu_v(f_0)$ for an exceptional prime. If $A_0<1$, \eqref{eq:terminalface} gives an original discrepancy at most $A_0$, so $A_0>1-\delta$. If $A_0=1$, it gives an original discrepancy that is an integer at most one. Positivity forces equality, excluded by the extraction. All remaining exceptional valuations over the complex resolution have log discrepancy greater than one by the normal crossings bound. The complex gap theorem makes this quotient canonical, and the exclusion of equality makes it terminal.

Consider now any isolated point at which a place of discrepancy below one persists. Its local index satisfies $r\ge5$, since for $r\le4$ such a discrepancy would be at most $3/4$. The regular cover case here is terminal by the characteristic-independent lattice calculation. Otherwise Mori's classification supplies the weights $(a,-a,1,0)$ with invariant equation. Its rank-one quadratic cases have index at most four, so the original quadratic part has rank at least two. Only the weight-zero variable can have a square in an invariant quadratic part. Hence there is a nondegenerate cross term of opposite unit characters. Splitting off this block by successive graded coordinate changes writes the equation as $xy+g(z^r,t)$, with $g(0,t)\ne0$ by freeness. Lemma~\ref{lem:splitcar} makes the original chart terminal, a contradiction.
\end{proof}

\begin{theorem}[Exceptional Fano threefolds with zero boundary]\label{th:allexceptionalbare}
There are absolute positive integers $P,N$ such that every exceptional Fano threefold with zero boundary over an algebraically closed field of characteristic $p>P$ has a klt $N$-complement. The class is arithmetically bounded in sufficiently large characteristic.
\end{theorem}

\begin{proof}
For noncanonical $X$, apply Proposition~\ref{pr:positivegap} and Theorem~\ref{th:nonempty} with $\gamma=1/13$. Suppose $X$ is canonical. If a smooth projective resolution has a nonzero regular two-form, use Theorem~\ref{th:twoformscaled}. Otherwise Lemma~\ref{lem:chieuler} gives $\chi(\cO_X)=1$.

The crepant extraction of the exceptional discrepancy-one places, as in the proof of Proposition~\ref{pr:positivegap}, gives a projective $\Q$-factorial terminalization of $X$. At its non-Gorenstein points Lemma~\ref{lem:localRRall} gives the cyclic baskets and the multiplicity bound with $M=9$. Proposition~\ref{pr:basketslope} now bounds the first effective anticanonical degree. A common multiple of these degrees and of the indices in the other cases gives $N$. An effective member in degree $N$ gives a klt complement by exceptionality. Corollary~\ref{cor:bare} gives arithmetic boundedness.
\end{proof}

\begin{proof}[Proof of Theorem~\ref{th:intcomp}]
For a non-exceptional pair apply \cref{th:nonexceptional}. For an exceptional pair with $B\ne0$, apply Theorem~\ref{th:nonempty}. These two results include the cases of non-big $L$.

It remains to treat an exceptional pair with $B=0$. Since $X$ is of Fano type and $-K_X$ is nef, it is big and semiample. Write $f:X\to Y$ for its ample contraction. Then $-K_Y$ is ample and $K_X=f^*K_Y$ with compatible choices of canonical divisors. Every effective anticanonical completion on $Y$ pulls back to one on $X$, so $Y$ is exceptional by crepancy. Pulling back the principal log divisor of a complement from Theorem~\ref{th:allexceptionalbare} gives a klt complement of the same index on $X$. For fixed $R$, take the largest characteristic threshold and the least common multiple of the indices from the three cases. A multiple retains the same boundary and satisfies \eqref{eq:definition}.
\end{proof}

\section{Boundedness of weak log Fano threefolds}\label{s:bab}
The anticanonical ring associates an $\epsilon$-lc Fano model to each $\epsilon$-lc weak log Fano pair. We use the family constructions of Subsection~\ref{s:arithmetic} to recover the possible weak models, and the estimates of Section~\ref{s:degrees} to choose one characteristic bound for all $\epsilon$.

\subsection{Recovering weak models from their anticanonical models}\label{s:anticanonical}
The boundary $C$ may have real coefficients. Since $K_X+C$ can be $\R$-Cartier even when $K_X$ has no Cartier multiple, we form the anticanonical section ring using integral Weil multiples of $-K_X$.

\begin{proposition}\label{pr:amplemodel}
In characteristic $p>5$, let $(X,C)$ be an $\epsilon$-lc weak log Fano threefold as in Theorem~\ref{th:intbab}. A divisible Veronese subring of its anticanonical ring defines a birational model $Y$ with $-K_Y$ ample and $\Q$-Cartier. For every divisorial valuation $v$,
\begin{equation}\label{eq:discrep}
                           A_{Y,0}(v)\ge A_{X,C}(v).
\end{equation}
Every prime of $X$ exceptional over $Y$ has log discrepancy at most one over $Y$. When $C=0$, the map $X\to Y$ is a crepant morphism.
\end{proposition}
\begin{proof}
Take a small projective $\Q$-factorialization $X'\to X$ of the klt pair \cite[Theorem~1.6]{Bpos}. Write $C'$ for the strict transform and put $L'=-(K_{X'}+C')$. Smallness gives $K_{X'}+C'=f^*(K_X+C)$ for $f:X'\to X$. Write $L'\simr A+G$ with $A$ ample and $G$ effective. For small $u>0$ the boundary $C'+uG$ is klt and its negative log divisor is equivalent to $(1-u)L'+uA$, which is ample. Rational approximation on its fixed support gives a rational log Fano boundary $\Theta$ on $X'$.

Put $D=-K_{X'}=L'+C'$. This is a big rational Cartier divisor. For small positive rational $t$, the class $tD-(K_{X'}+\Theta)$ is ample. Choose its representative so that the new boundary is klt \cite[Theorem~1]{Tan}. The resulting log MMP is a $D$-MMP. After the last semiample contraction the ample divisor on the target is $-K_Y$ up to rational linear equivalence, by pushing forward. On a common resolution with maps $u:U\to X'$ and $v:U\to Y$, negativity gives
\begin{equation}\label{eq:fixed}
       u^*D=v^*(-K_Y)+F,\qquad F\ge0,\quad F\text{ exceptional over }Y.
\end{equation}
Consequently divisible Veronese subrings agree. Smallness of $X'\to X$ identifies the ring on $X'$ with the Weil section ring on $X$.

Every effective representative $G'\simr u^*D$ contains $F$, since $G'=F+v^*v_*G'$. Indeed this follows from \eqref{eq:fixed} by writing out the real principal difference and pushing it to $Y$. Pullbacks of effective $\R$-Cartier divisors are effective, as can also be seen by rational approximation in the finite rational cone cut out by the coefficient inequalities. At each valuation, the infimum of orders among effective representatives of $L'$ is zero. To see this, write $L'\simr ((1-h)L'+hA)+hG$ for $0<h<1$. The first summand is ample and has an effective representative avoiding the centre. Taking $h$ to zero in $u^*(C'+L'_{\mathrm{eff}})\ge F$ gives $u^*C'\ge F$. The identity
\[
 K_U-v^*K_Y
   =K_U-u^*(K_{X'}+C')+u^*C'-F
\]
now proves \eqref{eq:discrep} on every higher model. On a prime of $X'$ contracted over $Y$, its log discrepancy for $(Y,0)$ is one minus its coefficient in $F$. Finally, if $C=0$, semiampleness of $-K_X$ gives the morphism to $Y$ directly and $F=0$.
\end{proof}

\begin{theorem}[Recovering the original weak model]\label{th:weak}
Fix \(\eps>0\). Arithmetic boundedness of the \(\eps\)-lc Fano threefolds with zero boundary in sufficiently large characteristic implies arithmetic boundedness, in sufficiently large characteristic, of the threefolds with an \(\eps\)-lc weak log Fano witness. The horizontal construction can raise the threshold depending on the chosen family and on \(\eps\). In any one fixed characteristic above the bound for Cartier descent of Lemma~\ref{lem:index}, the same implication holds for bounded families over that characteristic.
\end{theorem}

\begin{proof}
Proposition~\ref{pr:amplemodel} associates to each \(X\) an anticanonical model \(Y\) in the bounded class. Use Lemma~\ref{lem:index} with empty boundary. On each marked resolution, the exceptional primes over \(Y\) of log discrepancy at most one consist of primes on that resolution and the additional stratum valuations from Lemma~\ref{lem:toroidal}. The coefficients and the dimensions of the images are constant, so the list is obtained on the geometric generic fibre and spread.

Extract exactly these primes on the geometric generic fibre with a \(\Q\)-factorial source \(T_\eta\), by \cite[Corollary~1.4.3]{BCHM} in characteristic zero and \cite[Lemma~7.7]{Bpos} in positive characteristic. The crepant boundary is effective with coefficients at most \(1-\eps\), and its negative log divisor is the pullback of the ample \(-K_{Y_\eta}\). Kodaira perturbation gives a log Fano witness as in Proposition~\ref{pr:amplemodel}. Spread the extraction and a common resolution with the labelled stratum blowups, keeping the images of the specified primes divisorial on the source. Every prime divisor on \(X\) then occurs on \(T_s\), by the discrepancy comparison in Proposition~\ref{pr:amplemodel}. Apply Lemma~\ref{lem:chambers} to \(T_s\dashrightarrow X\). This bounds the original weak models, including those with only a Weil canonical divisor. On closed remainders perform the same construction starting from the marked resolutions of their geometric generic fibres.
\end{proof}

\subsection{An absolute characteristic bound}\label{s:vertical}
\begin{theorem}\label{th:nonexbab}
For each fixed $\epsilon>0$, non-exceptional $\epsilon$-lc Fano threefolds have klt complements of common index and are arithmetically bounded in sufficiently large characteristic. The arithmetic boundedness extends to the original weak log Fano models with these ample models.
\end{theorem}
\begin{proof}
Use the lc index $N$ from \cref{th:nonexceptional} with $B=0$, and the construction $(Y,B,R)$ in Proposition~\ref{pr:prepare}. Write $d$ for the first birational degree of $R$ and $V_R=R^3$. Apply Theorem~\ref{th:jetsurface} with $v=16$. If $d>Kl$, it gives a covering family with $R^2S/V_R<l/8$. The divisors in the family are Weil linearly equivalent by the Picard argument in Proposition~\ref{pr:moving}, since a smooth resolution of $Y$ is rationally chain connected \cite[Theorems~1.5 and~4.8]{GNT}. Two distinct members give the pencil of Theorem~\ref{th:smallpencil}, hence $d\le C_N(1+l/8)$. Either bound on $d$ gives
\[
 d^3V_R\le C'_N\max(1,V_R).
\]
The right hand side is bounded by Theorem~\ref{th:fullupper} since $V_R=q^3(-K_X)^3$ with $q=\ceil{1/\epsilon}$. Proposition~\ref{pr:prepare} gives a birational anticanonical system of degree $b=Nqd$ with $b^3(-K_X)^3=N^3d^3V_R$ bounded. Theorem~\ref{th:birvol} gives the klt complement and the arithmetic bound, and Theorem~\ref{th:weak} recovers the original varieties.
\end{proof}

The numerical estimates hold above an absolute characteristic bound. For fixed $\eps$, spreading from characteristic zero can omit finitely many further primes. At each such prime we construct the missing families over geometric generic points of the vertical strata.

\begin{proof}[Proof of Theorem~\ref{th:intbab}]
Choose \(P_*\) to include \(5\), the absolute bound in \cite[Corollary~4.4]{Bern}, and the characteristic bound for the klt index of the exceptional Fano threefolds of Theorem~\ref{th:allexceptionalbare}. Include also the bound for the empty-boundary index \(N\) of \cref{th:nonexceptional}, the characteristic bound in Theorem~\ref{th:jetsurface}\textup{(i)} for that \(N\), and \(\max(5,N)\) for Theorem~\ref{th:smallpencil}. The argument of Theorem~\ref{th:fullupper} uses the same bound. These choices are independent of \(\eps\).

Fix \(0<\eps\le1\). By Proposition~\ref{pr:amplemodel}, the ample model \(Y\) of \(X\) gives an \(\eps\)-lc pair \((Y,0)\). In the exceptional case this model has the absolute klt index and volume at most \(27\). In the non-exceptional case, the proof of Theorem~\ref{th:nonexbab} gives an lc \(N\)-complement and the bound \eqref{eq:birvol} with \(U=U(\eps)\), for every \(p>P_*\). In particular \(q=\lceil1/\eps\rceil\) in that proof changes the volume bound and the degrees, while the jet and pencil thresholds use only \(N\).

The exceptional ample models are bounded by Theorem~\ref{th:allexceptionalbare}. Theorem~\ref{th:birvol} bounds the non-exceptional ample models. Apply Theorem~\ref{th:weak} to the union to recover the original varieties \(X\). For every prime \(p>P_*\) omitted in the horizontal constructions, use the same numerical bounds on the parameter spaces over \(\mathbb F_p\). The marked resolutions are constructed on their geometric generic threefolds. The Frobenius estimates apply there and on the smooth fibres after spreading. The klt indices in Theorem~\ref{th:partner} are actual principal multiples on each geometric generic fibre, and the finite descent and spreading preserve these relations. Lemmas~\ref{lem:index} and~\ref{lem:chambers}, including the fixed-characteristic Picard and contraction arguments, then give the families of weak models over that prime.

For this fixed \(\eps\) only finitely many primes need these extra families. Their disjoint union with the horizontal families is projective over a finite-type \(\Z[1/(P_* !)]\)-scheme, after restricting the bases to that ring. When \(\eps>1\) the class is empty, since the log discrepancy at any prime on \(X\) is \(1-\operatorname{coeff}_P C\le1\).
\end{proof}

\section{Kawamata--Viehweg vanishing}\label{s:kv}

\subsection{Reduction by the adjoint minimal model program}\label{s:kv-adjoint}
Pulling back an integral Weil divisor can introduce fractional coefficients. Rounding on a common resolution gives a single line bundle that computes cohomology on both birational models.

On birational models we use the same function $u$ and the same denominator $m$ when writing the rational principal divisor $\divi(u)/m$.

\begin{lemma}\label{kv:transport}
Let \(v:U\to T\) be a log resolution of a three-dimensional klt pair \((T,C)\) in characteristic \(p>5\). Let \(A\) be a rational Cartier divisor on \(T\), with \(C=\ceil A-A\). Put \(D=K_T+\ceil A\). Then
\[
 Rv_*\cO_U(K_U+\ceil{v^*A})=\cO_T(D).
\]
\end{lemma}
\begin{proof}
Write \(C_U=v^*(K_T+C)-K_U\), and set \(Q=\ceil{v^*A}-v^*A\). The divisor upstairs is relatively \(\Q\)-linearly equivalent to \(K_U+Q\), all coefficients of \(Q\) are below one, and \(v_*Q=C\). Ordinary Grauert--Riemenschneider vanishing applies \cite[Definition~1 and Theorem~3]{BK}. The difference from \(v^*D\) is \(Q-C_U\), exceptional with coefficients greater than \(-1\). A section of \(\cO_T(D)\) pulls back to an effective rational divisor after its principal divisor is added to \(v^*D\). Adding \(Q-C_U\) gives an integral divisor with coefficients greater than \(-1\), hence effective. The converse follows by pushforward.
\end{proof}

Two applications of the formula will recur. If \(D=K_T+C+\divi(u)/m\) with \(D\) integral and \(C\) effective klt, take \(A=\divi(u)/m\). On a second model on which the corresponding correction is klt, a common resolution computes both cohomologies. A pullback from a common base can be added to \(A\). For a rounded-down divisor \(J=\floor A\), apply the formula to \(K_T-J\) with correction \(\{A\}=A-\floor A\), when \((T,\{A\})\) is klt and \(T\) is \(\Q\)-factorial. Proper duality gives the analogous comparison for \(\cO(J)\), also for derived images on a common base. Indeed the divisorial sheaves on these klt threefolds are Cohen--Macaulay \cite[Corollary~1.3]{ABL}. For a \(\Q\)-Cartier Weil divisor on a non-\(\Q\)-factorial klt model the local assertion follows from \cite[Theorems~3 and 19(2)]{BK}, since a Cartier multiple is locally principal. Their duals into the dualizing complex are their canonical dual sheaves shifted by three. Over a singular base one uses its dualizing complex in the same duality formula.

\begin{lemma}\label{kv:relative}
Let \(f:T\to Z\) be a contraction from a projective \(\Q\)-factorial threefold of Fano type to a positive-dimensional normal variety, and let \(L\) be an integral Weil divisor. Then \(R^i f_*\cO_T(L)=0\) for \(i>0\) in either of the following cases.
\begin{enumerate}[label=\textup{(\roman*)}]
\item \(L\equiv_Z0\).
\item Both \(-L\) and \(L-K_T\) are \(f\)-ample.
\end{enumerate}
If the general fibre has positive dimension, the ordinary direct image also vanishes in case \textup{(ii)}.
\end{lemma}
\begin{proof}
Work near a closed \(z\in Z\). In (i) some multiple of \(L\) is pulled back from \(Z\) by relative base point freeness. Its semiample contraction is constant on the fibres and its Stein factor is finite birational over the normal base. Thus locally \(A=\divi(u)/m=L\). Choose a klt log Fano witness for \(f\). In (ii) take a diluted general effective \(\Delta'\sim_{\Q,Z}-L\) retaining klt and avoiding the divisorial components of \(T_z\). Now \(-(K_T+\Delta')\) is relatively ample. Shrinking the base, write \(A=\divi(u)/m=L+\Delta'\). For (i) the fractional correction of \(A\) is zero.

Use the fibre model in \cite[Proposition~2.15]{GNT} for the chosen log Fano pair. On this model \(g:Y\to Z\) the reduced fibre above \(z\) is a prime \(S\), \(-S\) is relatively nef, and there is a plt boundary \(\Lambda\) with floor \(S\) and \(-(K_Y+\Lambda)\) relatively ample. Step~1 of the construction contracts the exceptional primes outside \(T_z\), the later steps extract none, and only \(S\) remains over \(z\).

In (ii) we need \(\Lambda-S\) to contain the strict transform \((\Delta')^s\). The effective sum at the end of Step~2 is relatively trivial, is supported above \(z\), and omits \(S\). The function in its principal relation is constant on the generic fibre, so the sum is a pullback. A nonzero effective pullback supported over \(z\) would contain \(S\). Hence the sum vanishes and Step~3 preserves a relatively trivial plt log divisor \(K_Y+S+(\Delta')^s+P\). Here \(P\) comes from a relatively ample divisor \(P_V\) on the resolution. For a general ample Cartier divisor \(H_Y\), the class \(P_V-\eta H_Y^s\) is still relatively ample for small \(\eta>0\). Mix \(P\) with the transform of a diluted member of this class. This preserves \((\Delta')^s\) and plt and makes the negative log divisor relatively ample.

Write \(Q=\{A_Y\}\), \(b=\coeff_S Q\), \(J=\floor{A_Y}\). In (i) \(Q=bS\), and in (ii) \(Q=(\Delta')^s+bS\), where the strict trace away from \(S\) is meant. Our choice of \(\Delta'\) avoids \(S\) if it was already a divisor on \(T\). Thus \((Y,Q)\) is klt and \(S\le\Lambda-Q+bS\le\Lambda\). We have
\[
 J-(K_Y+\Lambda-Q+bS)\sim_{\Q,Z}-(K_Y+\Lambda)-bS ,
\]
which is relatively ample. The plt vanishing theorem \cite[Theorem~1.1]{BernVan} gives the required higher direct image vanishing for \(J\). Dual rounded comparison for the single rational principal \(A\) on a common resolution identifies its derived image with that of \(L\) on this neighbourhood. The last assertion follows from negativity on the generic fibre and torsion freeness of the direct image.
\end{proof}

To run the MMP for the divisor in Theorem~\ref{th:kv}, write \(A_0=D-K_X-\Delta\simq H+E\) with \(H\) ample and \(E\ge0\). For small rational \(t>0\), add \(tE\) to \(\Delta\) and take a sufficiently diluted effective representative of the ample class \((1-t)A_0+tH\). This gives an effective big klt boundary \(C\) such that \(D\simq K_X+C\) \cite[Theorem~1]{Tan}. Take a small \(\Q\)-factorialization and run the klt adjoint MMP. Each model is of Fano type, by pushing a big klt Calabi--Yau completion. The transforms of \(C\) stay klt and big. \Cref{kv:transport} applied to \(D-K_X-C\) preserves all cohomology through these steps, including the small modification. The top group already vanishes on \(X\) by duality, since the negative of the big effective class \(D-K_X\) has no section.

At a negative Mori contraction onto a positive-dimensional base, \(D-K_X\) is relatively ample because its big effective representative has positive degree against a moving fibre curve and the relative Picard number is one. Apply \cref{kv:relative}. On a nef and big output, relative base point freeness gives a birational map to the ample model. Its pushforward is an integral \(\Q\)-Cartier ample divisor whose pullback is \(D\). The log divisor \(K+C\) is crepant here, so \cref{kv:transport} again preserves cohomology. The ample target itself need not be \(\Q\)-factorial.

Fix the index \(n=n(\{1\})\) of the zero-boundary complement from \cref{th:intcomp}. Every projective threefold \(X\) of Fano type has an effective lc boundary \(\Sigma\) with \(n(K_X+\Sigma)\sim0\). To see this, work first on a small \(\Q\)-factorialization and take the ample model \(Z\) of the big divisor \(-K_X\) by its MMP. This is a Fano variety of Fano type. On a common resolution \(u:U\to X\), \(v:U\to Z\), negativity gives
\begin{equation}\label{kv:eq:fixed}
 u^*(-K_X)=v^*(-K_Z)+E,\qquad E\ge0 .
\end{equation}
If \(K_Z+\Sigma_Z=\divi(s)/n\), its crepant trace on \(X\) is the pushforward of \(v^*\Sigma_Z+E\), hence effective. The same principal relation proves lc and descends from the small modification.

\begin{lemma}\label{kv:bounded}
Fix \(e>0\). In proving \cref{th:kv}, in characteristic sufficiently large in terms of \(e\), any ample integral \(\Q\)-Cartier Weil divisor on a threefold whose bare anticanonical ample model is \(e\)-lc has vanishing higher cohomology. On an \(e\)-lc Fano threefold of Picard number one the same holds for an integral \(\Q\)-Cartier Weil \(J\) with \(J-K_X\) ample.
\end{lemma}
\begin{proof}
We may work over algebraic closures of finite fields. A counterexample to \cref{th:kv} would spread over an integral algebra of finite type over \(\mathbb F_p\). In doing so one can make the adjoint difference ample by adding a small multiple of \(E\) in the Kodaira decomposition above to the boundary. Spread the log Fano witness, a log resolution, the rational Cartier multiples and the rounded line bundle in \cref{kv:transport}, with \(A=D-(K_X+\Delta)\) for the perturbed boundary. After shrinking, the resolution is fibrewise birational with the same SNC components and coefficients. The ample classes remain ample. Upper semicontinuity of cohomology on this resolution preserves nonvanishing on a geometric closed fibre. For the fixed-model assertion one can include the anticanonical maps and \eqref{kv:eq:fixed} among the spread data, keeping the crepant coefficients of a log resolution of \(Z\) and the dimensions of exceptional images.

By \cref{th:intbab} the \(e\)-lc Fano ample models belong to a family of finite type over \(\Z\). Use \cref{lem:index} to stratify the closure of the selected points, with klt Fano geometric generic fibre on every horizontal stratum. After a finite cover of such a stratum, the reduction theorem for characteristic-zero log Fano varieties gives global \(F\)-regularity at its closed points on a dense open \cite[Theorem~5.1]{SchwedeSmith}. Noetherian induction and omission of the finitely many vertical residue characteristics give one bound depending on \(e\). This is the use of arithmetic boundedness in \cite[Proposition~3.2 and Remark~3.3]{CTW}. Splittings extend to algebraic closures of finite fields. For a splitting along a divisor defined over a finite extension, first dominate that divisor by the pullback of a Cartier divisor over the original field and then extend the splitting.

For every effective Cartier \(H\), global \(F\)-regularity gives a splitting of \(\cO\to F_*^j\cO(H)\) for some \(j>0\). In particular \(Z\) is Frobenius split. The associated rational trace map has value one on \(1\). On a normal variety it is regular precisely when its associated anticanonical divisor is effective \cite[\S4]{SchwedeSmith}. By \eqref{kv:eq:fixed} that divisor on a small model of \(X\) is the trace of its pullback from \(Z\) plus \(E\). Thus the same map splits \(X\). Twisting on the regular locus and extending reflexively gives injections from \(H^i(X,\cO_X(L))\) into \(H^i(X,\cO_X(p^{aj}L))\). For ample \(L\), Serre vanishing in the finitely many Cartier residue classes kills these targets for large \(a\).

On the bounded rank-one model a negative \(J\) has vanishing in degrees one and two by the same injections and anti-ample Serre vanishing, using Cohen--Macaulay duality in each residue class. If \(J\equiv0\), a multiple is principal by \cref{lem:numerical-torsion}. Choose an effective ample Cartier \(H\) sufficiently positive for all the finitely many sheaves \(\cO(jJ)\) modulo this multiple. A splitting along \(H\) injects the required groups into those of \(\cO(p^jJ+H)\). Finally \(H^3(\cO(J))\) is dual to \(H^0(\cO(K_X-J))=0\) by the assumed positivity.
\end{proof}

On the generic fibre an integral divisor may represent a nontrivial torsion class. The next lemma treats this possibility when computing the direct image.
\begin{lemma}\label{kv:vertical}
Let \(g:T\to Z\) be a Mori fibre contraction with positive-dimensional base, \(T\) a \(\Q\)-factorial threefold of Fano type, and suppose \(n(K_T+\Sigma)\sim0\) with \((T,\Sigma)\) lc. If \(Z\) is a curve, assume the underlying generic surface is an \(e\)-lc Fano surface for a fixed \(e>0\). Let \(L\) be integral, and let \(0\le\Gamma\le \alpha\Sigma\), \(0\le\alpha<1\), be vertical. If \(L+\Gamma\simq g^*M\) for a nef rational Cartier divisor \(M\), then \(H^{>0}(T,\cO_T(L))=0\) in characteristic sufficiently large in terms of \(n,e\).
\end{lemma}
\begin{proof}
Every vertical prime \(V\) is numerically trivial over \(Z\). Indeed its degree on a general fibre curve is zero and \(\rho(T/Z)=1\). Relative base point freeness shows that it is rationally linearly a pullback. Once a vertical representative has been chosen this is an actual pullback after a principal adjustment on \(Z\), since a function with vertical divisor is constant on the proper generic fibre, whose field of constants is \(k(Z)\). In particular no prime lies over a codimension-two subset of \(Z\), and \(\Gamma=g^*\Gamma_Z\). By \cref{kv:relative} the higher direct images of \(\cO_T(L)\) vanish. Its generic restriction is torsion in the Weil class group. A nontrivial torsion class has no section, since an effective nonzero divisor has positive degree against an ample complete intersection. In that case its torsion-free direct image is zero. Otherwise an integral principal adjustment makes \(L=g^*N\) as a rational divisor, with \(N+\Gamma_Z\simq M\). Sections on the generic fibre are constants and effectivity is tested at prime divisors. Consequently \(g_*\cO_T(L)=\cO_Z(\floor N)\).

Suppose first that \(Z\) is a surface. It is of Fano type by \cite[Corollary~5.12\textup{(1)}]{Fib}, applied to \(\Sigma\) with its fixed coefficient set. Choose a small ample \(H\) on \(Z\) and a klt effective \(U\) on \(T\) with \(K_T+U\simq-g^*H\), by completing a log Fano witness. Use \(\Xi=(1-t)\Sigma+tU\) with rational \(0<t<\min\{1-\alpha,1/2\}\). The generic curve is a smooth genus-zero curve in characteristic \(p>5\). The marks in \(\Sigma\) have total unweighted geometric degree at most \(2n\), so for \(p>2n\) they split over a prime-to-\(p\) Galois extension. At such separable marks the klt coefficients of \(\Xi\) are below one. The remaining cycle has geometric weighted degree at most \(2t<1\). Thus the geometric generic pair is klt and the threshold bundle formula \cite[Theorem~5.9]{Fib} gives
\[
 K_Z+B_Z+\mathbf M_Z\simq-tH
\]
with b-nef b-\(\Q\)-Cartier moduli divisor. The generalized pair is klt since the thresholds are positive over every prime on every base model. For a prime \(P\) on \(Z\), a component \(V\) above it of multiplicity \(m\) has \(mN_P\in\Z\). Testing the threshold at \(V\) gives
\[
 (B_Z)_P\ \ge\ 1-\frac{1}{m}+\frac{\Xi_V}{m}
       \ \ge\ \{N_P\}+(\Gamma_Z)_P .
\]
Subtract \(\{N\}+\Gamma_Z\) from \(B_Z\). The boundary is still effective and the generalized pair is klt, since the base surface is \(\Q\)-factorial. Approximate the remaining boundary plus the nef part and \(tH/2\) by an ordinary effective klt boundary, using \cite[Lemma~5.10]{Fib}. Its adjoint difference from \(\floor N\) is rationally equivalent to \(N+\Gamma_Z+tH/2\), which is ample. The surface vanishing theorem \cite[Theorem~1.1]{ABL} finishes this case.

If \(Z\) is a curve, it is \(\PP^1\). Apply the threshold formula to \(\Sigma\) itself \cite[Theorem~5.11\textup{(ii)}]{Fib}. Its assumptions use the fixed \(e,n\), and \(\Sigma\simq-K_T\) is big. Nefness of the moduli part gives \(\deg B_Z\le2\). The same inequality with \(\Sigma_V\) in place of \(\Xi_V\) gives \(B_Z\ge\{N\}+\Gamma_Z\), strictly at any point where \(\Gamma_Z>0\). Thus \(\deg\floor N\ge-2\), with strict inequality unless \(\Gamma_Z=0\), \(\deg N=0\), and every inequality is an equality.

In this last case all multiplicities over a given point \(z\) equal \(m_z\), with \(\{N_z\}=1-1/m_z\). The latter fractions sum to two. The entries \(m_z>1\) are therefore \((2,2,2,2),(3,3,3),(2,4,4)\) or \((2,3,6)\). If \(q\) denotes their least common multiple, \(qN\) is principal. Taking a \(q\)-th root of its defining function gives a connected tame cover of \(\PP^1\), of genus one by Riemann--Hurwitz. Connectedness follows from the minimality of \(q\). The normalized base change to \(T\) is connected since \(k(Z)\) is algebraically closed in \(k(T)\). It is \'etale at every codimension-one point because the orders of the function are divisible by \(q\) on \(T\) and \(q\) is invertible. Tame discrepancy pullback makes this cover of Fano type, so it is rationally chain connected \cite[Theorem~4.8 and Proposition~4.11]{GNT}. Its morphism onto a genus-one curve is impossible. In every case \(\deg\floor N\ge-1\), as required.
\end{proof}

\begin{proposition}\label{kv:small}
In uniformly large characteristic an integral nef \(\Q\)-Cartier Weil divisor of numerical dimension at most two on a threefold of Fano type has zero higher cohomology. The two intermediate groups also vanish at every negative rank-one endpoint of the adjoint MMP above.
\end{proposition}
\begin{proof}
Take a small \(\Q\)-factorialization first. For nef \(D\), its semiampleness permits a klt adjoint representation by adding a diluted representative of \(D\) to a big klt Calabi--Yau completion. Thus \cref{kv:transport} applies across this small map. Write \(D\simq f^*A\) for the semiample contraction \(f:X\to V\), with \(A\) ample if \(V\) has positive dimension. For numerical dimension zero take \(V\) a point and \(A=0\). In each case fix \(u,m\) with \(D=f^*A+\divi(u)/m\).

Take a crepant dlt modification of \((X,\Sigma)\). It is of Fano type. For this, mix its crepant boundary with the trace of a big klt Calabi--Yau completion from \(X\). A small weight on the latter gives effective crepant coefficients at the extracted lc places and retains a big movable summand. Kodaira perturbation gives a log Fano witness.

On this dlt model write the full crepant boundary as \(B\). The pair with boundary \(B/2\) is \(1/(2n)\)-lc. Indeed a positive log discrepancy for \(B\) is at least \(1/n\). At a zero discrepancy the centre meets the SNC locus of the dlt pair, where the log discrepancy for the bare variety is at least one.

The log divisor \(K+B/2\simq-B/2\) is anti-big on the generic fibre over \(V\), since \(B\simq-K\) is big. Its MMP over \(V\) therefore ends with a Mori fibre space \(g:T\to Z\). The full boundary, now denoted \(\Sigma_T\), stays crepant, while the half-boundary MMP preserves the discrepancy lower bound. In particular \(-K_T\simq\Sigma_T\) is \(g\)-ample. If \(Z\) is a curve the generic surface is Fano and \(1/(2n)\)-lc: localization at the generic point preserves the discrepancy inequalities.

On all these models use the same rational divisor \(P=f^*A+\divi(u)/m\), using the morphism to \(V\) on each model. Its fractional part on the dlt model is supported on new floor components, and each of its coefficients is less than one. On \(T\) its trace \(\Gamma=\{P_T\}\) satisfies \(\Gamma\le\alpha\Sigma_T\) for some \(\alpha<1\). This inequality holds on each intermediate model and makes the correction klt, by convexity with the klt underlying variety. Dual rounded comparison from \cref{kv:transport} identifies the cohomology with that of \(L=\floor{P_T}\). Here \(L+\Gamma\) is rationally equivalent to the pullback of a nef divisor on \(Z\). If \(Z\) is a point, \(L-K_T\simq\Sigma_T-\Gamma\) is ample, so \cref{kv:bounded} applies. If \(\Gamma\) has a horizontal component over a positive-dimensional \(Z\), then \(-L\) and \(L-K_T\) are \(g\)-ample by relative Picard number one. Apply \cref{kv:relative}. For a vertical \(\Gamma\) use \cref{kv:vertical}. This proves the first assertion, including \(H^{>0}(X,\cO_X)=0\).

At a negative rank-one endpoint write \(D\simq d(-K_X)\). Here \(-1<d<0\), since \(D-K_X\) is big and \(D\) is anti-ample. The rational linear equivalence follows from \cref{lem:numerical-torsion}. Put \(c=1+d\). This time take \(P=\divi(u)/m=D-K_X-c\Sigma\) and transport \(K+\ceil P\) through the same half-boundary MMP over a point. Its correction is \(c\Sigma\) on the initial model. On the dlt model the only new primes are in the floor and their corrections are below one. Thus on \(T\) again \(\Gamma=\ceil{P_T}-P_T\le\alpha\Sigma_T\) with \(\alpha<1\). For \(J=K_T+\ceil{P_T}\) we have \(-J\simq\Sigma_T-\Gamma\) and \(J-K_T\simq\Gamma\). If \(Z\) is a point and \(\Gamma\ne0\), use \cref{kv:bounded}. If \(Z\) is a point and \(\Gamma=0\), the canonical dual is numerically trivial. If \(Z\) has positive dimension, use \cref{kv:relative} for a horizontal \(\Gamma\), and use \cref{kv:vertical} on \(K_T-J\) for a vertical \(\Gamma\), with nef class zero. Duality exchanges the two intermediate groups.
\end{proof}

\subsection{Vanishing by restriction to a plt surface}\label{s:kv-restriction}
\begin{lemma}\label{kv:surface}
Fix \(n>0\). Let \((S,\Theta)\) be a projective klt surface pair in sufficiently large characteristic with \(n(K_S+\Theta)\sim0\). Let \(G\) be an integral Weil divisor and let \(0\le C\le\Theta\) be rational. If \(G-(K_S+C)\) is nef and big then \(H^i(S,\cO_S(G))=0\) for \(i>0\).
\end{lemma}
\begin{proof}
The coefficients of \(\Theta\) lie in \(\{0,1/n,\ldots,(n-1)/n\}\). Kawakami's log lifting theorem for klt Calabi--Yau surfaces \cite[Proposition~3.5]{KawakamiLCY} gives a log resolution with the exceptional locus and the strict boundary lifting over \(W(k)\). Put \(A=G-K_S-C\). On this resolution \(\pi\), the fractional support of \(\pi^*A\) lies in the marked divisor. Testing effectivity of a rational Cartier divisor before and after pullback gives \(\pi_*\cO(-\ceil{\pi^*A})=\cO_S(-\ceil A) =\cO_S(K_S-G)\). The logarithmic vanishing theorem for the nef and big divisor \(\pi^*A\), with differential degree zero, kills \(H^j(\cO(-\ceil{\pi^*A}))\) for \(j=0,1\) \cite[Theorem~2.11]{KawBS}. Leray injects the two corresponding groups on \(S\) into these groups. Surface duality proves the claim.
\end{proof}

Let \(h:Y\to X\) be a projective birational morphism from a \(\Q\)-factorial threefold of Fano type and let \((Y,\Omega)\) be plt with prime floor \(S\) and \(K_Y+\Omega\simq0\). Assume that every exceptional prime of \(h\) equals \(S\). Let \(D\) be an ample integral \(\Q\)-Cartier Weil divisor on \(X\), set \(H=h^*D\), and put \(d=\coeff_S\{H\}\). We first identify the cohomology of \(D\) with that of \(H-dS\). Indeed, the pushforward \(\Omega_X=h_*\Omega\) is crepant. Take an effective klt log Calabi--Yau boundary \(C_0\) on \(X\). For small \(\lambda>0\) the crepant transform \(C_Y\) of \(C=(1-\lambda)\Omega_X+\lambda C_0\) is effective klt and \(\coeff_S C_Y>d\). Only at the exceptional \(S\) could the crepant transform of \(C_0\) have a negative coefficient. We have
\[
K_Y+\ceil{h^*(D-(K_X+C))}
=\ceil{H-C_Y}=H-dS.
\]
The two klt boundaries for the common rounded comparison in \cref{kv:transport} are \(C\) and \(C_Y-dS\). In particular this comparison applies when only \(K_X+C\) is \(\Q\)-Cartier.

First run the \(-S\)-MMP over \(X\) with ample scaling at parameter zero. It ends birationally with \(-S\) nef over \(X\). Polyhedrality then makes \(H-\delta S\) nef and big for small rational \(\delta>0\). For the positive parameters use the MMP for \(K_Y+\Omega-\epsilon S\), with scaling of a diluted representative of \((\epsilon/\delta)H\) on that output. Here \(\Omega-\epsilon S\) is effective klt and big for small \(\epsilon>0\), since \(-K_Y\) is big. The supports are \(H_r=H_i-rS_i\), nef between successive parameters as well as at each wall \cite[Theorems~1.2, 1.5 and 1.6]{BW}. At a wall \(H_r\) descends to the contraction target by relative base-point-freeness. Each chosen ray has positive \(S_i\)-degree, so \(S\) survives and the completed pair stays crepant plt. Termination gives a Mori fibre space \(g:W\to T\) at parameter \(r_*\), with \(S_W\) relatively ample and \(H_*=H_W-r_*S_W\) nef and \(\Q\)-linearly pulled back from \(T\). This is the ordered programme of \cref{prop:global-walls}, using just its construction of the walls.

At a birational wall with \(0<b\leq1\), transport the integral \(L_i=H_r-bS_i\) by \cref{kv:transport} with boundary \(\Omega_i-bS_i\). Indeed,
\[
L_i-(K_{X_i}+\Omega_i-bS_i)=H_r-(K_{X_i}+\Omega_i)
\]
is the same rational pullback up to the common principal relation on the two sides. For \(b=0\), use an effective crepant klt Calabi--Yau boundary obtained by mixing \(\Omega\) with a klt Calabi--Yau boundary on \(Y\) before taking transforms.

Put \(r_j=d+j\) and let \(r_J\) be the first of these values at least \(r_*\). At \(r_j<r_*\), use the model before all steps at \(r_j\). The restriction sequence of \cref{lem:restriction} is
\[
 \begin{gathered}
 0\to\cO(H_{r_j}-S_i)\to\cO(H_{r_j})
          \to\cO_{S_i}(G_j)\to0,\\
 G_j\simq H_{r_j}|_{S_i}-\delta_j,\qquad
 0\le\delta_j\le\Diff_{S_i}(0).
 \end{gathered}
\]
If \(G_j\) has no higher cohomology, the source surjects onto the middle term in degrees one and two. Transport it across each wall \(s\in[r_j,r_{j+1})\) with correction \(r_{j+1}-s\), which equals one at \(s=r_j\). In the last interval include every wall at \(r_*\). Before \(r_0=d\) use the correction \(d-s\), including the terminal wall if \(r_*\le d\). Thus the terminal groups for \(L_*=H_*-(r_J-r_*)S_W\) surject onto the two groups on \(X\). Write \(a=r_J-r_*\in[0,1)\).

For \(a=0\), \cref{kv:small} applies to the nef non-big divisor \(L_*\). For \(a>0\) over a positive-dimensional base, \(L_*\) is relatively anti-ample and \(L_*-K_W\sim_{\Q,T}(1-a)S_W+(\Omega_W-S_W)\) is relatively ample, so \cref{kv:relative} applies. Over a point both \(L_*\) and \(K_W-L_*\) are anti-ample. If \(\Omega_W-S_W\) is big, \cref{thm:rank-one-vanishing} gives vanishing in the two intermediate degrees. For a plt completion of fixed index we can instead use \cref{kv:surface}. For an anti-ample \(J\) and the quotient \(G_j\) along \(S_W\) for \(J-jS_W\), write
\[
G_j\simq (J-jS_W)|_{S_W}-\delta_j,
\qquad 0\leq\delta_j\leq\Diff_{S_W}(0).
\]
For \(K_{S_W}-G_j\) use the boundary \(\delta_j\). The difference is \(-(J-jS_W)|_{S_W}\), which is ample. Surface duality gives \(H^1(S_W,\cO(G_j))=0\). Descending the restriction sequences from a large Cartier multiple of \(S_W\) gives \(H^1(W,\cO(J))=0\), using duality and Serre vanishing for \(j\gg0\). Apply this to both anti-ample divisors and use threefold duality.

\begin{proposition}\label{kv:strict}
Let \(X\) be a projective threefold of Fano type over an algebraically closed field of characteristic \(p\) sufficiently large. If an effective rational \(P\) satisfies \(K_X+P\simq0\) and \((X,P)\) is non-lc, then \(H^i(X,\cO_X(D))=0\) for every ample integral \(\Q\)-Cartier Weil divisor \(D\) and \(i>0\).
\end{proposition}

\begin{proof}
Take a small \(\Q\)-factorialization, retaining the pullback of \(D\), which is nef and big. Denote this model by \(X_0\). Choose a klt log Fano boundary \(C\) with \(A=-(K_{X_0}+C)\) ample. Mix the strict transform of \(P\) and a klt log Calabi--Yau boundary separately with the same small weight \(s>0\) on \(C\). The mixtures are respectively non-lc and klt. For \(\mu>0\) small, add the same diluted effective representative \(J\simq sA-\mu D\) to both, preserving the klt endpoint. For \(v\geq0\) small, subtract \(vJ\) from both and take the boundary \(B_v\) at the first non-klt parameter on their segment from the klt endpoint. Then
\[
-(K_{X_0}+B_v)\simq \mu D+vJ .
\]
On one log resolution the lc inequalities are finite linear inequalities in \(v\) and the segment parameter. Thus \(B_v\) converges to \(B_0\), and every coefficient-one component on that resolution for small \(v>0\) has coefficient one for \(B_0\). The description of lc places of a log smooth subpair then gives the same inclusion for all lc places.

Use the extraction in the proof of \cref{prop:global-plt-preparation} for \(B_v\) with \(v>0\). On a crepant dlt modification run the canonical MMP with the reduced exceptional boundary over \(X_0\). The positive bare discrepancies and negativity contract all exceptional primes. Its last divisorial contraction extracts a single prime \(S\) with \((Y,S)\) plt and \(-S\) relatively ample. If there are no exceptional primes, take a component of \(\floor{B_v}\) on \(X_0\). Interpolating the crepant boundary for \(B_v\) slightly toward \(S\) gives a plt log Fano boundary \(\Gamma\) with floor \(S\), since \(\mu D+vJ\) is ample on \(X_0\) and \(-(K_Y+S)\) is relatively ample in the extraction case. The crepant transform of \(B_0\) has coefficient one on \(S\) and negative log divisor \(\Q\)-linearly equivalent to \(\mu H\).

Complete \(\Gamma\) to a plt log Calabi--Yau boundary by a diluted representative of \(-(K_Y+\Gamma)\). Apply the ordered programme to \(Y\to X_0\to X\) with the surface control in \cref{prop:global-walls}, scaling at zero by a diluted multiple of \(-(K_Y+\Gamma)\). That proof only uses the rational supporting divisor \(H\) up to positive rescaling. It supplies a klt log Fano boundary \(\Xi_i\geq\Diff_{S_i}(0)\) on each \(S_i\), and a crepant plt completion with big residual. At \(r=d+j<r_*\), \cref{lem:restriction} writes the quotient as \(G\simq H_r|_{S_i}-\delta'\), \(0\leq\delta'\leq\Diff_{S_i}(0)\). The difference \(G-(K_{S_i}+\Xi_i-\delta')\) is ample, so the quotient vanishes in positive degrees by \cite[Theorem~1.1]{ABL}. The terminal calculation proves the intermediate vanishing. For the top degree, \(K_X-D\) has negative degree against the square of an ample divisor, by comparison with a klt log Fano witness. Duality gives the result.
\end{proof}

\begin{proposition}\label{kv:single}
Fix a positive integer \(n\). Let \(X\) be a projective threefold of Fano type over an algebraically closed field of sufficiently large characteristic depending on \(n\). Suppose an effective lc boundary \(\Sigma\) satisfies \(n(K_X+\Sigma)\sim0\). If \((X,\Sigma)\) has exactly one lc place, then \(H^i(X,\cO_X(D))=0\) for every ample integral \(\Q\)-Cartier Weil divisor \(D\) and \(i>0\).
\end{proposition}

\begin{proof}
Take a crepant \(\Q\)-factorial dlt modification \(h:Y\to X\). Its completed boundary \(\Omega=S+B\) is plt with the single prime floor \(S\). Interpolating with the crepant transform of a klt log Fano witness gives an effective klt boundary with nef and big negative log divisor on \(Y\). A Kodaira perturbation gives Fano type. Adjunction to \(S\) and to its subsequent transforms gives klt log Calabi--Yau boundaries \(\Theta_i=\Diff_{S_i}(B_i)\) with coefficients in \(\{j/n\}\).

If \(d=0\) and \(S\) is exceptional, the first restriction \(H|_S\) comes from \(h(S)\) and need not be big. Use the fractional part \(B\) of the complement and first run its MMP over \(X\), making \(B\) relatively nef and semiample. Adding a small positive multiple of \(B\) to a klt Calabi--Yau boundary realizes this as a klt adjoint MMP. If \(C_X\) is a klt log Fano witness on \(X\) and \(a=a(S,X,C_X)>0\), then \(B\sim_{\Q,X}C_X^s-aS\), where \(C_X^s\) is the strict transform. Thus a divisorial step contracting \(S\) cannot be \(B\)-negative, by testing curves covering \(S\) in its fibres. There is no other exceptional prime over \(X\). The programme remains crepant for \(\Omega\) and the integral \(H\) is transported over \(X\) as at a wall with \(b=0\).

On this model \(-S|_S\) is relatively big. Over the local ring at the generic point of \(h(S)\), lift a transcendence basis of the residue field of \(S\) over \(k(h(S))\) to units of its divisorial valuation ring. Write the lifts as ratios of regular functions on the local base, with numerator and denominator of equal \(S\)-order. Clearing denominators, multiplying all numerators and the common denominator by a function of positive order if needed, and taking a common power gives sections of some Cartier \(-mS\) of equal order \(m\) along \(S\). Their ratios on \(S\) give that basis up to taking powers, proving the assertion.

Now \(B|_S\) is relatively big by \(B\sim_{\Q,X}C_X^s-aS\). Since \(B\) is relatively semiample, \((H+tB)|_S\) is nef and big for \(t>0\) sufficiently small, by the induced contraction over \(X\) and ampleness of \(D\). For the first quotient use the test boundary \(\Theta-\delta'-t B|_S\) in \cref{kv:surface}. If \(S\) is not exceptional, \(H|_S\) is already nef and big. After this initial quotient when \(d=0\), use the parameter-zero steps for \(-S\) with the offset \(H-S\). If \(d>0\), start with \(H-dS\) as in the ordered programme.

For \(0<r<r_*\), the supporting divisor \(H_r\) is big. Indeed \(H_i\) is big. To see that \(H_i-r_*S_i\) is pseudo-effective, observe that at every remaining wall \(s\le r_*\) its degree on the contracted ray is \((s-r_*)S_i\cdot R\le0\). Negativity therefore recovers this class from the nef terminal support by adding effective exceptional divisors on common resolutions. The class \(H_r\) lies on the open segment joining these two classes. Immediately before the steps at a nonterminal \(r>0\), nefness of \(H_r+\eta S_i\) gives \(S_i\cdot C\geq0\) on every \(H_r\)-trivial curve. The semiample birational contraction of \(H_r\) therefore cannot contract \(S_i\) by negativity. Consequently \(H_r|_{S_i}\) is nef and big. For the quotient \(G\simq H_r|_{S_i}-\delta'\), the boundary \(\Theta_i-\delta'\) gives \(G-(K_{S_i}+\Theta_i-\delta')\simq H_r|_{S_i}\). The quotient therefore has zero higher cohomology by \cref{kv:surface}.

The last parameter now gives the desired intermediate vanishing. The top degree follows by duality as in \cref{kv:strict}.
\end{proof}

\subsection{Frobenius splitting and the strata of a complement}\label{s:kv-strata}
\begin{lemma}\label{kv:positive}
Let $X$ be a normal projective variety over a perfect field of positive characteristic. If $X$ is quasi-$F$-split and $L$ is an ample integral $\mathbb Q$-Cartier Weil divisor, then $H^i(X,\mathcal O_X(L))=0$ for $i>0$.
\end{lemma}

\begin{proof}
Write $Q_{L,j}$ and $B_{L,j}$ for the sheaves in \cite[Lemma~3.8]{KTTWYY}. For $j$ at least the splitting height $h$, $\mathcal O_X(L)$ and $B_{L,j}$ are direct summands of $Q_{L,j}$ by \cite[Lemma~3.11]{KTTWYY}. Take $m>0$ with $mL$ Cartier. For fixed $j=h$, the sheaves $Q_{p^rL,h}$ are obtained from the finitely many $Q_{aL,h}$, $0\leq a<m$, by tensoring with growing multiples of $mL$. They have no higher cohomology for $r$ sufficiently large by Serre vanishing. Descend from $i=r$ to $i=0$ using
\[
 0\longrightarrow F_*B_{p^{i+1}L,h+r-i-1}
 \longrightarrow Q_{p^iL,h+r-i}
 \longrightarrow F_*\mathcal O_X(p^{i+1}L)\longrightarrow0 .
\]
Both end terms are Frobenius pushforwards of summands of the sheaf at the previous step. This proves the vanishing for $Q_{L,h+r}$ and hence for its summand $\mathcal O_X(L)$.
\end{proof}

A quasi-$F$-splitting descends through a proper birational contraction between normal varieties by pushing forward its Witt retraction. Indeed, Witt functions are tuples of ordinary functions, and the direct image of the structure sheaf is the structure sheaf. A retraction also extends across a subset of codimension at least two by normality, applied to each Witt coordinate \cite[Remark~3.12]{KTTWYY}. Thus quasi-$F$-splitting persists through birational steps of the adjoint MMP and on the ample model.

We now consider a projective threefold $X$ of Fano type over an algebraically closed field $k$ for which $(X,P)$ is lc for every effective rational divisor $P$ with $K_X+P\sim_{\mathbb Q}0$. Fix one such boundary $\Sigma$ for which $(X,\Sigma)$ is non-klt and $n(K_X+\Sigma)\sim0$. We take $p$ sufficiently large in terms of $n$, in particular $p\nmid n$. A global F-splitting suffices for ample Weil vanishing by \cref{kv:positive}.

\begin{lemma}\label{kv:residual}
On a crepant $\mathbb Q$-factorial dlt modification $(V,F+B)\to(X,\Sigma)$, write $F$ for the reduced floor. There is a birational $B$-MMP to a Fano-type model $(W,F_W+B_W)$ on which $B_W$ is semiample, preserving the generic points of all lc strata. Put $b=\kappa(V,B)$. Then $0\le b\le2$. If $b=0$, then $B_W=0$. If $b>0$, write $g:W\to T$ for the contraction of $B_W$. There is an effective ample rational Cartier divisor $A$ on $T$ such that $B_W=g^*A$. Every lc stratum of $(W,F_W)$ dominates $T$.
\end{lemma}
\begin{proof}
The dlt modification is of Fano type by the crepant-boundary perturbation in the proof of \cref{kv:small}. Its birational contractions remain of Fano type by pushing a big klt Calabi--Yau completion.

For small rational $\eta>0$, the pair $(V,F+(1+\eta)B)$ is dlt, as is seen on a log resolution witnessing dlt for $(V,F+B)$. Its log divisor is rationally equivalent to $\eta B$. Run its MMP to a nef model, since $B\ge0$. Equivalently, one can mix with a big klt log Calabi--Yau boundary to run the same MMP as a klt adjoint MMP \cite[Theorem~1.6]{BW}. The resulting pair with the increased boundary remains dlt, hence so does the pair with the original boundary. All contracted curves lie in the support of $B$, which avoids the generic point of each lc stratum. The contractions and flips are isomorphisms at these generic points. Discrepancies of the increased pair cannot decrease, so a new lc stratum would come from an old lc place already in this isomorphism locus. The original $n$-principal log divisor transfers crepantly at each step. If $B_W\equiv0$, effectivity gives $B_W=0$. Otherwise the nef endpoint is semiample by base-point-freeness applied to a Cartier multiple on a Fano-type model \cite[Corollary~1.5]{ABL}.

On a common resolution with maps $u$ to $V$ and $v$ to $W$, negativity gives $u^*B=v^*B_W+E$ for $E\ge0$ exceptional over $W$. This also preserves the spaces of sections at divisible multiples and hence the Iitaka dimension. For every $B'_W\ge0$ rationally equivalent to $B_W$, replace $B$ by $u_*(v^*B'_W+E)$. The two replacements change the log divisors by the same rational principal divisor. The resulting pair on $V$ pushes to an effective anticanonical completion on $X$, crepant to the pair on $W$. Thus $(W,F_W+B'_W)$ is lc. In particular $B_W$ is not big, since otherwise an effective replacement could include a positive multiple of a component of $F_W$.

If $b>0$, write $B_W\sim_{\mathbb Q}g^*A_0$ for an ample rational Cartier divisor $A_0$ on $T$. The effective divisor $B_W$ has no horizontal prime, since its restriction to the proper normal generic fibre is rationally trivial. A rational function giving the equivalence has zero divisor on that fibre and is therefore constant on it. The generic fibre has constant field $k(T)$ because $g_*\mathcal O_W=\mathcal O_T$. Absorb the function in $A_0$ to get equality $B_W=g^*A$. Its effectivity implies that $A$ is effective by considering primes above each prime of $T$. If an lc centre of $(W,F_W)$ had proper image in $T$, take an effective rational representative of $A$ containing a positive multiple of a very ample divisor through that image. Adding its pullback makes the discrepancy of an lc place over the centre negative, contrary to the replacement argument.
\end{proof}

If the fractional part of the boundary has Iitaka dimension zero and the smallest strata are curves, the Witt cohomology of the double curves gives a quasi-Frobenius splitting.

\begin{proposition}\label{kv:witt}
In sufficiently large characteristic let $(X,\Sigma)$ be a projective lc threefold pair with $X$ of Fano type and $n(K_X+\Sigma)\sim0$ for a positive integer $n$ prime to $p$. Suppose that on a crepant $\mathbb Q$-factorial dlt modification $(V,F+B)$ the floor $F$ has a curve stratum and has no point stratum, and $\kappa(B)=0$. Then $X$ is quasi-$F$-split.
\end{proposition}

\begin{proof}
Run the $B$-MMP of \cref{kv:residual}, using only the construction of its semiample model. The trace of $B$ has Iitaka dimension zero, so it vanishes. This MMP is crepant for the log divisor. Its steps are isomorphisms at the generic points of the lc strata, since every $B$-negative curve lies in $\Supp B$ and $B$ contains no stratum. Take the tame cyclic index cover $W'\to W$ of $K_W+F_W$ and let $G$ be the reduced inverse image of $F_W$. The cover is quasi-\'etale and is \'etale on the snc locus by purity. It preserves the dlt property and the dimensions of the minimal strata, by the tame discrepancy formula. We have $K_{W'}+G\sim0$, and $W'$ is of Fano type by pulling back a klt log Fano witness.

We first compute the Witt cohomology of the boundary. By \cref{kv:small}, $H^{>0}(W',\mathcal O_{W'})=0$. The sequence
\[
 0\longrightarrow\omega_{W'}\longrightarrow\mathcal O_{W'}
       \longrightarrow\mathcal O_G\longrightarrow0
\]
shows that $G$ is Cohen--Macaulay by the depth lemma. Duality then shows that $G$ is connected, $H^1(G,\mathcal O_G)=0$, and $H^2(G,\mathcal O_G)=k$. The components of $G$ are normal and they meet in smooth curves. There are no triple strata. On each double curve $E$, adjunction gives $K_E+\operatorname{Diff}_E\sim0$. The different is integral, and any coefficient one would give a point stratum. Hence $E$ is elliptic.

On a component $G_i$ the sum of the double curves is a reduced anticanonical divisor. On its minimal resolution $T_i$ the strict transforms of these curves are still disjoint and anticanonical. Indeed, the adjoint is Cartier and the induced pair on $G_i$ is plt. The exceptional coefficients are nonnegative by the intersection matrix on the minimal resolution and less than one by plt, so they vanish by integrality. For a rational $T_i$, the anticanonical sequence and $H^1(T_i,\omega_{T_i})=0$ give just one double curve. Otherwise $T_i$ is ruled over a curve of positive genus. The elliptic double curves are horizontal and their degrees sum to two on the general ruling fibre. The base is elliptic by Riemann--Hurwitz. In this case there are two sections. A single bisection would give a nonzero map on $H^1(\mathcal O)$ by trace of degree two, whereas the anticanonical sequence and duality make that restriction map zero.

Thus the dual graph is a path or cycle, with a single curve for each edge and with rational end components in the path case. The normalization sequence
\[
 0\longrightarrow\mathcal O_G\longrightarrow
 \bigoplus_i\mathcal O_{G_i}\longrightarrow
 \bigoplus_E\mathcal O_E\longrightarrow0
\]
is exact. The difference map is surjective locally, since only two components meet. Its kernel is $S_2$ and agrees with $\mathcal O_G$ in codimension one. The sequence on constants and $H^1(G,\mathcal O_G)=0$ exclude a cycle. The same sequence holds for the additive sheaves $W_j\mathcal O$, since Witt vectors preserve fibre products and restriction onto an intersection is surjective. Rational end components have zero higher Witt cohomology in finite length. On an inner component, each of its two sections induces an isomorphism on $H^1(W_j\mathcal O)$, and $H^2(W_j\mathcal O)$ vanishes. These assertions follow first for $\mathcal O$ by resolution and the elliptic ruling, then for $W_j\mathcal O$ by truncation. Taking the cokernel of the component-to-curve difference map along the path identifies $H^2(G,W_j\mathcal O_G)$ with $H^1(E,W_j\mathcal O_E)$ for a double curve $E$, compatibly with Frobenius and truncation. Frobenius is injective in the inverse limit, also when $E$ is supersingular, by \cite[Lemma~6.9]{STYQuasi} with zero twist.

To transfer the calculation to \(X\), take a smooth common resolution \(Y\) of \(X,V,W\). Resolve the crepant boundary \(\Sigma_Y\), the exceptional locus and the divisor of a function \(u\) trivializing the adjoint. Write
\[
 K_Y+\Sigma_Y=D_Y=\divi_Y(u)/r,\qquad F_Y=\Sigma_Y^{=1},
\]
where \(r\mid n\) is the least principal index. Normalize \(Y\) in \(z^r=u\) to get \(a:U\to Y\). This extension also defines \(W'\), so \(U\) maps to \(W'\). Locally the cover is a tame root of an SNC monomial after removing a unit by an \'etale extension. Take a smooth toroidal resolution \(b:Z\to U\). We use \(Z\to W'\) to describe the boundary on \(Z\), and \(Z\to U\to Y\) to recover the twist by \(D_Y\).

Let \(F_U=(a^{-1}F_Y)_{\rm red}\) and let \(F_Z\) be the coefficient-one part of the crepant boundary. By the toroidal discrepancy formula, a ray belongs to \(F_Z\) precisely when all its positive coordinates correspond to components of \(F_U\).

We claim that $Rb_*\mathcal O_Z(-F_Z)=\mathcal O_U(-F_U)$. For the zeroth image, vanishing along the original floor forces vanishing along each new floor prime, whose centre is in that floor. To prove the higher vanishing work in a toroidal chart. Write $P_U,P_Z$ for the full toroidal boundaries and $C=P_U-F_U$. The cones on $U$ are simplicial, so $C$ is $\mathbb Q$-Cartier. Every ray off the floor has positive order along $b^*C$, and the floor rays have order zero. For sufficiently small $t>0$, the divisor $P_Z-F_Z-tb^*C$ therefore has coefficients in $[0,1)$, pushes to $(1-t)C$, and its adjoint is relatively equivalent to $-F_Z$. Apply Grauert--Riemenschneider vanishing over the klt pair $(U,(1-t)C)$ \cite[Theorem~3]{BK}. The same direct-image identity for $W_j\mathcal I_{F_Z}$ follows from the truncation sequences.

The morphism $Z\to W'$ resolves a klt threefold, so $H^{>0}(Z,\mathcal O_Z)=0$ by rationality and \cref{kv:small}. Every additional floor divisor on $Z$ dominates a double curve on $W'$. At the generic point of that curve the induced birational morphism of regular surfaces factors into point blowups. Blowing up a crossing of two floor components subdivides the edge and creates a projective line over the residue field of the double curve. No other blowup creates a floor prime. Thus the inserted surfaces are birationally ruled over the same elliptic curve, with double curves that map isomorphically to it. There are no extra floor intersections over points, since an snc floor stratum there would give a zero-dimensional lc centre on $W'$. The strict transforms of the old components preserve the rational or ruled cohomology calculation. Accordingly $F_Z$ has the same path calculation as $G$.

Since $W_j\mathcal O_Z$ has no higher cohomology, the Witt exact sequence for the ideal of $F_Z$ gives
\[
 H^3(Z,W_j\mathcal I_{F_Z})\simeq H^2(F_Z,W_j\mathcal O_{F_Z}).
\]
These isomorphisms commute with Frobenius and truncation. The finite-length cohomology modules satisfy the Mittag--Leffler condition \cite[Proposition~2.32]{STYQuasi}. In the inverse limit the calculation on the double curve therefore gives injectivity of Frobenius on $H^3(Z,W\mathcal I_{F_Z})$.

Multiplication by \([z]\) identifies \(W_j\mathcal I_{F_Y}(D_Y)\) with the summand of \(a_*W_j\mathcal I_{F_U}\) on which a deck transformation \(\sigma\) acts by \([\sigma(z)/z]\). Indeed, for a prime \(P'\) over \(P\) the order of the \(i\)-th coordinate is
\[
 \ord_{P'}(z^{p^i}f_i)=e(P'/P)
       \bigl(\ord_P f_i+p^i\coeff_P D_Y\bigr).
\]
At a floor prime \(D_Y\) is integral and \(e(P'/P)=1\), so the coordinate must vanish there. Elsewhere its required order is nonnegative. These conditions give exactly the rounded coordinate bounds of Subsection~\ref{pre:kv}. Normality extends the identification globally and tame averaging projects onto the summand. For \(pD_Y\) use \([z^p]\). Frobenius takes the first character to its $p$-th power. Together with the toroidal direct-image comparison this proves injectivity of
\[
 H^3(Y,W\mathcal I_{F_Y}(D_Y))
 \xrightarrow{F} H^3(Y,W\mathcal I_{F_Y}(pD_Y)).
\]
The finite-length Frobenius kernels form a Mittag--Leffler system with zero inverse limit. At some length their image in length one is therefore zero. The pushout calculation in the proof of \cite[Proposition~3.9]{STYQuasi} uses only the defining exact row and top cohomology, so it gives the required injectivity on $Y$ also for the crepant subboundary. Apply \cite[Proposition~6.1(1)]{STYQuasi} to obtain a pure quasi-$F$-splitting of $(X,\Sigma)$. The underlying $X$ is quasi-$F$-split by \cite[Remark~3.5]{STYQuasi} and restriction of the retraction.
\end{proof}

For a positive-dimensional base, the following lemma descends a splitting from a stratum that is generically finite over the base.

\begin{lemma}\label{kv:multisection}
With the notation of \cref{kv:residual}, suppose $b>0$ and $M$ is a normal dlt lc stratum of dimension $b$ dominating $T$. Factor $M\to T$ as $\mu:M\to U$ followed by a finite map $\gamma:U\to T$. Write $\Theta_M^0$ for the iterated different on $M$ from the reduced pair $(W,F_W)$ along a flag of lc strata, and put $\Theta_U^0=\mu_*\Theta_M^0$. If $\deg\gamma$ is prime to $p$ and $(U,\Theta_U^0)$ is globally sharply F-split, then $X$ is globally F-split.
\end{lemma}
\begin{proof}
Use $q=p^e$ with $n\mid(q-1)$, choosing $e$ a multiple of an exponent of a splitting along $\Theta_U^0$. A splitting along a boundary iterates at all multiples of its exponent. The principal multiple of the log divisor defines a nonzero rational $p^{-e}$-linear map $\varphi$ with attached divisor $F_W+B_W$ on $W$. In particular it is regular on $F_*^e\mathcal O_W$. On open subsets of $T$, applying $\varphi$ to pullbacks of regular functions and using $g_*\mathcal O_W=\mathcal O_T$ defines a map $\varphi_T$ on $T$. At each coefficient-one divisor of the flag the map preserves the ideal of that divisor and induces a nonzero rational map on it. The boundary of the induced map is the ordinary different \cite[Theorem~5.3]{Das}. Consequently the induced map $\varphi_U$ on $k(U)=k(M)$ has boundary $\Theta_U^0+\gamma^*A$, by taking the trace of the divisor on the birational Stein factor. On elements of $k(T)$ this map agrees with $\varphi_T$, by restriction on the flag.

For a finite separable extension $L/K$, the map $L\otimes_K F_*^e K\to F_*^e L$ sending $c\otimes a$ to $F_*^e(c^q a)$ is an isomorphism. Thus $\varphi_T$ is nonzero and $\varphi_U$ is its unique extension. A splitting along $\Theta_U^0$ at exponent $e$ has the form $\varphi_U(f'\,\cdot)$ for $f'\in k(U)$ and gives
\[
 \varphi_U(f')=1,\qquad \operatorname{div}(f')+(q-1)\gamma^*A\ge0 .
\]
Put $f=\operatorname{Tr}_{k(U)/k(T)}(f')/\deg\gamma$. The trace commutes with the maps, since $\operatorname{Tr}(c^q)=(\operatorname{Tr}c)^q$. Hence $\varphi_T(f)=1$. The valuation inequality descends by summing the conjugates in a normal closure. It gives $\operatorname{div}(f)+(q-1)A\ge0$. This effective rational Cartier divisor pulls back effectively on $W$ at every prime, including primes above points. The negativity comparison in \cref{kv:residual} gives $\operatorname{div}_V(f)+(q-1)B\ge0$ as well. The rational map $\varphi$ is the same on $V$ by crepancy. Therefore $\varphi(f\,\cdot)$ is a regular Frobenius retraction on $V$ and pushes to $X$.
\end{proof}

First, if there is a zero-dimensional stratum on $V$, apply the same divisorial adjunction successively down a flag to it. A nonzero $p^{-e}$-linear map on the perfect field $k$ is nonzero on $1$. Accordingly, the global map on $V$ has nonzero value on $1$ and can be rescaled to a splitting. Push to $X$.

Suppose next that $b=2$. On the proper normal generic curve of $g$, adjunction from the generic point of the base gives $K+F_W|_{\eta}\sim_{\mathbb Q}0$. Its constant field is $k(T)$ and the restriction of $F_W$ has positive degree because all its components dominate $T$. The genus formula gives genus zero and total boundary degree two. Any floor surface $S$ therefore has generic degree at most two and is separable over $T$. By \cref{kv:residual}, distinct floor components are disjoint, so $(W,F_W+B_W)$ is plt near $S$. The klt adjunction boundary on $S$ is $\Theta_S^0+(g|_S)^*A$, where $\Theta_S^0$ has standard coefficients.

On the Stein factor $U$ we have $n(K_U+\Theta_U^0+\gamma^*A)\sim0$, and adjunction and pushforward are crepant for this full log divisor. The full boundary is klt, with log discrepancies at least $1/n$. Removing $\gamma^*A$ preserves this lower bound. Since $-(K_U+\Theta_U^0)\sim_{\mathbb Q}\gamma^*A$ is ample and the standard coefficients bounded above by $1-1/n$ form a finite set, \cite[Proposition~3.2 and Remark~3.3]{CTW} gives global F-regularity for sufficiently large $p$ in terms of $n$. Apply \cref{kv:multisection} with $M=S$.

Suppose now that $b=1$ and write $C=T$, $K=k(C)$, and $(P,G)=(W_\eta,F_W|_{\eta})$. The generic surface $P$ is of del Pezzo type over $K$. Its reduced dlt pair satisfies $n(K_P+G)\sim0$. The connected tame cyclic cover of the log torsion class is quasi-\'etale, since $K_P+G$ is integral. The cover is also of del Pezzo type. Over its own finite field of constants, \cite[Theorem~1.7]{BT} gives vanishing of higher cohomology of its structure sheaf. The torsion eigensheaf therefore gives $H^1(P,\mathcal O_P(K_P+G))=0$. Duality on the normal surface then gives $H^1(P,\mathcal O_P(-G))=0$. As $G$ is nonempty, the restriction sequence gives
\[
 H^0(G,\mathcal O_G)=K.
\]

If there is a curve stratum, the generic reduced boundary $G$ has a node. Over $K^{\mathrm{sep}}$ its components are normal and meet at SNC nodes. Indeed, normality and the dlt resolution criterion are preserved under separable field extensions. There are no self-edges between components of the base change of the same normal component, since those components are disjoint. The graph is connected by base change of the displayed equality. Let $K_i$ be the constants of one component and let $L$ be the residue field of an incident node. Both are finite purely inseparable over $K^{\mathrm{sep}}$. Adjunction from the generic log Calabi--Yau pair gives coefficient one at each node. The canonical class calculated over $K_i$ agrees with that calculated over $K^{\mathrm{sep}}$, since the dualizing module of the finite extension is a one-dimensional $K_i$-space. Taking degrees gives $[L:K_i]\le2$ and valence at most two. Thus $L=K_i$ for $p>2$. The constants agree across every node. They glue to global functions on the nodal curve and hence are all $K^{\mathrm{sep}}$. In particular each original node is separable over $K$.

The graph is a path or a cycle, allowing a cycle with two parallel edges. Its Galois action has finite image $J$ in a dihedral group. A stabilizer of an edge cannot interchange its vertices, since they come from distinct original components. Fixing an edge and its two vertices fixes the whole graph. Consequently the closure $C'$ of a chosen generic node is a normal lc stratum and $\gamma:C'\to C$ is a Galois cover with group $J$. Adjunction twice gives a klt different $\Theta_{C'}^0+\gamma^*A$ with $n$-principal log divisor, where $\Theta_{C'}^0$ is the different for $(W,F_W)$. There is no zero-dimensional stratum on $W$. The bare different $\Theta_{C'}^0$ has standard coefficients. The different of a reduced lc boundary is standard, and iteration preserves the standard set in the lc range. In the usual formula $(m-1+\sum l_i b_i)/m$, a sum of standard nonzero $b_i$ at most one is either a single term or two half-coefficients. Since $-(K_{C'}+\Theta_{C'}^0)\sim_{\mathbb Q}\gamma^*A$ is ample, $C'\simeq\mathbb P^1$, and $C\simeq\mathbb P^1$ as well by L\"uroth's theorem.

The construction of the maps in the proof of \cref{kv:multisection} uses only $n\mid(p^e-1)$ and separability up to the multiplier step. Uniqueness of extension from $k(C)$ shows that $\varphi_{C'}$, and hence $\Theta_{C'}^0$, is $J$-invariant. If $p\mid |J|$, a nontrivial normal cyclic $p$-Sylow of the dihedral subgroup has a unique fixed point $x$ on $\mathbb P^1$, fixed by all of $J$. Outside $x$ its orbits have size at least $p$. As the different has degree less than two and nonzero coefficients at least $1/2$, it would be supported at $x$. Replace $A$ by $(\deg A)\gamma(x)$. The fibre of $\gamma$ there is $|J|x$, so the iterated different for the replacement has coefficient
\[
 \operatorname{coeff}_x\Theta_{C'}^0+|J|\deg A=2
\]
at $x$. This contradicts the log canonicity of all replacements from \cref{kv:residual}. Thus $p\nmid |J|$. Global F-regularity of the standard klt log Fano curve follows from \cite[Proposition~5.5]{SatoTakagi}, and \cref{kv:multisection} gives a splitting of $X$.

It remains to consider $b\le1$ when all lc strata have dimension two, as the case $b=0$ with a curve stratum is treated in \cref{kv:witt}. For $b=1$ the generic boundary is connected by $H^0(\mathcal O_G)=K$, so the disjoint floor surfaces consist of one prime. For $b=0$, use \cref{kv:small} on the numerically trivial integral Weil divisor $K_W+F_W$. Duality gives $H^1(W,\mathcal O_W(-F_W))=0$, and the restriction sequence gives $H^0(F_W,\mathcal O_{F_W})=k$. Again there is just one prime. All floor primes survived the $B$-MMP, so on the original dlt modification the sole floor is a prime $S$ and the pair is plt. Thus $S$ is the unique lc place, and the ample Weil vanishing in this case follows from \cref{kv:single}.

\begin{proof}[Proof of \cref{th:kv}]
Run the adjoint MMP of Subsection~\ref{s:kv-adjoint}. At a negative Mori contraction over a curve or a surface, apply \cref{kv:relative}. The non-big and negative rank-one outputs are covered by \cref{kv:small}. For the remaining ample divisor take the fixed-index zero-boundary complement constructed in Subsection~\ref{s:kv-adjoint}. If it is klt, the same completed pair on the anticanonical ample model is klt and has index dividing \(n\). All its log discrepancies are positive multiples of \(1/n\). Apply \cref{kv:bounded}. A non-lc effective anticanonical representative gives \cref{kv:strict}. Otherwise every effective anticanonical representative is lc, and the preceding analysis of the strata applies to the fixed non-klt completion. The characteristic bounds used here depend only on \(n\), the finite coefficient sets it determines and the fixed discrepancy \(1/(2n)\). The zero-boundary complement theorem fixes \(n\) before \(X,D,\Delta\) are chosen. Taking the maximum proves the asserted uniformity.
\end{proof}

\end{document}